%% file: bott2025-2.tex
\input{RVpreamble.tex}

\newcommand{\SmArt}{\text{(SmArt)}}
\newcommand{\pp}{\mathbb{P}}
\newcommand{\et}{\text{\'et}}

\newcommand{\Bun}{\operatorname{Bun}}

\renewcommand{\Maps}{\operatorname{Maps}}

\newcommand{\nicely}{\text{(nicely-$k$-conn)}}
\newcommand{\connected}{\text{($k$-conn)}}
\newcommand{\cGsp}{\cG_{\text{sp}}}
\newcommand{\cGw}{Ho\mathscr{G}}
\newcommand{\cGst}{\cG}
\newcommand{\alg}{\text{alg}}
\newcommand{\Top}{(\text{Top})}
\renewcommand{\top}{\operatorname{top}}
\newcommand{\HoTop}{(\text{HoTop})}
\renewcommand{\tU}{T}
\renewcommand{\cX}{\Bun^\bullet}
\renewcommand{\cY}{\Bun}

\input{BCpreamble}

\begin{document}
\pagestyle{plain}
\title{\Large{Complex Bott Periodicity in algebraic geometry}}
\author{Hannah Larson and Ravi Vakil \\ appendix by Benjamin Church}
\address{Dept.\  of Mathematics 970 Evans Hall, MC 3840 Berkeley, CA 94720-3840; Dept.\ of Mathematics, Stanford University, Stanford CA~94305--2125}
\email{hlarson@berkeley.edu; rvakil@stanford.edu; bvchurch@stanford.edu}
\begin{abstract}We state and prove a form of Bott periodicity (for
  $U(n)$) in an algebraic setting (so, $GL(n)$) which makes sense over
  $\Z$, which also specializes to  Bott periodicity in
  the usual sense (hence giving yet another proof of classical Bott periodicity).  An appendix by B. Church gives a specialization of the constructions and results to motivic homotopy theory, which may be of independent interest.
 \end{abstract}
\maketitle
\tableofcontents

{

  \begin{figure}[ht]
  $$
  \xymatrix@R=1.5em{
\boxed{
  {\shortstack{composition of affine bundles \\  minus high codimension subsets}}} \ar[d] \ar@{<->}[r]^-\sim & \boxed{
  {\shortstack{affine bundle, minus \\ high codimension subset}}}  \ar[d] \\
\cX_{r,d} \ar@{<~>}[r]^{\text{$\therefore$ approximate isomorphism}}_{\text{in an appropriate homotopy category}} & BGL(d)  
  }
  $$
  \caption{(Complex) Bott periodicity is at heart an approximate isomorphism of moduli spaces}\label{f:1}
  \end{figure}

  \section{Introduction}

  Bott periodicity is a central result of topology which has been
  foundational for a number of later developments, including K-theory.
  In this paper, we prove a simple (in many senses)  algebro-geometric version of Bott
  periodicity (for $U(n)$), which extends classical Bott periodicity.
  This argument is (in multiple ways) inspired by and based on  work of Segal, Kirwan, Cohen, Mitchell,
  Lupercio, and other topologists (see in particular \cite[\S 8.8]{ps}, as well as
  \cite{cohensegal, cls, lupercio, mitchell, mitchell2, segal}, and more).

  Our philosophy is this:  Bott periodicity is at its root an approximate isomorphism of two spaces of a certain sort --- the moduli space $\cX_{r,d}$ of degree $d$ rank $r$ vector bundles on $\proj^1$ (framed at infinity), and the moduli space $BGL(d)$ of rank $d$ vector spaces --- which ``plays well'' with increasing $r$ and $d$, and becomes closer and closer to an isomorphism as $r$ and $d$ become larger (see Figure~\ref{f:1}).  The top row of Figure~\ref{f:1} is 
  an open subset of the space of representations of the quiver  of Figure~\ref{f:2}.

  \begin{figure}[h]
    \centering
  \[
  \begin{tikzcd}[
      row sep=huge, column sep=huge,
      cells={shape=circle,draw,inner sep=2pt},
      arrows={-Stealth}
    ]
    A
      \ar[loop above,"\alpha"]
      \ar[r,"j"]
    & U
  \end{tikzcd}
  \]
  \caption{The affine bundle over $BGL(d)$ in Figure~\ref{f:1} is the space of representations of this   quiver, where $A$ is a rank $d$ vector space (parametrized by the $BGL(d)$) and $U$ is a rank $r$ vector space with basis.  See Remark~\ref{r:qi} and \cite{bryanvakil} for more.}\label{f:2}
  \end{figure}
  
  Traditionally, Bott periodicity is presented as  related to complex geometry, but the argument here suggests
  that it is a more fundamental statement, which works over the integers and is purely algebraic, from which the
  complex manifestation is obtained by specializing to $\C$.    This is of course slightly
  misleading. To get from the algebraic statement to the topological statement, some topological or holomorphic input is needed
  in order to even make contact with the language of the traditional statement; but this connection just involves standard tools and methods (see \S \ref{pf:trad}).

  One goal, to be hopefully used in future work (such as \cite{bryanvakil}), is to create a simple but versatile setting in which
  homological ``limiting'' arguments can be made in algebraic geometry ---  a setting differently flexible than stacks
  where many of the standard tools (cohomology, Chow rings, Grothendieck rings) can apply.  The specialization to motivic homotopy theory (in Appendix~\ref{appendix} by B. Church) is less trivial; see \S \ref{appendixintroduction} for a friendly introduction to these ideas.

For algebraic geometers: the idea behind Bott periodicity is that
there is a natural homotopy equivalence
$\be:  BU \overset \sim \longrightarrow \Omega^2_d U$, which
(reversing the direction, and using the homotopy equivalence
$U(n)\overset \sim \longrightarrow GL(n)$) translates to
$\al: \Omega^2_d (BGL) \overset \sim \longrightarrow BGL$ for any $d \in \Z$.   (The subscript $d$ refers to the component of ``degree $d$".) Now
$\Omega^2X$ is the space of (pointed) maps of
the  sphere to $X$,  suggesting that for each vector bundle $\cE$ on
$\proj^1 \times S$ over a base $S$, of indeterminate rank, there is
associated (by $\al$) a vector
bundle on the base $S$.  Clearly the first thing to try is 
pushing forward the bundle $\cE$ by $\pi: \proj^1 \times S \rightarrow S$.
But $\cE$ needs to be sufficiently positive
on the fibers of $\pi$ to make this work; so we imagine either (i) twisting the
bundle to make it positive, and {\em then} pushing forward, or (ii)
taking a ``virtual'' pushforward.  Bott and Atiyah's miracle is that
{\em both} work, and yield the same thing (see \cite{bott}, \cite{atiyah} and
\cite{cls}, where Atiyah's map $\al$ is called $\overline{\partial}$).

\bpoint{Desiderata}
\label{s:desiderata}In the process of making this precise, we describe
a category $\cGw$ (some sort of naive ``homotopy category of geometric spaces'') that is a natural home for earlier stabilization results in
algebraic geometry.  Our desiderata are as follows.  Notice that most
are satisfied by  finite type algebraic stacks with affine diagonal. 
\begin{itemize}
  \item {\em (over $\Z$)} The spaces should be ``defined over $\Z$'',
    and ``base change well''.
    \item {\em (complex realization)} Upon specializing to $\C$, we should be able to analytify to obtain an appropriate
      complex analytic notion. 
      \item {\em (homotopy axiom)}  We want ``vector bundle morphisms'' (morphisms $\pi: E
\rightarrow X$ where $\pi$ expresses $E$ as a vector bundle over $X$) and more generally affine bundle morphisms
to be isomorphisms.  
    \item  {\em (topological realization)} The analytic object should have a homotopy type (and in
      particular, cohomology and homotopy groups).
    \item {\em (cohomology groups with structure)}  Different types of cohomology groups,
      with various structures (Galois structure, Hodge structure) should be defined. (Because of the homotopy axiom,
      Tate twists will have to be trivial, however.)
    \item {\em (Chow realization)} The spaces should have Chow groups and Chow rings.
      \item {\em (Grothendieck ring realization)} The spaces should
        have a class in the (completed
        localized) Grothendieck ring of varieties.
        \item {\em (motivic spaces)}  Over a perfect field, there should be a functor to the homotopy category of  motivic spaces.
          \end{itemize}

          The reader is welcome to work out how other potential desiderata
          are also obtained by this construction (\'etale cohomology, \'etale homotopy type, etc.).
          
We make no claim about our category  being the only setting (or
          even the ``right setting'') for these arguments.  If
          anything, {\em we emphasize that the arguments are robust
            (or simple-minded) enough that we don't need much care in
            choosing the setting for them.}  In particular, the spaces relevant
          for Bott periodicity are smooth, but one may want to allow singularities
          of various sorts.   See \S \ref{s:variations} for brief discussion of potential variations.

          \bpoint{Main result}  
          Bott periodicity is about isomorphisms of homotopy types, so we work in a robust and simple  category $\cGw$ which extends this notion to an algebraic setting:  there is a functor to the category of homotopy types.
          We will define $\Omega^2_{\alg ,d}(BGL(r))$  (\S \ref{d:final}) as, roughly, the moduli space of rank
          $r$ degree $d$ vector bundles on $\proj^1$, trivialized at the
          point $\infty \in \proj^1$.  Its ``homotopy type'' over $\C$ is indeed $\Omega^2_d(BGL(r))$.
          ($\Omega^2_{\alg, 0}(BGL(r))$ is  ``essentially'' the affine Grassmannian for $SL(r)$.)      It is straightforward to show
          $\Omega^2_{\alg, d}(BGL(r)) \cong  \Omega^2_{\alg, d'}(BGL(r))$ for all $d, d' \in \Z$ (see \S \ref{isor}). We have a map $\Omega^2_{\alg, d}(BGL(r)) \rightarrow \Omega^2_{\alg, d}(BGL(r+1))$
           obtained by sending the vector bundle $\cE$ to
          $\cE \oplus \oh$.  $BGL$ and $\Omega^2_{\alg, d}(BGL)$  similarly will be elements of $\cGw$ directly extending the corresponding notions in topology (see \S \ref{d:BGL} and \S \ref{ex:bgl}).

          \tpoint{Theorem (Bott periodicity, informal statement)}  {\em
            \label{t:bpis}We have an isomorphism 
            $$
             \Omega^2_{\alg, d}(BGL) \overset \sim
            \longrightarrow BGL$$  
            in $\cGw$. Upon specializing to $\C$,
            analytifying, and taking homotopy type, this specializes
            to the traditional Bott periodicity isomorphism. Over a perfect field, this induces an isomorphism of homotopy types of motivic spaces, which extends Bott perodicity in motivic spectra.}

In particular, Theorem~\ref{t:bpis} gives (yet) another proof of the classical statement of Bott periodicity.
For precise statements and proofs, see Theorem~\ref{t:bpnew} for the first part,
 Theorem~\ref{t:bp2} for the second, XXX for the third.

 This approach
has  additional consequences in fixed rank  $r$ ({\em i.e.}, without taking the limit $r
\rightarrow \infty$).    
The map from (the homotopy type of) $\Omega^2_{\alg ,d}(BGL(r))$ to $\Omega_d^2 BGL(r)$  is a homotopy equivalence (Theorem~\ref{hequiv}). 
(This approximation  of topologically-defined spaces by
 algebraically- or holomorphically-defined spaces is part of a
venerable tradition, for example \cite{segal}, \cite{cohensegal}, 
        \cite[Thm.~1.2]{milgram}, etc.)
        Furthermore, 
the (integral) cohomology ring
of $\Omega_d^2 BGL(r)$ is isomorphic to the Chow ring of the corresponding algebraic object, the moduli space of degree $d$ rank $r$ vector bundles on $\pp^1$ trivialized at $\infty$ (Theorem~\ref{algcoh}).

\epoint{What happened to $\Omega^1 BU \overset \sim \longleftrightarrow U$ and $\Omega^1 U \overset \sim \longleftrightarrow BU$?}
The ``single-loop'' statements are implicit here,  as well as one can
reasonably hope.
Philosophically:  the theory of motivic spaces teaches us that there
are two different types of ``motivic'' $\Omega^1$:  $S^1$ and $\G_m$, whose smash
product is $\proj^1$ (see for example \cite[Ex.~3.1.10]{morel}).
Then $\xymatrix{\Omega^1(BGL) \ar[r]^{\quad \sim} & GL}$, and more generally
$\xymatrix{\Omega^1 (BG) \ar[r]^{\quad \sim} & G}$ for any group $G$, should
be interpreted as a tautological statement corresponding to the
``$S^1$'' loop: $\xymatrix{\Hom(S^1, BG) \ar[r]^{\quad  \quad \sim} & G }$.

The other map 
$\xymatrix{\Omega^1(GL) \ar[r]^\sim & BGL}$ is 
special to $GL$, and can
be interpreted as follows.  A
(pointed) map $\G_m \rightarrow GL(r)$, interpreted as a transition or
clutching function, yields a rank $r$ vector bundle on $\proj^1$,
trivialized at $1 \in \G_m$.  For those bundles of degree $d$ (those
maps $\G_m \rightarrow GL(r)$ where the determinant has valuation $d$), we 
obtain  $\Hom_d(\G_m, GL(r)) \rightarrow \Omega^2_{\alg, d}(BGL(r))$, which
  composed with the map $\Omega^2_{\alg, d}(BGL(r)) \rightarrow  BGL$ of Theorem~\ref{t:bpis}, yields
  $\Hom_d(\G_m, GL(r)) \rightarrow BGL$.  Taking 
  $r \rightarrow \infty$ appropriately (in the same way as in Theorem~\ref{t:bpis}, translated into this language of clutching functions) yields the desired isomorphism.  (Caution is needed in taking colimits, cf.\ \S \ref{false}--\ref{d:final}.)

        \bpoint{Real Bott periodicity}
        Analogous generalizations of real Bott periodicity (with  the unitary groups
         replaced by orthogonal and symplectic groups) are the
        subject of forthcoming work with J. Bryan \cite{bryanvakil}.  Bryan and Sanders
        \cite{bryansanders} provide a geometric explanation for key
        parts of Bott periodicity which indicate how to proceed algebraically.

\bpoint{Standing assumptions} \label{s:staass}We work over an  arbitrary base  $B$.  (The cases of most of interest are presumably when 
$B$ is  the spectrum of a field or $\Z$.)    All statements are over $B$ (e.g., ``smooth'' means
``smooth over $B$''; ``morphisms $X \rightarrow Y$'' means morphisms
over $B$, etc.), and will obviously be preserved by base change by any $B'
\rightarrow B$.

  We use the language of locally free sheaves and vector
bundles interchangeably.  Stacks are assumed to be algebraic, and
locally of finite type.  All quotients are stack (or homotopy)
quotients, but to be clear, we will write $[X/G]$ rather than $X/G$.

By the category of topological spaces $\Top$, we mean the category of topological spaces homotopy equivalent to CW-complexes; and $\HoTop$ is
 the homotopy category of $\Top$ (the category of unstable homotopy types)  as usual.

\epoint{The starting point:  sufficiently nice smooth algebraic stacks}  \label{s:startingpoint}Let
$\SmArt$ be the $2$-category of Artin stacks
that are smooth, irreducible, have affine diagonal (their diagonal morphism is an affine morphism), and are
finite type (all over $B$, \S \ref{s:staass}).  (Notice:  for fixed $r>0$, the moduli space of rank $r$ vector bundles
on $\proj^1$ is not quasicompact and hence not in $\SmArt$, despite being in some sense the
main topic of this paper!  It is however the countable increasing union of open substacks, each  in $\SmArt$, which is all we will need, see Definition~\ref{d:xrdbetter}.)

        \bpoint{Acknowledgments} \label{ack}We thank Dave Anderson, Jim Bryan, Linda Chen,  Ralph Cohen, 
 Dan Dugger,        Aaron Landesman, Arpon Raksit, Emily Riehl, Isabel Vogt, and Ben Williams for many helpful
        conversations.  We thank Aravind Asok, Joel Kamnitzer, Thomas
        Lam, and Zhiwei Yun for useful comments.   We thank Leonardo
        Mihalcea for alerting us to the key reference \cite{stromme}.
We thank Maria Yakerson, Kirsten Wickelgren, and  Elden Elmanto  for patiently  explaining just enough motivic
homotopy theory to be dangerous; our        remaining
misunderstandings are due only to ourselves.

\section{Two classes of ``highly-connected'' morphisms:  $\connected$ and $\nicely$}

\label{s:niceconn}Suppose $k$ is a positive integer.
Define $\nicely$ to be the class of morphisms $X \rightarrow Y$ in $\SmArt$ that
can be written as a composition of morphisms of the form 
$$
\xymatrix{X \ar@{^(->}[r]^i & Z \ar[r]^a & Y}$$
where $i$ is an open embedding,  and $a$ is an affine bundle ($\A^n$-bundle for some $n$, in the smooth topology),
and over every point of $Y$, $i$ is a dense open embedding with complement of codimension greater than  $k$.  
Clearly ${\text{(nicely-$(k+1)$-conn)}} \subset \nicely$, and the class $\nicely$ is preserved by pullback and composition.  It is not clear whether the notion of $\nicely$  is local on the target.

\point \label{d:isocodimk}Define the  class $\connected$  of  morphisms in $\SmArt$
as the smallest collection of morphisms satisfying:
\begin{enumerate}
\item[(i)]  {\em (affine bundles)} The class $\connected$ contains all ``$\A^n$-bundle morphisms'' $\pi:  E \rightarrow
      X$.
   \item[(ii)] {\em (high codimension open subsets)}  The class $\connected$ contains all {\em open embeddings $X \hookrightarrow Y$
     where over every point $b \in B$, $\codim_{(Y \setminus X) \subset Y} > k$} (each irreducible
     component of the fiber of $Y \setminus X$ over $b$ is codim $>k$ in  the fiber of $Y$ over $b$).
   \item [(iii)] {\em (two-of-three)}  If $$\xymatrix{ X \ar[rr] \ar[dr] & & Y \ar[dl] \\ &
       Z}$$ is a  (2-)commuting triangle in $\SmArt$, and {\em two of
     the three} morphisms are $\connected$, so is the third.    
   \end{enumerate}
      The class $\connected$ is closed under composition (by (iii)),
      and preserved by change of base ring $B$ (\S \ref{s:staass}).
      Clearly
      $\text{($(k+1)$-conn)} \subset \connected$ and $\nicely \subset \connected$.
      Because of the two-of-three rule (iii), it is not clear whether the notion of $\connected$ is preserved by arbitrary pullbacks, unlike the case of $\nicely$.

Here is a first important example.

\tpoint{Proposition} \label{claimBGL} {\em The map $BGL(n) \rightarrow
  BGL(n+m)$,
given by sending the
vector bundle $\cV$ to $\cV \oplus \oh^{\oplus m}$, is 
${\text{($n$-conn)}}$.}

  \bpf  We deal first with $m=1$.
Let $e_1  \in \A^{n+1}$ be the first vector in the standard
basis. Acting on $e_1$ gives a map $GL(n+1) \to \A^{n+1} \smallsetminus
\{0\}$ with fiber the stabilizer $\operatorname{Stab}(e_1)$ of $e_1$. Thus the morphism $B
\operatorname{Stab}(e_1) \to BGL(n+1)$ has fiber $\A^{n+1} \smallsetminus \{0\}$,
and is hence is in
${\text{($n$-conn)}}$.  But
$\operatorname{Stab}(e_1) \cong GL(n) \ltimes \A^{n}$, so $BGL(n) \to
B\operatorname{Stab}(e_1)$ has fiber $\A^{n}$, and is thus also in
${\text{($n$-conn)}}$.   Thus their composition $BGL(n) \to BGL(n+1)$ is also in
${\text{($n$-conn)}}$.  Informal translation:
$$\! \! \xymatrix@C=33pt{
  BGL(n) \ar[r]^-{\text{rank $n$}}_-{\text{affine bundle}} & \, B
  \operatorname{Stab}(e_1)  \,
  \ar@{<->}[r]^--{\text{orbit-stabilizer}}_--{\text{theorem}} & \, [(\A^{n+1}
  \setminus \{ 0 \} )/ GL(n+1)] \\
\ar@{^(->}[r]^--{\text{
      open}}_--{\text{``codim
      $n$''}} & \, [\A^{n+1} / GL(n+1)] \,  \ar[r]^{\text{\quad rank
      $n+1$}}_{\text{\quad vector bundle}} & \, BGL(n+1)}$$
Proposition~\ref{claimBGL} then  follows by induction on $m$. \epf

\section{The categories $\cGsp \rightarrow \cGst \rightarrow \cGw$:  Cauchy sequences of Artin stacks, their spatial incarnations, and inverting ``infinitely connected morphisms''}

We  first  define the objects in two of our  categories of interest, $\cGsp$ and $\cGst$.

\epoint{Definition}
\label{d:gsp}The \underline{objects   of $\cGsp$} are diagrams of the form
\begin{equation}
  \label{eq:Gspdef}
\xymatrix{
  Y_1 \ar@{^(->}[r] \ar[d] &
  Y_2 \ar@{^(->}[r] \ar[d] &
  Y_3 \ar@{^(->}[r] \ar[d] &  \cdots \\
  X_1 \ar[r] &
  X_2 \ar[r] &
  X_3 \ar[r] &  \cdots
}
\end{equation}
(denoted $(Y_\bullet \rightarrow X_\bullet)$) 
where 
\begin{enumerate}
\item[(i)] the $Y_i$ are smooth finite type algebraic spaces over $B$ (e.g., smooth varieties if $B$ is a field), and $X_i \in \SmArt$; 
\item[(ii)] the morphisms in the top row  $Y_\bullet$ are {\em closed embeddings};
\item[(iii)] the vertical morphisms $Y_\ell \rightarrow X_\ell$ are $\text{(nicely-$k(\ell)$-conn)}$,
  with $k: \Z^+ \rightarrow \Z^+$ weakly increasing and unbounded; and
  \item[(iv)] 
  the morphisms $X_\ell \rightarrow X_{\ell+1}$
  in the bottom row are {\em separated} $\text{($k(\ell)$-conn)}$   morphisms in $\SmArt$ 
   with $k: \Z^+ \rightarrow \Z^+$ weakly increasing and unbounded.
\end{enumerate}
By the ``two-of-three'' axiom for $\connected$ (\S \ref{d:isocodimk}(iii)), the top morphisms $Y_\ell \rightarrow Y_{\ell+1}$ are also 
$\text{($k(\ell)$-conn)}$   with $k$ weakly increasing and unbounded.  ``All the arrows in \eqref{eq:Gspdef} get more highly connected the further to the right you go.''

The \underline{objects  of $\cGst$} are defined to be  diagrams \begin{equation}\label{eq:Gdef}
  \xymatrix{X_1 \ar[r] &  X_2 \ar[r] &  X_3 \ar[r] &  \cdots}
\end{equation}
(denoted $X_\bullet$)
in $\SmArt$ which can be extended to a diagram \eqref{eq:Gspdef}.
In particular, 
the morphisms $X_\ell \rightarrow X_{\ell+1}$
  in the bottom row are {\em separated}   $\text{($k(\ell)$-conn)}$ morphisms 
   with $k$ weakly increasing and unbounded.
  We will call
  \eqref{eq:Gspdef}  a ``spatial incarnation of \eqref{eq:Gdef}''.  We might call 
  a sequence \eqref{eq:Gdef} in $\SmArt$ of separated $\text{($k(\ell)$-conn)}$ morphisms 
     with $k$ weakly increasing and unbounded a ``Cauchy sequence of Artin stacks'',
  and call objects of $\cG$ ``Cauchy sequences of Artin stacks admitting a spatial incarnation''.

  The top row $Y_\bullet$ of \eqref{eq:Gspdef} will soon be naturally interpreted as a particularly nice ind-variety (or more correctly,
  ind-algebraic space over $B$).  
  We remark that the top row $Y_\bullet$ of \eqref{eq:Gspdef} is also in $\cGst$ because it has itself as a spatial incarnation.

  \epoint{Examples}  \label{s:threxa}The reader will readily think through the following  examples of elements  $\cGsp$. They 
  may help motivate  the definition.

  ``The ind-variety $\proj^\infty$ as a spatial incarnation of $B \G_m$'':
$$
\xymatrix{
  \proj^0  \ar@{^(->}[r] \ar[d] &
  \proj^1 \ar@{^(->}[r] \ar[d] &
  \proj^2 \ar@{^(->}[r] \ar[d] &  \cdots \\
  B \G_m \ar[r] &
  B \G_m \ar[r] &
  B \G_m \ar[r] &  \cdots
.}$$
More generally, ``the ind-variety $G(n, \infty)$ as a spatial incarnation of $BGL(n)$'':
\begin{equation}\label{eq:BGLnsp}
\xymatrix{
  G(n,n) \ar@{^(->}[r] \ar[d] &
  G(n,n+1) \ar@{^(->}[r] \ar[d] &
  G(n,n+2) \ar@{^(->}[r] \ar[d] &  \cdots \\
  BGL(n) \ar[r] &
  BGL(n) \ar[r] &
  BGL(n) \ar[r] &  \cdots
.}
\end{equation}

If $X$ is a smooth algebraic space with an action of $GL(n)$, the reader can readily extend the above ``mixing construction'' to exhibit $[X/GL(n)]$ as an element of $\cGst$, although we will not need this in what follows.  

\bpoint{Morphisms in $\cGsp$}
We now define the \underline{morphisms in $\cGsp$}.  Morphisms  in $\cGsp$ are, roughly speaking, morphisms of ``ladders'' of the form \eqref{eq:Gspdef}, except we allow ``rungs of the ladder'' of the source to map to ``later rungs'' of the target.  
More precisely,  the following data will (soon) determine such a morphism:
a weakly increasing unbounded function $m: \Z^{\geq 0} \rightarrow \Z^{\geq 0}$; and  morphisms
$Y_n \rightarrow Y'_{m(n)}$ and $X_n \rightarrow X'_{m(n)}$ such that the following cubes commute:
\begin{equation} \label{cube} \xymatrix{
\cdots \ar[r] & Y_n \ar[ddr] \ar[r] \ar[d] & Y_{n+1} \ar[d]  \ar[ddr]  \ar[r] & \cdots \\
\cdots \ar[r] & X_n  \ar[ddr] \ar[r] & X_{n+1}  \ar[ddr] \ar[r] & \cdots  \\
& \cdots \ar[r] &  Y'_{m(n)} \ar[r] \ar[d] & Y'_{m(n+1)} \ar[d] \ar[r] & \cdots  \\
& \cdots \ar[r] & X'_{m(n)} \ar[r] & X'_{m(n+1)} \ar[r] & \cdots.}
\end{equation}
If $M: \Z^{\geq 0} \rightarrow \Z^{\geq 0}$ is another weakly increasing function dominating $m$ ({\em i.e.}, such that $M(n) \geq m(n)$ for all $n$), then the above data determines similar data with $m$ replaced by $M$:  simply compose $Y_n \rightarrow Y'_{m(n)}$ with $Y'_{m(n)} \rightarrow Y'_{M(n)}$.

We define morphisms in $\cGsp$ as the above data, where the morphism ``with weakly increasing function $m$'' {\em is identified with} the  morphism ``with weakly increasing function $M$''.
Consequently, objects in $\cGsp$ depend only their on their ``asymptotic behavior''.
For example,  if $M: \Z^{\geq 0} \rightarrow \Z^{\geq 0}$ is any {\em increasing} function, then we have isomorphisms $( Y_n \rightarrow X_n)_n \overset \sim \longrightarrow (Y_{M(n)} \rightarrow X_{M(n)})_n$ in $\cGsp$.

\epoint{Comparing spatial incarnations} \label{comp} Given two spatial incarnations \eqref{eq:Gspdef} of an element $X_\bullet \in \cGst$, their ``fibered product'' over $X_\bullet$ is a ``spatial incarnation dominating both''.  More precisely, we have the following.

\tpoint{Proposition}  {\em \label{p:twocof}If  $(Y_\bullet \rightarrow X_\bullet)$ and $(Y'_\bullet \rightarrow X_\bullet)$ are two spatial incarnations of $X_{\bullet}$, then
\begin{equation} \label{spi}
  \xymatrix{
  Y_1\times_{X_1} Y'_1 \ar@{^(->}[r] \ar[d] &
  Y_2 \times_{X_2} Y'_2 \ar@{^(->}[r] \ar[d] &
  Y_3 \times_{X_3} Y'_3 \ar@{^(->}[r] \ar[d] &  \cdots \\
  X_1 \ar[r] &
  X_2 \ar[r] &
  X_3 \ar[r] &  \cdots
}
\end{equation}
is also a spatial incarnation of $X_{\bullet}$, {\em i.e.} is an element of $\cGsp$ as well.
}

\bpf  First, $Y_n \times_{X_n} Y'_n$ is a finite type algebraic space over $B$ (using that $X_n \in \SmArt$).  Next, suppose
$Y_n \rightarrow X_n$ is $\text{(nicely-$k(n)$-conn)}$
and
$Y'_n \rightarrow X_n$ is $\text{(nicely-$k'(n)$-conn)}$.
Then each vertical map in \eqref{spi} is a composition
$$
\xymatrix{
Y_n \times_{X_n} Y'_n \ar[rrr]^-{\text{(nicely-$k'(n)$-conn)}} & &  & Y_n \ar[rrr]^-{\text{(nicely-$k(n)$-conn)}} &  & & X_n,
}
$$
which is $\text{(nicely-$\min(k(n), k'(n))$-conn)}$. Moreover, $\min( k(n), k'(n))$ is weakly increasing and unbounded.
Finally, we show that each horizontal map in the top row of \eqref{spi} is a closed embedding by observing it as a combination of three closed embeddings, from the following diagram (using that closed embeddings are preserved by base change):
$$
\xymatrix{
  & Y_n \times_{X_n} Y'_n \ar[dl] \ar[r] & Y_n \times_{X_{n+1}} Y'_n \ar[dl] \ar[d] \ar[r] & Y_{n+1} \times_{X_{n+1}} Y'_n \ar[d] \ar[dr] \ar[r] & Y_{n+1} \times_{X_{n+1}} Y'_{n+1} \ar[dr] \\
  X_n \ar@{^(->}[r] & X_n \times_{X_{n+1}} X_n & Y_n \ar@{^(->}[r] & Y_{n+1} & Y'_n  \ar@{^(->}[r] & Y'_{n+1}  .}$$
\epf

Clearly this construction gives a morphism  $(Y_\bullet \times_{X_\bullet} Y'_\bullet \rightarrow X_\bullet) \rightarrow
(Y_\bullet  \rightarrow X_\bullet)$ in $\cGsp$.
The spatial incarnations of $X_{\bullet}$ form a 
partially ordered set with this ``domination'' order, and Proposition~\ref{p:twocof} shows that this
 partially ordered set is cofinal:  given any two spatial incarnations of $X_{\bullet}$, there is one dominating both.

More generally, given a morphism $\si$ in $\cGsp$ from a
spatial incarnation $(Y_\bullet \rightarrow X_\bullet)$ of
$X_\bullet$ to a spatial incarnation
$(Y'_\bullet \rightarrow X'_\bullet)$ of $X'_\bullet$, and any {\em
  other} spatial incarnations $(Z_\bullet \rightarrow X_\bullet)$ and
$(Z'_\bullet \rightarrow X'_\bullet)$ of $X_\bullet$ and $X'_\bullet$
respectively, there is not clearly a corresponding map of spatial incarnations
$(Z_\bullet \rightarrow X_\bullet) \rightarrow (Z'_\bullet \rightarrow
X'_\bullet)$. Nevertheless, there is a corresponding map $\si'$ of spatial incarnations dominating both.  This is made precise in the following diagram, in which all nonrectangular parallelograms are fibered products, and all dotted arrows are  $\nicely$ for appropriate $k$.
$$
\xymatrix{
  &  (Z_n \times_{X_n}  Y_n) \times_{X'_{m(n)}}   Z'_{m(n)} \ar@{.>}[dl] \ar@{.>}[dr] \ar[drrr]^{\si'_n} \\
  Z_n \times_{X_n} Y_n  \ar@{.>}[d] \ar@{.>}[dr] & & Y_n \times_{X'_{m(n)}}  Z'_{m(n)} \ar[rr] \ar@{.>}[dl]  & & Y'_{m(n)} \times_{X'_{m(n)}} Z'_{m(n)} \ar@{.>}[d] \ar@{.>}[dl] \\
  Z_n  \ar@{.>}[dr] & Y_n \ar@{.>}[d] \ar[rr]^{\si_n} & & Y'_{m(n)} \ar@{.>}[d] & Z'_{m(n)} \ar@{.>}[dl] \\
  & X_n \ar[rr] & & X'_{m(n)} 
}
$$

\bpoint{Morphisms in $\cGst$} We now define \underline{morphisms in $\cGst$}. A morphism $X_\bullet \rightarrow X'_\bullet$ is the data of a weakly increasing unbounded  $m: \Z^+ \rightarrow \Z^+$ along with morphisms $a_n: X_n \rightarrow X'_{m(n)}$ such that the squares
$$
\xymatrix{ X_n \ar[r] \ar[d]_{a_n} & X_{n+1} \ar[d]^{a_{n+1}} \\
  X'_{m(n)} \ar[r] & X'_{m(n+1)}}$$ commute, with the similar
identification of when two data are the same (in terms of the cofinal
partially order set of weakly increasing functions). With this definition, there is an obvious functor $\cGsp \to \cGst$ that sends $(Y_\bullet \to X_\bullet)$ to $X_\bullet$.
We define a \emph{lift} of a morphism $X_\bullet \to X'_\bullet$ in $\cGst$ to be a morphism in $\cGsp$ that is sent to $X_\bullet \to X'_\bullet$. In other words, a lift is the data of 
spatial incarnations $Y_\bullet \to X_\bullet$ and $Y_\bullet' \to X_\bullet'$ such that the cubes in \eqref{cube} commute.
Note that every morphism in $\cGst$ admits a lift. Indeed, given a morphism $X_\bullet \to X'_\bullet$, choose spatial incarnations $Y_\bullet \to X_\bullet$ and $Y_\bullet' \to X_\bullet'$ and consider the fiber product
$Y_\bullet \times_{X_{m(\bullet)}'} Y_{m(\bullet)}' \to X_\bullet$. Because (nicely $k$-conn) morphism are preserved by base change and composition, this 
is a spatial incarnation of $X_\bullet$. Thus, the diagram
\begin{center}
\begin{tikzcd}
Y_\bullet \times_{X_{m(\bullet)}'} Y_{m(\bullet)}' \arrow{d} \arrow{r} & Y_{m(\bullet)}' \arrow{d} 
\\
X_\bullet \arrow{r} & X_{m(\bullet)}'
\end{tikzcd}
\end{center}
defines a lift of our morphism.
Note that such lifts are not unique, but any two choices of lifts are dominated by another. Indeed, If $(Y_\bullet \to X_\bullet) \to (Y'_\bullet \to X_\bullet')$ and $(Z_\bullet \to X_\bullet) \to (Z'_\bullet \to X'_\bullet)$ are two lifts of the same morphism, then $Y_\bullet \times_{X_\bullet} Z_\bullet \to Y'_{m(\bullet)} \times_{X'_{m(\bullet)}} Z'_{m(\bullet)}$ is a lift dominating both.

Just as in $\cGsp$, objects in $\cGst$
depend only their on their ``asymptotic behavior''.

\epoint{Definition} \label{d:BGL}
Define $BGL \in \cGst$ by
$BGL(1) \rightarrow BGL(2) \rightarrow \cdots$, which has the following spatial incarnation:
\begin{equation}\label{eq:BGLsp}
\xymatrix{
  G(1, 2) \ar@{^(->}[r] \ar[d] &
  G(2, 4) \ar@{^(->}[r] \ar[d] &
  G(3, 6) \ar@{^(->}[r] \ar[d] &  \cdots \\
  BGL(1) \ar[r] &
  BGL(2) \ar[r] &
  BGL(3) \ar[r] &  \cdots .
}\end{equation}
More generally, 
  if  $a: \Z^+ \rightarrow \Z^+$ and $b: \Z^+ \rightarrow \Z^+$ are both increasing unbounded functions, then 
\begin{equation}\label{eq:BGLsp2}
\xymatrix{
  G(a_1,  a_1+b_1) \ar@{^(->}[r] \ar[d] &
  G(a_2,  a_2+b_2) \ar@{^(->}[r] \ar[d] &
  G(a_3, a_3+b_3) \ar@{^(->}[r] \ar[d] &  \cdots \\
  BGL(a_1) \ar[r] &
  BGL(a_2) \ar[r] &
  BGL(a_3) \ar[r] &  \cdots
}
\end{equation} is a spatial incarnation of (an object canonically isomorphic to) $BGL$.
The corresponding ind-variety is independent of $(a_i, b_i)$, and deserves the name  $G(\infty, \infty + \infty)$.

\bpoint{The final category $\cGw$}
Let $(WHE)$ be the class of morphisms $X_\bullet \rightarrow X'_\bullet$ in $\cGst$ such that for any (equivalently, one) representative of the morphism, the maps 
$X_n \to X_{m(n)}'$ are 
$\text{($b(n)$-conn)}$ for all $n$ for some
 unbounded weakly increasing  function $b: \Z^+ \rightarrow \Z^+$. 
 We note that in such a situation, there is a map of appropriate spatial incarnations $$\xymatrix{  Y_n \ar[r]^-{a_n} \ar[d] & Y'_{m(n)} \ar[d] \\
 X_n \ar[r] & X'_{m(n)}}$$
 where $a_n$ is $\text{($b(n)$-conn)}$ for some unbounded weakly increasing $b$.

When $B= \C$, we will define a functor $\top: \cG \to (HoTop)$  to the homotopy category of topological spaces and show (Lemma~\ref{l:twofour}) 
that if $\si \in (WHE)$, then $\top(\si)$ 
induces isomorphisms of all homotopy groups, and thus is a homotopy
equivalence by the Whitehead theorem.  Thus motivated, we define $\cGw$ as $(WHE)^{-1} \cGst$, the category defined by inverting the morphisms in $WHE$. (Concretely, the objects of $\cGw$ are the same as the objects of $\cGst$, and the morphisms in $\cGw$ are compositions of forwards and backwards arrows in $\cGst$ where backwards arrows are required to lie in $WHE$, and the composition of an arrow with its backwards is the identity.)

\epoint{Example}  As objects of $\cGw$, $BGL$ is isomorphic to the top row of \eqref{eq:BGLsp2}. More generally,
as objects of $\cGw$, any object of $\cGst$ is isomorphic to (the top row of) any spatial incarnation.\label{ex:bgl}

\section{Desiderata} \label{sec:desiderata}

We now show that our category $\cGw$ cheaply satisfies the Desiderata~\ref{s:desiderata}.
The directly algebro-geometric desiderata will be particularly easy, so we dispatch them first (\S \ref{s:agdes}).  The topological
desiderata necessarily require input from outside of algebraic geometry, and require more work (\S \ref{s:topdes}).  See
\S \ref{s:variations} for potentially useful extensions, for example to singular settings.

\bpoint{Algebro-geometric desiderata}
\label{s:agdes}
The following straightforward fact will yield all the algebro-geometric desiderata.

\epoint{Definition/Observation} Fix a positive integer $k$.   Suppose $F$ is a contravariant functor from the category of smooth algebraic spaces over $B$ to another category $\cC$,
  such that
  \begin{enumerate}
    \item [(i)] if $\pi: E \rightarrow X$ is an affine bundle, then $F(\pi): F(X) \rightarrow F(E)$ is an isomorphism, and 
    \item[(ii)]  if $i: X \hookrightarrow Y$ is an open embedding where $\codim_{(Y \setminus X) \subset Y}> k$, then $F(i): F(Y) \rightarrow F(X)$ is an isomorphism
    \end{enumerate}
    (cf.\ Definition~\ref{d:isocodimk}(i) and (ii)).
 We claim that $F$ can be extended in a natural way to a contravariant functor $F_{\mathrm{WHE}}: \cGw \rightarrow \cC$.  (This may be simpler for the reader to think through independently than to read.)

First we define a functor $F_{\mathrm{sp}}: \cGsp \to \cC$.
For any object $(Y_\bullet \to X_\bullet)$ of $\cGsp$, we define $F_{\mathrm{sp}}(Y_\bullet \to X_\bullet) = F(Y_i)$ for some choice of
$i$ such that the morphisms $Y_j \rightarrow Y_{j+1}$ are
$\connected$ for all $j>i$.
(Precisely and pedantically, the object $F_{\mathrm{sp}}(Y_\bullet \to X_\bullet) \in \cC$ is only defined up to canonical isomorphism, and we choose representatives simultaneously for all objects of $\cGsp$.  Detail-oriented readers can worry about foundational issues of choice.)
Having made such a choice of index $i$ for every object $(Y_\bullet \to X_\bullet)$, we then define $F_{\mathrm{sp}}: \cGsp \to \cC$ on morphisms
by sending $(Y_\bullet \to X_\bullet) \to (Y_\bullet' \to X_\bullet')$ to the
morphism which is the composition $$\xymatrix{ 
F(Y'_{i'})
 \ar@{<->}[r]^-\sim & F(Y'_{m(j)})  \ar[r] &  F(Y_j)\ar@{<->}[r]^-\sim & 
 F(Y_i)}
 $$
for some choice of $j$ such that $j > i$ and $m(j) > i'$.
    
Next we define a functor $F': \cGst \to \cC$ by choosing a spatial incarnation $Y_\bullet \to X_\bullet$ for each object $X_\bullet$, and then defining $F'(X_\bullet) = F_{\mathrm{sp}}(Y_\bullet \to X_\bullet)$. Note that for any other choice of spatial incarnation $Z_\bullet \to X_\bullet$, there is a canonical isomorphism $F_{\mathrm{sp}}(Y_\bullet \to X_\bullet) \overset \sim \longleftrightarrow F_{\mathrm{sp}}(Z_\bullet \to X_\bullet)$ induced by 
    \[ \xymatrix{ F(Y_i) \ar@{<->}[r]^\sim &  F(Y_I) \ar@{<->}[r]^-\sim &  F(Y_I \times Z_I) \ar@{<->}[r]^-\sim &  F(Z_I) \ar@{<->}[r]^\sim &  F(Z_j) }\] where $I > i, j$ is any index such that all vertical and horizontal arrows after the $I$th column are $\connected$ in both spatial incarnations. (By definition of $i$ and $j$, we have that the horizontal arrows in the top row after $i$ and $j$ are already $\connected$.)
    To define the functor $F'$ on morphisms in $\cGst$, we choose a lift $(Z_\bullet \to X_\bullet) \to (Z'_\bullet \to X'_\bullet)$ of every morphism.
    The source and target of the lift may not be the same as our chosen spatial incarnation, but we can compose with the canonical isomorphisms between alternate spatial incarnations to obtain a morphism between the desired source and target.

Finally, we claim that this $F': \cGst \rightarrow \cC$ descends to a functor $F_{\mathrm{WHE}}: \cGw \rightarrow \cC$. We must show that if $X_\bullet \to X_\bullet'$ lies in $(WHE)$, then $F'(X_\bullet \to X_\bullet')$ is an isomorphism.
Indeed, if $X_\bullet \to X_\bullet'$ lies in $(WHE)$,
then for our chosen lift $(Z_\bullet \to X_\bullet) \to (Z_\bullet' \to X_\bullet')$ the morphisms $Z_n \to Z_{m(n)}'$ are eventually $\connected$. In particular, if $j$ is the index we used to define $F_{\mathrm{sp}}$ of this morphism, then we can choose some $I \geq j$ such that $F(Z_I) \leftarrow F(Z_{m(I)}')$ is an isomorphism. Then the square below commutes
\begin{center}
\begin{tikzcd}
F(X_\bullet) \arrow[<->]{r} & F(Z_j) \arrow{d} & F(Z_{m(j)}) \arrow{d} \arrow{l} \arrow[<->]{r} & F(X_\bullet')\\
& F(Z_I)  & \arrow{l} F(Z_{m(I)}')
\end{tikzcd}
\end{center}
and the vertical arrows and bottom horizontal arrows are isomorphisms, so the top is too. It follows that $F'(X_\bullet \to X_\bullet')$ is an isomorphism, as desired.

Suppose $B$ is a field.     We can then apply this to ``Chow rings up to degree $k$'', and then by taking $k \rightarrow \infty$, we easily define a ``Chow ring'' contravariant functor from $\cGw$ to the category of graded rings.

    Similarly (with $B$ still a field), in the completed Grothendieck ring of stacks over a field  (introduced in
     \cite{e}), the class of an irreducible finite type stack $X$ with
     affine diagonal, normalized by dimension, $[X]/ \L^{\dim X}$,
     modulo classes of dimension less than $-k$, 
    is preserved by morphisms in $\connected$.
     (This can be made more precise: Morphisms in $\connected$ come
     with an additional invariant of ``relative dimension''.)

(With $B$ still a field, we can apply this to \'etale cohomology, using \cite[Exp.~XVI, Thme.~3.7]{sga4}; we omit the short argument because we shall not need it.)
 
\bpoint{Topological desiderata}
\label{s:topdes}
We first 
show that elements in $\connected$ are (topologically) highly
     connected, which motivates the terminology ``$\connected$''.  

       \tpoint{Lemma}  {\em \label{l:twofour}Suppose  $f: X \hookrightarrow Y$ is an open
       embedding in $\SmArt$, with $Z = Y \setminus X$, and
       $\codim_{Z \subset Y}
       >k$. 
       If $B = \C$, then $f$ induces an isomorphism  $\xymatrix{\pi_n(X^{an, top} ) \ar[r]^\sim &  \pi_n(Y^{an,
         top})}$  on
       homotopy groups for $n \leq 2k$. }

We remark that excision-type results are more difficult in
$\A^1$-homotopy theory; see \cite{asokdoran} for  examples of 
nontrivial results that are much more challenging to prove  than what we need here.

\bpf 
%
{\em (i)} First we show the result where $Z$ is smooth, following the suggestion of J. Bryan.

{\em (i.1)} We begin with the case where 
$Y$ is a smooth complex variety.  This case is
undoubtedly
well-known, but because we could not
find a reference, we include a sketch.
 The detailed explanation in \cite{benwilliams} is particularly satisfactory and is highly recommended to the reader. 
Suppose we have a map of the $n$-disc $f: D^n \to Y$ so that the boundary $S^{n-1}$ is sent to $X$.
If $n + \dim Z < \dim Y$, then the transversality theorem, in the form \cite[Corollary 2.5]{K}, implies there is a homotopy of $f$ which does not meet $Z$. Since $Z$ is closed and does not meet the boundary, we can arrange that this homotopy is constant on the boundary. 
This establishes that $\pi_{n-1}(X) \to \pi_{n-1}(Y)$ is injective. Meanwhile, given any map $S^n \to Y$, by considering the precomposition with the map $D^{n} \to S^{n}$ which contracts the boundary, this argument also establishes that $\pi_{n}(X) \to \pi_{n}(Y)$ is surjective.

{\em (i.2)}  Next, suppose $Y$ is the quotient of a smooth complex
variety $Y'$ by a smooth group $G$.  Define $X' = X \times_Y Y'
\subset Y'$, so 
$X=[X'/G]$.  By {\em (i.1)}, $\pi_n(X') \rightarrow \pi_n(Y')$ is an
isomorphism for $n \leq 2k$.   We compare the long exact sequence of
homotopy groups for  the fibrations $X = [X'/G]$ and $Y=[Y'/G]$ and induct on $n$:
$$
\xymatrix{
\pi_n(G) \ar[r] \ar[d]^= & \pi_n(X') \ar[d]^\sim \ar[r] & \pi_n(X)
\ar[d] \ar[r] & \pi_{n-1}(G) \ar[d]^= \ar[r] & \pi_{n-1}(X')
\ar[d]^\sim \\
\pi_n(G) \ar[r]  & \pi_n(Y') \ar[r] & \pi_n(Y)
\ar[r] & \pi_{n-1}(G) \ar[r] & \pi_{n-1}(Y')
}
$$
By the Five Lemma, we have that $\pi_n(X) \rightarrow \pi_n(Y)$ is an
isomorphism.  (The case $n=0$ is slightly different; it  follows because $\pi_0(X')$ surjects onto $\pi_0(X)$, and similarly for $Y'$ and $Y$.)

{\em (i.3)}  Next we deal with the case when $Y$ is an increasing union
of stacks of the form described in {\em (i.2)}.  Because homotopy
groups commute with filtered colimits, this case follows.

{\em (ii)}
   Next, we do the general case (where $Z$ is not necessarily smooth).  Let $Z_0=Z$, and inductively define
   $W_m = Z_{m-1}^{sm}$ (here ``sm''
   refers to the smooth locus of largest dimension) and $Z_m = Z_{m-1} \setminus W_m$.
   Let $X_m = Y \setminus Z_m$.  Then by {\em (i)}, the
   results hold for each of the inclusions $X \hookrightarrow X_1 \hookrightarrow X_2
   \hookrightarrow X_3  \hookrightarrow \cdots$, and hence $X \rightarrow X_m$ induces
   an isomorphism on $\pi_n$.
   Because homotopy groups commute with filtered colimits,  $X \rightarrow \cup X_m=Y$ also induces an isomorphism on
    $\pi_n$. \epf

 \tpoint{Corollary} \label{c:twosix}  {\em  
   If $B= \C$, then every $f \in \connected$ induces an isomorphism
   of homotopy groups $\pi_n$ for $n \leq 2k$. In particular,    the induced map of unstable homotopy types $f^{an,
    top}$ is $2k$-connected.
}
 
\bpf This is because such properties hold for the ``atomic'' morphisms
in the Definition~\ref{d:isocodimk} of $\connected$. \epf

\tpoint{Proposition}  {\em Suppose $B = \C$.

(a) Then we have a functor $\top: \cGsp \rightarrow \Top$ by sending \eqref{eq:Gspdef} to  the colimit (=union)  of topological realizations of the $Y_i$. 

(b)   This induces a functor $\top: \cGst \rightarrow \HoTop$ (to which we abusively give the same name).}

\bpf (a)
(Part (a)  motivates why $\cGsp$ was  defined to contain the data of a spatial incarnation.) Because our morphisms in spatial incarnations are closed embeddings, the colimit is an actual topological space; here we implicitly use the fact that finite type complex algebraic spaces are CW complexes, see \cite{kyle}.
(The colimit as a topological space is nothing exotic.  The colimit as a set is just the union of the sets in the sequence; and a subset of this set is open if its restriction to each set in the sequence is open.)
The functor $\top$ on morphisms is determined just by the universal property of colimit.

(b)  
For each object of $\cGst$, we choose any spatial incarnation to get a corresponding topological space, and use Lemma~\ref{l:twofour} to show that any two are canonically homotopic. The key facts that we use are that $\pi_n(\colim Y_\bullet) = \colim \pi_n(Y_\bullet)$ and the Whitehead theorem. 
  Moreover, given any morphism in $\cGst$, we can choose any lift to determine a corresponding map of topological spaces. Using the fact that any two choices of lift of a morphism are dominated by a third (\S \ref{comp}), we see that any two lifts determine the same map of homotopy types. The reader will readily verify that once we make the choice of where the objects go, even though there are many lifts of a morphism to spatial incarnations, the resulting morphism between homotopy types is uniquely determined. 
 \epf

By Lemma \ref{l:twofour} and the Whitehead theorem, the functor $\top$ sends any morphism in $WHE$ to an isomorphism in $(HoTop)$. Thus, the functor $\top: \cGst \rightarrow (HoTop)$  factors through $\top: \cGw \rightarrow (HoTop)$ (again abusively using the same name $\top$ for both functors).
These relationships are summarized in the following diagram.
\begin{equation}
\label{eq:RaviAug30}
\xymatrix{\cGsp \ar[r] \ar[d]_{\top} & \cGst \ar[r] \ar[dr]_-{\top} & \cGw \ar[d]^{\top} \\
  \Top \ar[rr] & & \HoTop.}\end{equation}

\epoint{Hodge theory}
None of the spaces in this paper have nontrivial Hodge structures, but the above arguments show that the cohomology groups on elements of $\cGw$ can be endowed with Hodge structures.

\bpoint{Desiderata on motivic spaces}

We record here properties proved by Church in Appendix~\ref{appendix}.

\tpoint{Theorem} {\em  Suppose $B= \Spec k$, where $k$ is a perfect field. \label{t:49}
\begin{enumerate}
\item[(a)] The usual functor from varieties to motivic spaces $Spc(k)$, postcomposed with  $\Sigma^2$,  extends to  a commuting diagram (cf.\ \eqref{eq:RaviAug30})
$$\xymatrix{
\cGst \ar[r] \ar[d] & \cGw  \ar[d] \\
 {\mathrm{Spc}}(k) \ar[r]  & {\mathrm{HoSpc}}(k)
}$$
where ${\mathrm{Spc}}(k)$ is the ($\infty$-)category of motivic spaces over $k$, and ${\mathrm{HoSpc}}(k)$ is its homotopy category.
\item[(b)] If we restrict to the subcategory built only out of $\A^1$-$1$-connected objects (i.e., the $X_i$ and $Y_i$ are $\A^1$-$1$-connected), the same holds without $\Sigma^2$; we do not need the double-suspension.
\end{enumerate}}

The left vertical arrow of Part (a) is explained in Appendix~\ref{appendix} starting with \S \ref{s:lvcs}.  The right vertical arrow follows by universal property of ($\infty$-category) localization, since $\cGw$  obtained from $\cGst$ by inverting morphisms which are equivalences in ${\mathrm{HoSpc}}(k)$, cf.\ \cite[Prop.~5.2.7.12]{Lur09}.
Part (b) follows from \S \ref{section:refinement}.

\bpoint{Possible variations on this theme}
\label{s:variations}
To extend these ideas to encompass singular spaces, the right
definition is determined by what properties are desired for the
application in mind.  Losing smoothness may lose the ``Chow ring''
functor for example.  But the following setting is still suitable for
many algebro-geometric invariants.  Instead of $\SmArt$, work with
finite type Artin stacks (over the base $B$), allowing a chosen class
of singularities.  Use the same definition of $\nicely$.  For the
right definition of $\connected$, in Definition~\ref{d:isocodimk},
keep the affine bundle condition (i) and the two-of-three rule (iii),
but replace (ii) judiciously to include only open embeddings
$U \hookrightarrow Y$ whose complement $Z \rightarrow Y$ meets ``the
singularities of $Y$ in codimension $>k$'' (tailored to the type of
singularities allowed).

\section{Moduli spaces related to bundles on $\proj^1$}

\label{s:three}We first define the central spaces we will discuss.

\bpoint{Important setting of definitions and coordinates:  $[x:y] \in \proj^1$; $H_0$ and $H_\infty$; and  $\psi_{r,d}: \cX_{r,d} \rightarrow \cY_{r,d}$} \label{d:XY}Define $\cX_{r,d}$ to be the moduli stack of rank $r$, degree $d$ locally free
sheaves on $\proj^1$, trivialized at $\infty \in \proj^1$.  More
precisely, the Artin stack $\cX_{r,d}$ parametrizes the following: for any $B$-scheme $S$, the maps $S
\rightarrow \cX_{r,d}$ correspond to the data 
(i)  a locally free sheaf $\cE$  of rank $r$ and fiber degree $d$ on 
$\pi:\proj^1\times S \rightarrow S$, and (ii)  a
trivialization  $\tau: \oh^{\oplus r} \overset \sim \longrightarrow
\si_{\infty}^* \cE$, where $\si_{\infty}$ is the section ``$\infty$'' of $\pi$:
$$
\xymatrix{
&    \cE \ar@{-}[d] \\
H_0,   H_\infty \ar@{^(->}[r] &  \proj^1 \times S \ar[d]_\pi  \ar[r]^{pr_1} & \proj^1 \ar[d]\\
 &   S   \ar[r] & B 
.}$$
We choose coordinates $[x:y]$ on $\pp^1$, writing $\infty \coloneq [1:0]$;  and
let $H_\infty = pr_1^* (\infty) = V(y)$ be the $\infty$-section on $\proj^1 \times S$, and (for future use) let $H_0= pr_1^* (0) = V(x)$ be the $0$-section on $\proj^1 \times S$.

Similarly, define $\cY_{r,d}$ to be the moduli stack of rank $r$ degree $d$ locally free
sheaves on $\proj^1$ ({\em without} trivalization at $\infty$).

\epoint{Remark} Notice that $\cX_{r,d}$ is naturally interpreted as
the moduli space of pointed morphisms
$\proj^1 \rightarrow BGL(r)$ of ``degree $d$'' (where
$\infty$ is mapped to the trivial(ized) rank $r$ bundle).  One might
write this as $\Maps_d^{\bullet}(\proj^1, BGL(r))$.
The space $\cX_{r,d}$ is clearly related to the affine Grassmannian
(for $GL(r)$).
Similarly, $\cY_{r,d} = \Maps_d(\proj^1, BGL(r))$.    

\epoint{Tautological (iso)morphisms between these spaces} \label{isor}The forgetful map 
$\psi_{r,d}: \cX_{r,d} \rightarrow \cY_{r,d}$ is a $GL(r)$-torsor.

Sending $\cE \mapsto \cE(H_0)$ (where $H_0$ is the zero section, \S \ref{d:XY}) defines isomorphisms $\cY_{r,d} \overset \sim \longrightarrow \cY_{r,d+r}$ and $\cX_{r,d} \overset \sim \longrightarrow \cX_{r,d+r}$ and the diagram
$$\xymatrix{
  \cY_{r,d} \ar[d] \ar[r]^-\sim &  \cY_{r,d+r} \ar[d] \\
    \cX_{r,d}  \ar[r]^-\sim &   \cX_{r,d+r} 
  }$$
clearly commutes. Similarly, given a bundle $\cE$ with a trivialization at $\infty$, we may define $\cE'$ to be the bundle whose sections are rational sections of $\cE$ that are regular except they are allowed a simple pole along the first component of the trivialization at $\infty$ (the ``positive modification" along the first vector of our trivialization).
Then $\cE \mapsto \cE'$ defines an isomorphism $\cX_{r,d} \overset \sim \longrightarrow \cX_{r,d+1}$. (We remark that this {\em cannot} lift to an isomorphism $\cY_{r,d} \overset \sim \longrightarrow \cY_{r,d+1}$ in general.  For example,   $\cY_{2,0} \not\cong \cY_{2,1}$, see
\cite[Cor.~5.6]{hannah}.  In \cite{hannah}, $\cY_{r,d}$ is denoted by $\cB^\dagger_{r,d}$.)

The morphism $\cX_{r,d} \rightarrow \cX_{r+1,d}$, sending a bundle $\cE$ to $\cE \oplus \oh$, will also be important in what follows.

\point
We recall the following useful result.

\tpoint{Str\o{}mme's Lemma \cite[Prop.\ (1.1)]{stromme}}
{\em \label{l:stromme}Suppose $\cE$ is a locally free sheaf on $\proj^1 \times S$, and
  $\pi: \proj^1 \times S \rightarrow S$ is the projection.  If $R^1
  \pi_* ( \cE(-1)) = 0$, then we have an exact sequence of (locally free) coherent sheaves on $\proj^1 \times S$:
\begin{equation}\label{eq:stromme}
\xymatrix{0 \ar[r] &  (\pi^* \pi_* \cE(-1)) (-1) \ar[r] &  \pi^* \pi_* \cE
\ar[r] &  \cE \ar[r] &  0.}\end{equation}}

The proof is a straightforward local calculation.  (It readily
generalizes to the case where the $\mathbb{P}^1$-bundle is not
trivial \cite[Lemma 3.1]{degen}, although we will not need this.)

\epoint{Remark}
Quillen \cite[Prop.~1.11]{Quillen} has in some sense a generalization of  the Str{\o}mme sequence to regular vector bundles on $\pp^r$-bundles.
(Given a $\pp^{r-1}$-bundle $X$ with structure map $\pi: X \to S$, a vector bundle $\cE$ on $X$
is called \emph{regular} if $R^q\pi_* \cE(-q) = 0$ for all $q > 0$.)   In particular, in $K$-theory, any regular vector bundle on $X$ is an alternating sum of twists of $\pi^*$'s of vector bundles on $S$. This is the key ingredient in proving the isomorphism $K(X) \cong K(S)^{\oplus r} $, called the ``Projective Bundle Theorem" in \cite[Theorem 2.1]{Quillen}.
The $r = 1$ case says $K(\pp^1 \times S) \cong K(S) \oplus K(S)$, a statement called the ``Fundamental Product Theorem'' in \cite[\S 2.1]{Hatcher}, which is the core of the proof of the statement of ``Bott periodicity in $K$-theory'' given there (which is not  the same as the Bott periodicity discussed here, which deals with spaces).

\point We return to our discussion of Str\o{}mme's Lemma.
By the Cohomology and Base Change Theorem, the exact sequence \eqref{eq:stromme} is
functorial in $S$ (precisely: given $S' \rightarrow S$, the
construction applied to the pullback of $\cE$ is canonically isomorphic to the
pullback of the construction applied to $\cE$).

If $\cE$ is rank $r$ and fiber degree $D$, then
$\cA_{-1} \coloneq \pi_* \cE(-1)$ has rank $D$ and $\cA_0 \coloneq
\pi_* \cE$ has rank $D+r$, and  \eqref{eq:stromme} tells us that  $\cE$ can be recovered from the
cokernel of the full rank ({\em i.e.} injective with locally free quotient)
map $a: \pi^* \cA_{-1}(-1) \rightarrow \pi^*\cA_0$.

This construction is reversible:

\tpoint{Reverse Str\o{}mme Lemma \cite[Lem.~3.2]{hannah}} \label{l:rsl}
{\em  Let $\pi: \proj^1 \times S \rightarrow S$ be the projection.
  Suppose $\cA_{-1}$ and $\cA_0$ are rank $D$ and $D+r$ (respectively)
  locally free sheaves on $S$, and $a: \pi^* \cA_{-1}(-1) \rightarrow
  \pi^* \cA_0$ is ``full rank'' (injective with locally free quotient).  Let $\cE = \coker a$.  Then
  $R^1 \pi_* \cE(-1)=0$, and we have canonical isomorphisms $\cA_{-1}
  \overset \sim \longleftrightarrow \pi_* (\cE(-1))$ and $\cA_0
  \overset \sim \longleftrightarrow \pi_* \cE$, making the following diagram commute.
\begin{equation}\label{eq:stromme2}\xymatrix{0  \ar[r] &  (\pi^*
    \cA_{-1} ) (-1) \ar[r]^{\quad a}  \ar[d]^\sim & \pi^* \cA_0 \ar[r] \ar[d]^\sim  &
    \cE \ar[r] \ar@{=}[d] & 0 \\
0 \ar[r] &  (\pi^* \pi_* \cE(-1)) (-1) \ar[r] &  \pi^* \pi_* \cE
\ar[r] &  \cE \ar[r] &  0 & \text{ (exact sequence \eqref{eq:stromme})}
  }\end{equation}}

\noindent {\em Proof (following \cite{hannah}).}  Let $\cE = \coker(a)$.  We use Cohomology and Base Change
repeatedly.
From
$$
\xymatrix{0 \ar[r] & ( \pi^* \cA_{-1})(-2) \ar[r]^{a(-1)} & (\pi^*
  \cA_0 )(-1) \ar[r] & \cE (-1) \ar[r] & 0}$$
and applying $\pi_*$, we see $R^1\pi_* (\cE(-1)) = 0$.
Applying
$\pi_*$ to \eqref{eq:stromme2}, we obtain
$$\xymatrix{ 0 \ar[r] & \cA_0 \ar[r]^\sim & \pi_* \cE \ar[r] & 0.}$$
Then \eqref{eq:stromme2} may be rewritten as
\begin{equation}\xymatrix{0  \ar[r] &  (\pi^* \cA_{-1} ) (-1)
    \ar[r]^{\quad a}
    & \pi^* \pi_* \cE  \ar[r] &
    \cE \ar[r] & 0 }\end{equation}
with the morphism $\pi^* \pi_* \cE \rightarrow \cE$ the same as in
\eqref{eq:stromme}, thereby identifying
$$\xymatrix{(\pi^* \cA_{-1})(-1) \ar@{<->}[r]^--\sim & \pi^* ( \pi_*
\cE(-1)) (-1),}$$ as they are both kernels of the same morphism.
Twisting this isomorphism by $\oh(1)$ and applying $\pi_*$, we are
done.
\epf

\bpoint{The open substacks $U_{r,d}(m) \subset \cX_{r,d}$ and $V_{r,d}(m) \subset \cY_{r,d}$ (following Bolognesi and Vistoli \cite{bv})} 
\label{s:bv}Suppose $m \in \Z$.  Define open substacks
$$
\xymatrix{
  U_{r,d}(m) \ar@{^(->}[d]_{\text{open}} \ar[r] &  V_{r,d}(m) \ar@{^(->}[d]^{\text{open}} \\
  \cX_{r,d} \ar[r] & \cY_{r,d}}
$$
as the substacks where the corresponding bundle $\cE(m)$ (of degree $d+mr$) is globally generated.
Equivalently, $U_{r,d}(m) \subset \cX_{r,d}$ and $V_{r,d}(m) \subset \cY_{r,d}$ are the open substacks corresponding to
vector bundles whose splitting type has summands of degree at least $-m$.
Equivalently again, they are the loci where $R^1 \pi_* \cE(m-1)=0$.  Clearly
$\cX_{r,d} = \cup_m U_{r,d}(m)$ and $\cY_{r,d} = \cup_m V_{r,d}(m)$.

\tpoint{Proposition} {\em   Suppose $d + mr \geq 0$.
  The open embedding   $U_{r,d}(m) \hookrightarrow \cX_{r,d}$ is 
  $\text{($(d+mr)$-conn)}$.
Hence  $U_{r,d}(m) \hookrightarrow U_{r,d}(m+1)$ is 
$\text{($(d+mr)$-conn)}$.
}\label{codimc}

The condition $d+mr \geq 0$ ensures that $U_{r,d}(m)$ is nonempty.

\bpf
The locus of vector bundles whose splitting type has a summand of
degree $-m-1$ or less is the closure of the universal splitting locus
for splitting type $$\left( -m - 1, \left\lfloor \frac{d+m+1}{r-1} \right\rfloor,
  \ldots,\left\lceil \frac{d+m+1}{r-1} \right\rceil \right).$$ The codimension of this universal splitting locus is readily seen to be
\[  h^1\left(\mathbb{P}^1, \mathrm{End}(\oh(-m-1) \oplus \oh(\lfloor
  \tfrac{d+m+1}{r-1} \rfloor) \oplus \cdots \oplus \oh( \lceil
  \tfrac{d+m+1}{r-1} \rceil) ) \right) = d + mr + 1. \]
The second sentence of the proposition follows because this splitting locus still has codimension $d + mr+1$ when we restrict to the open $U_{r,d}(m+1) \subset \cX_{r,d}$.
\epf

Let $\xymatrix{FR_{r,D}  \ar@{^(->}[r]^{\text{open} \quad} &  \A^{2D(D+r)}_{B}}$
parametrize the $D \times (D+r)$ matrices whose are entries are linear
forms (over $B$)  in the coordinates $x$ and $y$ of $\proj^1$, and
which are full rank at all points of $\proj^1$.   ($FR$ is for ``Full Rank''.)
The Reverse Str\o{}mme Lemma~\ref{l:rsl} (using $D=d +mr$) immediately
gives the following.

\tpoint{Proposition} \label{vprop}  {\em  We have an isomorphism $$V_{r,d}(m) \overset \sim \longleftrightarrow
\left[  FR_{r,d+mr}/  ( GL(d+mr) \times GL(d+(m+1)r)) \right].$$}\label{p:37}

\epoint{Aside}  The following are easily seen to hold, but we omit proofs as we will not use them.
\begin{enumerate}\item[(i)] The closed subset $\A^{2D(D+r)} \smallsetminus
  FR_{r,D}$ has codimension $r$ in $\A^{2D(d+r)}$.
\item[(ii])  If $d+mr \geq 0$, the induced map
  $$V_{r,d}(m) \rightarrow
  BGL(d+mr) \times BGL(d+ (m+1)r)$$ is
$\text{(nicely-$(r-1)$-conn)}$.
(The hypothesis $d+mr \geq 0$ ensures that $V_{r,d}(m)$ is nonempty.)
\end{enumerate}

\bpoint{Definition:  $\tU_{r,d}(m)$}
\label{d:urd}Next, we define an affine bundle  $\tU_{r,d}(m) \to U_{r,d}(m)$.
An object of $\tU_{r,d}(m)$ over $S$ is the data of an object of $U_{r,d}(m)$ over $S$ together with a lift of the trivialization $\tau$ as in the diagram below.
\begin{center}
\begin{tikzcd} \\
& \sigma_{\infty}^*\pi^* \pi_* \cE(m) \arrow{d}
\\
\oh^{\oplus r} \arrow[dashed]{ur}{\tilde{\tau}} \arrow{r}[swap]{\tau} & \sigma_\infty^* \cE(m).
\end{tikzcd}
\end{center}
(Above, we have a natural identification $\sigma_\infty^*\cE \cong \sigma_\infty^*\cE(mH_0)$, so the
trivialization of $\sigma_\infty^* \cE$ canonically gives a trivialization of $\sigma_\infty^*\cE(m)$, which we are also denoting by $\tau$ above.)
For each of the $r$ vectors in $\tau$, any two lifts differ by an element in the $(d+mr)$-dimensional kernel of $\pi^*\pi_* \cE(m) \to \cE(m)$ at the fiber over $\infty$.
Therefore, $\tU_{r,d}(m) \to U_{r,d}(m)$ is an
$\mathbb{A}^{r(d+mr)}$ bundle.

\point \label{s:algspa}Note that the Artin stack $\tU_{r,d}(m)$ has trivial isotropy groups (the lifts $\tilde{\tau}$ of the sections $\tau$ kill all automorphisms), so it is actually an algebraic space.

\epoint{$\tU_{r,d}(m)$ as a quotient}
We now describe $\tU_{r,d}(m)$ as a quotient of an open subset of affine space.
Let $A$ be a free $B$-module of rank $d+mr$.
Consider the free $B$-module \begin{equation} \label{vsp} \Hom(A, A) \oplus \Hom(A, \oh^{\oplus r}).
\end{equation} Now  $GL(A)= GL(d+mr)$ acts on $\Hom(A, A)$ by conjugation and on $\Hom(A, \oh^{\oplus r})$ by right multiplication.
Given an element $(\alpha, j)$ with $\alpha \in \Hom(A, A)$ and $j \in \Hom(A, \oh^{\oplus r})$, we associate a block matrix of linear forms
\begin{equation} \label{matform} \left( \begin{matrix}Ix - \alpha y \\  jy \end{matrix} \right) \in \Hom(A \otimes \oh(-1), A \otimes \oh \oplus \oh^{\oplus r}).
\end{equation}

\point
\label{Xcodim}Let $X \subset \Hom(A, A) \oplus \Hom(A, \oh^{\oplus r})$ be the open subset of $(\alpha, j)$ so that the associated matrix of linear forms does not drop rank anywhere along $\mathbb{P}^1$.  A straightforward dimension count shows that the complement of $X \subset \Hom(A, A) \oplus \Hom(A, \oh^{\oplus r})$ is codimension $r$.

\tpoint{Lemma} {\em We have $\tU_{r,d}(m) = [X/GL(A)]$. In particular the pushforward map $\tU_{r,d}(m) \to BGL(d+mr)$ is $\text{(nicely-$(r-1)$-conn)}$.} 
\label{l:xgla}

\bpf
This is similar to Proposition \ref{vprop}, except now the $r$ sections that generate the fiber at infinity give  a splitting $\pi_* \cE(m) \cong \pi_* \cE(m-1) \oplus \oh^{\oplus r}$. Setting $\cA = \pi_*\cE(m-1)$, the matrix in the left-hand map of the Str{\o}mme sequence lies in \[\Hom(\pi^*\cA(-1), \pi^*\cA \oplus \oh^{\oplus r}).\]
This is the inclusion of the sections that vanish at infinity. So, 
at infinity ({\em i.e.} $y = 0$), the matrix is the identity between $\pi^*\cA$ and $\pi^*\cA$, and zero on the $\Hom(\cA, \oh^{\oplus r})$ factor. 
In other words, the matrix takes the form in \eqref{matform}.

The second sentence  then follows from \S \ref{Xcodim}.
\epf

\epoint{Remark}  If $B$ is a field of characteristic $0$, e.g., $B = \C$, then $\tU_{r,d}(m)$ is in fact a variety.   This may be seen  using geometric invariant theory:  under the action of $GL(A)$ on 
$\Hom(A, A) \oplus \Hom(A, \oh^{\oplus r})$, $GL(A)$ acts on $X$ freely (which is equivalent to the objects parametrized having no nontrivial automorphisms, \S \ref{s:algspa}).
Thus $X$ is contained in the stable locus of the action.   Geometric invariant theory quotients are schemes.  (It can be checked that  $X$ is precisely the GIT-stable locus, \cite{bryanpc}.)

\epoint{Remark:  Quiver interpretation}
\label{r:qi}The key construction of $\tU_{r,d}(m)$ may seem ad hoc and unmotivated, but can be cleanly understood in terms of representations of the following quiver (cf.\ Figure~\ref{f:2}).
This was explained to us by Jim Bryan, and may be discussed further in \cite{bryanvakil}.$$
\xymatrix{
  A &  \bullet \ar[r]^j
   \ar@(ul,ur)^{\al}
  & \bullet &   B^{\oplus r}
  }
$$

\section{Interpretation of moduli spaces of bundles in  $\cGw$, and proof of Bott periodicity}

We now interpret the statements of the previous section in terms of our new category, and discover that the statement of complex Bott periodicity just falls out.

We begin by recalling that we have defined  $BGL(n)$ (and
in particular $B \G_m$) and $BGL$ as elements of $\cGst$, see \S \ref{s:threxa} and \ref{d:BGL}.  We emphasize that we cannot just lazily
define $BGL(n)$ as the sequence
$BGL(n) \rightarrow BGL(n) \rightarrow \cdots$; we need also to verify
that there exists a spatial incarnation (which was exhibited in
\eqref{eq:BGLnsp}).

Next, we wish to describe the space of rank $r$ degree $d$ bundles on $\proj^1$ as an object of $\cGst$, and our set-up helps us make the right choice.
The {\em wrong} definition is $\cX_{r,d} \rightarrow \cX_{r,d} \rightarrow \cdots$, because such a sequence does not have a spatial incarnation, because $\cX_{r,d}$ is not quasicompact.  Instead we are led (in precise analogy with intuition from topology) to consider the increasing union of open subsets of $\cX_{r,d}$
\begin{equation}
  \label{eq:xrdbetter}
  \xymatrix{
  U_{r,d}(1) \ar@{^(->}[r] &
  U_{r,d}(2) \ar@{^(->}[r] &  U_{r,d}(3) \ar[r] &
  U_{r,d}(4) \ar@{^(->}[r] &  \cdots
}
\end{equation}
whose union is $\cX_{r,d}$ (\S \ref{s:bv}).   By Proposition~\ref{codimc}, $U_{r,d}(m) \hookrightarrow U_{r,d}(m+1)$ is
$\text{($(d+mr)$-conn)}$.

\tpoint{Proposition}\label{p:srxrdbetter}
{\em $$
\xymatrix{
  \tU_{r,d}(1) \ar@{^(->}[r] \ar[d] &
  \tU_{r,d}(2) \ar@{^(->}[r] \ar[d] &
  \tU_{r,d}(3) \ar@{^(->}[r] \ar[d] &  \cdots \\
  U_{r,d}(1) \ar@{^(->}[r] &
  U_{r,d}(2) \ar@{^(->}[r] &
  U_{r,d}(3) \ar@{^(->}[r] &  \cdots
}
$$
is a spatial incarnation of \eqref{eq:xrdbetter}.
}

\bpf The bottom row is \eqref{eq:xrdbetter} and the vertical rows are the affine bundles described in \S \ref{d:urd}. 
Recall that $\tU_{r,d}(m)$ is the moduli space of vector bundles $\cE$ on $\pp^1$ such that $\cE(m)$ is globally generated, together with a choice of sections $s_1, \ldots, s_r$ of $\pi_*\cE(m)$ lifting the trivialization $\cE$ at infinity.
The  horizontal morphisms in the top row are defined by sending $(\cE, s_1, \ldots, s_r)$ to $(\cE, x s_1, \ldots, xs_r)$. These horizontal morphisms are clearly closed embeddings:  $\tU_{r,d}(m)$ is identified in $\tU_{r,d}(m+1)$ by the condition that the sections all vanish at $H_0 = V(x)$.
\epf

\tpoint{Corollary} {\em The morphism  $\tU_{r,d}(m) \rightarrow \tU_{r,d}(m+1)$ (in Proposition~\ref{p:srxrdbetter}) is $\text{($(d+mr)$-conn)}$.}

\bpf  By Proposition~\ref{codimc}, $U_{r,d}(m) \hookrightarrow U_{r,d}(m+1)$ is
$\text{($(d+mr)$-conn)}$. Now apply the two-of-three rule (\S \ref{d:isocodimk}(iii)) to the squares in the diagram in Proposition~\ref{p:srxrdbetter}. \epf

\epoint{Definition}  Define $\Omega^2_{\alg, d}(BGL(r))$ (``algebraic-loops-2''), an object in $\cGst$, by
\eqref{eq:xrdbetter}:
\begin{equation}
  \label{eq:twelveprime}
  \xymatrix{\Omega^2_{\alg, d}(BGL(r)): & 
  U_{r,d}(1) \ar@{^(->}[r] &
  U_{r,d}(2) \ar@{^(->}[r] &  U_{r,d}(3) \ar[r] &
  U_{r,d}(4) \ar@{^(->}[r] &  \cdots
}
\end{equation}
Proposition~\ref{p:srxrdbetter} gives a spatial incarnation.
$\Omega^2_{\alg, d}(BGL(r))$ is the better version of $\cX_{r,d}$ in our situation.  We avoid calling it $\cX_{r,d}$ to prevent confusion.\label{d:xrdbetter}

\epoint{Remark} It will soon be relevant to notice that any family of degree $d$ rank $r$ vector bundles on $\proj^1$ parametrized by a {\em quasicompact} base $S$, thus corresponding to $S \rightarrow \cX_{r,d}$, necessarily factors through $U_{r,d}(m)$ for some $m$.  \label{qccomment}

\epoint{Examples} The reader is invited to describe the morphisms in $\cGst$ \begin{eqnarray}
  \Omega^2_{\alg, d}(BGL(r)) & \overset \sim \longleftrightarrow &  \Omega^2_{\alg, d+r}(BGL(r)), \nonumber \\
\Omega^2_{\alg, d}(BGL(r)) & \overset \sim \longleftrightarrow &  \Omega^2_{\alg, d+1}(BGL(r)),  \text{ and } \nonumber \\
\label{eq:plusone} \Omega^2_{\alg, d}(BGL(r)) &  \longrightarrow &  \Omega^2_{\alg, d}(BGL(r+1)) \end{eqnarray}
corresponding to the morphisms of stacks $\cX_{r,d} \overset \sim \longleftrightarrow \cX_{r, d+r}$, $\cX_{r,d} \overset \sim \longleftrightarrow \cX_{r,d+1}$, and $\cX_{r,d} \rightarrow \cX_{r+1,d}$ respectively, described in \S \ref{isor}.

\epoint{The sequence $Z$ (soon to be an object of $\cG$)} Fix a positive integer $m$.
The key player in what follows is the sequence of algebraic spaces
\begin{equation}
  \label{eq:fav}\xymatrix{Z: & \tU_{r,d}(m)  \ar[r]  & \tU_{r+1,d}(m) \ar[r]  & \tU_{r+2,d}(m)  \ar[r]  & \cdots .}
  \end{equation}
  The morphism
$\iota: \tU_{r+\ell,d}(m) \rightarrow  \tU_{r+\ell+1,d}(m)$ in \eqref{eq:fav}
corresponds to   adding a trivial summand $\oh$ (and adding its canonical section $1$ to the $r+\ell$ sections that generate the fiber at infinity), just as with  \eqref{eq:plusone}.

We wish to show that $Z$ is an object of $\cGst$.
(This goal will be completed by Proposition~\ref{p:ZinG}.)
We first show that the morphisms in $Z$ become increasingly connected.
Then, we construct a spatial incarnation of $Z$ in \S \ref{spcon}.

Consider the following commutative diagram where the top row is $Z$ (from \eqref{eq:fav}):
\begin{equation} \label{gs1}
\begin{tikzcd}
 \tU_{r,d}(m) \arrow{d} \arrow{r} & \tU_{r+1,d}(m) \arrow{d} \arrow{r} & \tU_{r+2,d}(m) \arrow{d} \arrow{r} & \cdots \\
 BGL(d+mr) \arrow{r} & BGL(d+m(r+1)) \arrow{r} & BGL(d+ m(r+2) ) \arrow{r} & \cdots.
\end{tikzcd}
\end{equation}
Above, the vertical maps are given by $\cE \mapsto \pi_*\cE(m-1)$. The maps on the bottom row are given by $V \mapsto V \oplus \oh^{\oplus m}$. The squares commute via the isomorphisms \[\pi_*((\cE \oplus \oh)(m-1)) \overset \sim \longleftrightarrow  \pi_*\cE(m-1) \oplus \pi_*\oh(m-1) \overset \sim \longleftrightarrow \pi_*\cE(m-1) \oplus \oh^{\oplus m}.\]
Each  bottom horizontal map $BGL(d+m(r+\ell)) \to BGL(d+m(r+\ell+1))$ of \eqref{gs1} is
${\text{($(d+m(r+\ell))$-conn)}}$ (by Proposition~\ref{claimBGL}).
The vertical maps $T_{r+\ell,d}(m) \to BGL(d + m(r+\ell))$ are (nicely-$(r+\ell-1)$-conn) (by Lemma~\ref{l:xgla}).
By the two-of-three rule (\S \ref{d:isocodimk}(iii)),  we therefore have the following.

\tpoint{Claim} {\em Assume $d + m(r+\ell) \geq r+\ell-1$. Then the top horizontal maps $\tU_{r+\ell,d}(m) \to \tU_{r+\ell+1,d}(m)$ of \eqref{gs1}, {\em i.e.} the morphisms in $Z$, are ${\text{($(r+\ell-1)$-conn)}}$.}\label{cl:conn}

\tpoint{Corollary to Claim~\ref{cl:conn}}  \label{c:c45}    {\em Assume $d + m(r+\ell) \geq r+\ell-1$. The morphism $U_{r+\ell,d}(m) \rightarrow U_{r+\ell+1,d}(m)$ is  ${\text{($r+\ell-1$-conn)}}$.}

\epoint{Constructing a spatial incarnation of $Z$} \label{spcon}
Although the $T_{r,d}(m)$ are algebraic spaces, $Z$ cannot serve as its own spatial incarnation because
the horizontal maps $T_{r+\ell,d}(m) \to T_{r+\ell+1,d}(m)$ are not closed embeddings. Nevertheless, they are \emph{locally} closed embeddings, as the next proposition will show. Thus by throwing out appropriate closed subsets, we will be able to construct our desired spatial incarnation.

Recall that a ``section of  $T_{r+1,d}(m)$ over $S$'' is  the data of $(\cE', s_1', \ldots, s_{r+1}')$ where $\cE'$ is a rank $r+1$ vector bundle on $\pp^1 \times S$ so that $\cE'(m)$ is relatively globally generated and $s_1', \ldots, s_{r+1}'$ are global sections of $\cE'(m)$ generating the fiber at infinity. (For clarity we write primes on notation for ``rank $r+1$ data'' and no primes on notation for ``rank $r$ data''.)
Let $S_{r+1,d}(m) \subset T_{r+1,d}(m)$ be the locally closed subspace defined by the following two locally closed conditions:

(i) $V(s_{r+1}') = m  H_0$  ({\em i.e.} the last section vanishes to order exactly $m$ along $H_0$ and nowhere else). Requiring that $s_{r+1}'$ vanish to at least a certain order would be a closed condition; asking to have a vanishing of an exact order $m$ is a locally closed condition.

(ii)
$V(s_1' \wedge \cdots \wedge s_r')$ (a closed
subscheme of $\proj^1 \times S$, the vanishing scheme of a section of a rank $r+1$ vector bundle) is a relative effective Cartier divisor of relative degree $d+mr$.
This is a locally closed condition on $T_{r+1,d}(m)$:  above $T_{r+1,d}(m)$, $V(s_1' \wedge \cdots \wedge s_r')$ gives a  closed subspace $V$ of $\proj^1 \times T_{r+1, d}(m)$. Consider the flattening stratification for this projective scheme $V$ over $T_{r+1,d}(m)$;
$S_{r+1,d}(m)$ is precisely the stratum corresponding to Hilbert polynomial $d+mr$.

\tpoint{Claim}
{\em  The morphism $\iota: T_{r,d}(m) \to T_{r+1,d}(m)$  in \eqref{eq:fav} is a locally closed embedding;  it induces an isomorphism  $T_{r,d}(m) \xrightarrow{\sim} S_{r+1,d}(m) \subset T_{r+1,d}(m)$.}

\bpf Let $(\cE, s_1, \ldots, s_r)$ be an object of $T_{r,d}(m)$ over $S$.
By definition, $\iota((\cE, s_1, \ldots, s_r)) = (\cE \oplus \oh, (s_1, 0), \ldots, (s_r,0), (0, x^m))$. 
Thus, any object $(\cE', s'_1, \dots, s'_{r+1})$ of $\iota(T_{r,d}(m))$ over $S$ satisfies $V(s'_{r+1}) = m H_0$.  Additionally, if $(\cE, s_1, \dots, s_r)$ is an object of $T_{r,d}(m)$ over $S$, then $s_1$, \dots, $s_r$ are global sections of $\cE(m)$ spanning $\cE(m)|_{\infty}$ (in every fiber), so $s_1 \wedge \cdots \wedge s_r$ defines a global section of $\det \cE(m)$ which is nonzero at infinity in every fiber, and is hence not the zero section on any fiber. Because $\det \cE(m)$ is a line bundle of relative degree $d+mr$, $V(s_1 \wedge \cdots \wedge s_r)$ is a relative effective Cartier divisor of relative degree $d+mr$. 
Thus,
$\iota((\cE, s_1, \ldots, s_r))$ is an object of $S_{r+1,d}(m)$ over $S$.  We have thus shown that $\iota: T_{r,d}(m) \rightarrow T_{r+1,d}(m)$ factors through the locally closed subspace $S_{r+1,d}(m)$.  It remains to describe the inverse $S_{r+1,d}(m) \rightarrow T_{r,d}(m)$.

Now suppose $(\cE', s'_1, \dots, s'_{r+1})$ is an object of $S_{r+1,d}(m)$ over $S$.
Let $\cF = \cE'(m)$, which is rank $r+1$, degree $d + m(r+1)$.
Let $D = V(s_1' \wedge \cdots \wedge s_r') \subset \pp^1 \times S$,
so $s_1' \wedge \cdots \wedge s_r'$ defines a section $s: \oh \to \wedge^r \cF$ that vanishes along $D$. By the definition of $S_{r+1,d}(m)$, we have that $D$ is a relative effective Cartier divisor of relative degree $d + mr$.
Recall the  natural identification $\wedge^r \cF \overset \sim \longleftrightarrow \det \cF \otimes \cF^\vee$. Taking the dual of our section $s$ and tensoring with $\det \cF$, we find that the map
\[\cF \to \det \cF \qquad \alpha \mapsto s_1' \wedge \cdots \wedge s_r' \wedge \alpha \]
vanishes precisely along $D$.
This implies that there is a surjection
\[\cF \to \det \cF(-D).\]
We define $\cE(m)$ to be the kernel, which is locally free. 
By construction, the sections $s_1', \ldots, s_r'$ of $\cF$ lie in the subbundle $\cE(m) \subset \cF$ and span the fibers along $H_\infty$. (In fact, the total space of the subbundle $\cE(m)$ is the closure of the span of these sections: $\cE(m) = \overline{\span(s_1', \ldots, s_r')} \subset \cF$.)
To identify the line bundle $\det \cF(-D)$, consider the section $s_1' \wedge \cdots \wedge s_r' \wedge s_{r+1}'$ of $\det \cF$. We certainly have the containment 
\begin{equation} \label{cont} D + mH_0 = V(s_1' \wedge \cdots \wedge s_r') + V(s_{r+1}') \subseteq V(s_1' \wedge \cdots \wedge s_r' \wedge s_{r+1}').
\end{equation}
Moreover, we know that $V(s_1' \wedge \cdots \wedge s_{r+1}')$ never contains a full fiber of $\pp^1 \times S \to S$ because $s_1' \wedge \cdots \wedge s_{r+1}'$ is non-vanishing along $H_\infty$. Since $\det \cF$ is relative degree $d + mr + m$, it follows that equality in \eqref{cont} holds in each fiber of $\pp^1 \times S \to S$ and thus equality holds in \eqref{cont}. We conclude that $\det \cF(-D) \cong \oh(mH_0)$ and we have a short exact sequence
\[\xymatrix{ 0 \ar[r] &  \cE(m) \ar[r] &  \cE'(m) \ar[r] &  \oh(m) \ar[r] &  0.} \]
Twisting down by $mH_0$ we obtain
\begin{equation} \label{ts}
\xymatrix{0 \ar[r] &  \cE \ar[r] &  \cE' \ar[r] &  \oh \ar[r] & 0.}
\end{equation}
We now see that $\cE \subset \cE'$ is a degree $d$ subbundle.

Next, we claim that \eqref{ts} splits. The splitting comes from $s'_{r+1}$.  Indeed, because $V(s'_{r+1}) = mH_0$,  the twist  $s'_{r+1} (- mH_0)$ defines a section $\oh \rightarrow \cE'$.  The composition
$$
\xymatrix{\oh \ar[rr]^-{s'_{r+1}(-mH_0)} & & \cE' \ar[r] & \cE' / \cE \cong  \oh}
$$ is nonzero because it is nonzero 
in the fiber at $\infty$.  For degree reasons, it must be an isomorphism, giving the desired splitting $\cE' = \cE \oplus \oh$. Correspondingly, we also have $\cE'(m) = \cE(m) \oplus \oh(m)$.
With respect to this splitting, we can write each $s_i' = (s_i, 0)$ where $s_i$ is a section of $\cE(m)$ (for $i \leq r$)  and $s'_{r+1} = (0, x^m)$.

Because $s_1', \ldots, s_r'$ span the fiber of the subbundle $\cE(m) \subset \cE'(m)$, we have that
$(\cE, s_1, \ldots, s_r)$ defines an object of $T_{r,d}(m)$ over $S$, which $\iota$ clearly sends to $(\cE',s_1', \ldots, s_{r+1}')$.

The composition of the maps described above
$$\xymatrix{  T_{r,d}(m) \ar[r] & S_{r+1, d}(m) \ar[r] & T_{r,d}(m) \\
(\cE, s_1, \ldots, s_r) \ar@{|->}[r]&  (\cE \oplus \oh, (s_1, 0), \ldots, (s_r, 0), (0, x^m)) \ar@{|->}[r] &  (\cE, s_1, \ldots, s_r)}$$
is clearly the identity.
\epf

\tpoint{Claim} {\em   For any nonempty open subset $T \subset T_{r,d}(m)$, 
$$\overline{\iota(T)} \setminus \iota(T) \subset T_{r+1,d}(m)$$ is closed of codimension at least $m(r+1) + d+mr+1$.}\label{claim2}

\bpf  Because $T_{r,d}(m)$ is irreducible, and $\iota$ is a locally closed embedding, for any nonempty open subset $T$,
$\overline{\iota(T)}$ is closed of dimension $\dim T_{r,d}(m) = r(d+mr)$.  Furthermore, $\iota(T)$ is a dense open subset of $\overline{\iota(T)}$.  Thus $\overline{\iota(T)} \setminus \iota(T)$ is closed of dimension at most $r(d+mr)-1$.  Meanwhile,
$\dim T_{r+1,d}(m) = (r+1)(d+m(r+1))$.  The codimension is thus as claimed. \epf

\epoint{Definition}  Let $T'_{1,d}(m) = T_{1,d}(m)$.  Inductively (in $r$) define $T'_{r,d}(m)$ by
$$
T'_{r+1, d}(m) \coloneq T_{r+1,d}(m) \setminus \left(
\overline{   \iota(T'_{r,d}(m))} \setminus
\iota(T'_{r,d}(m)) \right).$$
We will show that $T_{r,d}'(m) \subset T_{r,d}(m)$ is an open  substack using induction on $r$. It holds when $r = 1$. Assume we know $T_{r,d}'(m) \subset T_{r,d}(m)$ is open.
Since $\iota$ is a locally closed embedding $\iota(T'_{r,d}(m))$ (represented in blue below) is locally closed, so $\overline{   \iota(T'_{r,d}(m))} \setminus
\iota(T'_{r,d}(m))$ (represented in red below) is closed. Hence, $T'_{r+1,d}(m) \subset T_{r+1,d}(m)$ (both represented in grey) is open.
\begin{center}
\includegraphics[width=4in]{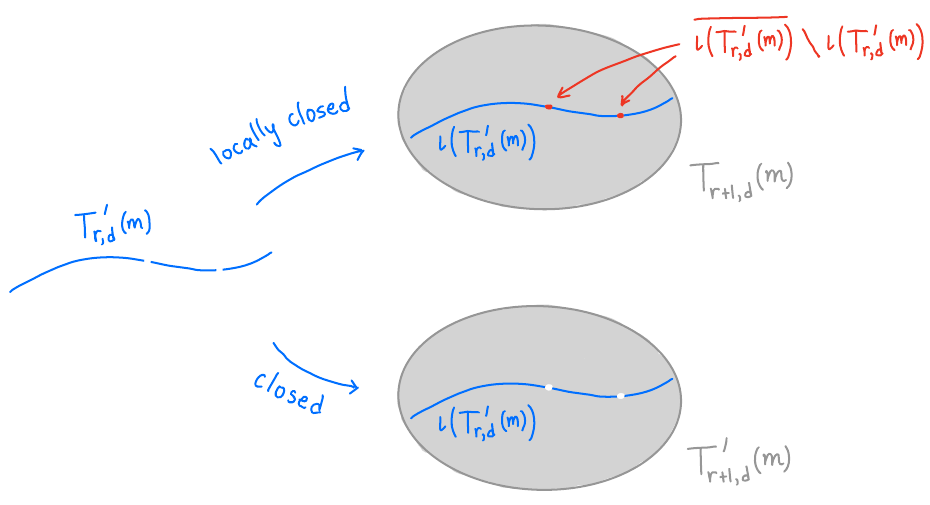}
\end{center}
By construction,
$T'_{r,d}(m) \rightarrow T'_{r+1,d}(m)$ are closed embeddings.
By Claim~\ref{claim2}, we have that $T'_{r+\ell,d}(m) \rightarrow T_{r+\ell,d}(m)$ is nicely-$(m(r+\ell)+d+m(r+\ell-1))$-conn.    Hence
$$
\xymatrix{\cdots \ar[r] & T'_{r,d}(m) \ar[dd] \ar[r] & T'_{r+1,d}(m) \ar[dd] \ar[r] & \cdots \\
& & & & \in \cG_{sp}. \\
\cdots \ar[r] & T_{r,d}(m) \ar[r] & T_{r+1,d}(m) \ar[r] & \cdots}
$$

We have thus finally proved what we have desired:

\tpoint{Proposition}{\em  $Z$ is an object of $\cGst$.}  \label{p:ZinG}

\point \label{s:fur}Furthermore, the bottom row of \eqref{gs1} is isomorphic to $BGL \in \cGst$.  
Recalling that the vertical maps of \eqref{gs1} are ${\text{(nicely-$(r+\ell-1)$-conn)}}$ (by Lemma~\ref{l:xgla}), we see that \eqref{gs1} thus describes  a morphism  $Z \rightarrow BGL$ in $\cGst$ that lies in $WHE$, {\em i.e.} defines an isomorphism in $\cGw$.

\bpoint{A false but motivating statement and proof of Bott periodicity} \label{false}
Complex Bott periodicity is a statement in the topological setting:  $\Omega^2 BU \overset \sim \longrightarrow \Z \times BU$, or equivalently (fixing a degree $d$ and using $U \sim GL$), the statement $\Omega^2_d BGL \overset \sim \longrightarrow BGL$.

Consider (with $m \in \Z^+$ as above) the commutative diagram
\begin{equation} \label{gs2}
\begin{tikzcd}
Z: & \tU_{r,d}(m) \arrow{d} \arrow{r} & \tU_{r+1,d}(m) \arrow{d} \arrow{r} & \tU_{r+2,d}(m) \arrow{d} \arrow{r} & \cdots \\
& \cX_{r,d} \arrow{r} &  \cX_{r+1,d}  \arrow{r} & \cX_{r+2,d} \arrow{r} & \cdots
\end{tikzcd}
\end{equation}
Each vertical map is a composition of an affine bundle $\tU_{r,d}(m) \to U_{r,d}(m)$ (\S \ref{d:urd}) with an open inclusion $U_{r,d}(m) \subset \cX_{r,d}$ whose complement has codimension $d+mr+1$ (Proposition~\ref{codimc}). (We need $m \geq 1$ to ensure that this codimension is unbounded as $r \rightarrow \infty$!)  Hence, the vertical map $\tU_{r+\ell}(m) \rightarrow \cX_{r+\ell, m}$ is
${\text{($(d+m(r+\ell))$-conn)}}$.
We already showed (Claim~\ref{cl:conn}) that the map $\tU_{r+\ell,d}(m) \rightarrow \tU_{r+\ell+1,d}(m)$ in the top row is ${\text{($(r+\ell-1)$-conn)}}$, so the map $\cX_{r+\ell, d} \rightarrow \cX_{r+\ell+1,d}$  in the bottom row is ${\text{($(r+\ell-1)$-conn)}}$
as well (by the two-of-three rule, \S \ref{d:isocodimk}(iii)):
$$
\xymatrix{ \tU_{r+\ell,d}(m)
  \ar[d]_{{\text{($(d+m(r+\ell))$-conn)}}}
  \ar[rr]^-{ {\text{($(r+\ell)$-conn)}}}
  & & \tU_{r+\ell+1,d}(m)
  \ar[d]^{{\text{($(d+m(r+\ell+1))$-conn)}}}
  \\
  \cX_{r+\ell,d}
  \ar[rr]^-{\therefore \;  {\text{($(r+\ell)$-conn)}}}
  & &   \cX_{r+\ell+1,d}.
}
$$
The top row $Z$ of \eqref{gs2} is isomorphic to $BGL$ in $\cGw$ (\S \ref{s:fur}), and
were the bottom row of \eqref{gs2}  an element of $\cGst$, \eqref{gs2} would give an isomorphism of $BGL$ with 
a plausible interpretation of $\colim\limits_{r \rightarrow \infty} \Omega^2_d(BGL(r))$.

Two things are wrong with this argument.  The bottom row is not an element of $\cGst$
because $\cX_{r, d}$ is not quasicompact,  and the phrase
$\colim\limits_{r \rightarrow \infty} \cX_{r,d}$ has no clear relation with $\Omega^2_d BGL$.

Instead, because we need to interpret $\cX_{r,d}$ as itself a sequence (a countable union of open subsets, \eqref{eq:twelveprime}), {\em i.e.} as $\Omega^2_{\alg, d}(BGL(r))$, we need to interpret/replace the bottom row of \eqref{gs2} as/by some sort of sequence obtained from the following two-dimensional diagram
\begin{equation} \label{grid}
\xymatrix{
  U_{r,d}(1) \ar@{^(->}[d] \ar[r] &   U_{r+1,d}(1) \ar@{^(->}[d] \ar[r] & U_{r+2,d}(1) \ar@{^(->}[d] \ar[r] & \cdots \\
  U_{r,d}(2) \ar@{^(->}[d] \ar[r]  &   U_{r+1,d}(2) \ar@{^(->}[d] \ar[r]  &  U_{r+2,d}(2) \ar@{^(->}[d] \ar[r]  & \cdots \\
  U_{r,d}(3) \ar@{^(->}[d]   \ar[r] &   U_{r+1,d}(3) \ar@{^(->}[d]   \ar[r] &  U_{r+2,d}(3) \ar@{^(->}[d]   \ar[r] & \cdots \\
\vdots & \vdots & \vdots & 
}
\end{equation}
where the $\ell$th column ($\ell \geq 0$) corresponds to $\cX_{r+\ell, d}$.
Now  all the arrows in this diagram are $\connected$ for various $k$, and roughly speaking they get more and more highly connected the further to the right and down you go:
\begin{equation}\label{eq:star}
\xymatrix{ U_{r+\ell, d}(m) \ar[rr]^{{\text{($(r+\ell)$-conn.)}}}_{\text{(Cor.\ \ref{c:c45})}}  \ar@{^(->}[d]_{\text{($(d+m(r+\ell))$-conn.)}}^{\text{(Prop.~\ref{codimc})}}
  & & U_{r+\ell+1, d}(m) \\
  U_{r+\ell, d}(m+1)
}
\end{equation}

\tpoint{Proposition}
{\em  \label{p:final}For any two sequences $(a_1, b_1), (a_2, b_2), (a_3, b_3), \dots$ with $a_i, b_i \in \Z_+$, weakly increasing, and $\lim a_i = \infty$,
  the elements
\begin{equation}
\label{eq:final} \xymatrix{ U_{r+a_1, d}(b_1) \ar[r] &  U_{r+a_2, d}(b_2) \ar[r] &  U_{r+a_3, d}(b_3) \ar[r] &  \cdots}
\end{equation}
  are canonically isomorphic as elements of $\cGw$.}

It is perhaps surprising that Proposition \ref{p:final} only requires that  $\lim a_i = \infty$ and not that $\lim b_i = \infty$. Moreover, we will actually use this fact in the proof of Theorem \ref{thm:bpg} below to take a sequence where $b_i = m$ is constant!

\bpf Such sequences form a partially ordered set where we say $(a_i, b_i) \leq (a_i', b_i')$ if $a_i \leq a_i'$ and $b_i \leq b_i'$ for all $i$. Given any two weakly increasing sequences $(a_i, b_i)$ and $(a_i', b_i')$ with $\lim a_i = \infty$ and $\lim a_i' = \infty$, there exists a third such sequence dominating both, namely $(\max\{a_i, a_i'\}, \{b_i, b_i'\})$. It thus suffices to prove the result when $(a_i, b_i) \leq (a_i', b_i')$. In this case, we define a morphism
\begin{center}
    \begin{tikzcd}
    U_{r+a_1,d}(b_1) \arrow{r} \arrow{d} & U_{r+a_2,d}(b_2) \arrow{r} \arrow{d} & U_{r+a_3,d}(b_3) \arrow{r} \arrow{d} & \cdots \\
        U_{r+a_1',d}(b_1') \arrow{r} & U_{r+a_2',d}(b_2') \arrow{r} & U_{r+a_3',d}(b_3') \arrow{r} & \cdots 
    \end{tikzcd}
\end{center}
where the vertical arrows are
obtained from composing appropriate arrows in the grid \eqref{grid}.
The above discussion shows that 
the  arrows of \eqref{eq:final} are ``increasingly connected'' the further rightward one goes (see \eqref{eq:star}). Consequently, the vertical arrows above become increasingly connected as we move rightward, so they define a morphism in $WHE$.\epf

Based on this (and \S \ref{qccomment}), we make the following definition.

\epoint{Definition}
Define $\Omega^2_{\alg, d}(BGL)$ to be any element of the isomorphism class described in Proposition \ref{p:final} in $\cGw$.  (It is defined uniquely up to unique isomorphism.)\label{d:final}

\tpoint{Theorem (Bott periodicity, in $\cGw$)} \label{thm:bpg} {\em We have an isomorphism in $\cGw$:\label{t:bpnew}
  $$\mu: \Omega^2_{\alg, d}(BGL) \overset \sim \longrightarrow BGL.$$ }

\bpf
Consider (with $m \in \Z^+$ as above) the commutative diagram
\begin{equation} 
\begin{tikzcd}
\tU_{r,d}(m) \arrow{d} \arrow{r} & \tU_{r+1,d}(m) \arrow{d} \arrow{r} & \tU_{r+2,d}(m) \arrow{d} \arrow{r} & \cdots \\
U_{r,d}(m)  \arrow{r} &  U_{r+1,d}(m)   \arrow{r} & U_{r+2,d}(m) \arrow{r} & \cdots
\end{tikzcd}
\end{equation}
which is a morphism in $\cGst$.  The bottom row
is $\Omega^2_{\alg, d}(BGL)$ (Proposition~\ref{p:final}, Definition~\ref{d:final}), and the top row is isomorphic to $BGL \in \cGst$ (\S \ref{s:fur}).  \epf

We record here further consequences proved in \S \ref{appendix:section:comparison}.

\tpoint{Theorem (Bott periodicity in motivic spaces, in 
${\mathrm{HoSpc}}(k)$)} 
{\em Suppose $B=\Spec k$, where $k$ is a perfect field.  Then the isomorphism of Theorem~\ref{thm:bpg}, after applying the functor of Theorem~\ref{t:49}(b), yields an isomorphism in the homotopy category of motivic spaces over $k$, ${\mathrm{HoSpc}}(k)$.    The resulting isomorphism in motivic spectra is the previously-known Bott periodicity for motivic spectra.}

\section{Proof of   ``traditional'' Bott periodicity}

\label{pf:trad}
We now show that this algebraic construction 
indeed generalizes Bott periodicity in topology: this is a
generalization, not merely an analogue.

\point First, we point out that Theorem~\ref{t:bpnew} immediately gives a homotopy equivalence
\begin{equation} \label{eq:pribp}\top(\Omega^2_{\alg, d}(BGL)) \overset \sim \longrightarrow  BU\end{equation}
by taking $B = \C$ and applying the functor $\top: \cGw \rightarrow \HoTop$ (for all $d \in \Z$).

\point  To do this, we quickly explain why there are homotopy equivalences  $$\top(BGL) \overset \sim \longrightarrow BGL \overset \sim \longrightarrow BU.$$  
 Warning for the first arrow:  the first $BGL$ is an object of $\cGw$ (or $\cGst$), and the second $BGL$ is  an object in $\HoTop$.   This first  homotopy equivalence comes from the fact that both are represented by the infinite Grassmannian $G(\infty, \infty + \infty)$ (see \S \ref{d:BGL} for the first; the second is standard). The second homotopy equivalence $BGL \overset \sim \longrightarrow BU$ is standard. 

\point 
It is not hard to show that there is a map
$\top(\Omega^2_{\alg, d}(BGL)) \rightarrow \Omega^2_d(BGL)$.  But it is not clear that this latter map is an isomorphism
(in fact, this will only be clear by the end of this section), and even if that were known, it would not be clear that
the map of \eqref{eq:pribp} is the same as that in (classical) Bott periodicity.

\point So we settle these issues now.  (We note that a different argument in the style of Mann and Milgram --- see \cite{bhmm} and the references and discussion therein --- can be given.  This might be  discussed in \cite{bryanvakil}.)

\tpoint{Theorem}  {\em \label{t:bp2}Take $B= \C$.  The isomorphism $\top(\mu)$ induced by Theorem~\ref{t:bpnew}
agrees with Atiyah's isomorphism  $\al$ in his proof of Bott periodicity.
Translation:  we have the following commuting diagram in $\HoTop$.
\begin{equation}\label{eq:bp2}
            \xymatrix{ \top( \Omega^2_{\alg, d}(BGL)) \ar[r]_-\sim^-{\top(\mu)} \ar[d]^\sim & \top(BGL) \ar[d]^\sim \\
\Omega_d^2 BGL \ar[d]^\sim & BGL \ar[d]^\sim \\
\Omega_d^2( BU) \ar[r]_-\sim^-\alpha & BU            }
\end{equation}}

\bpf We identify our map from spaces of rank $r$ bundles on $\proj^1$ to $BGL$ with Atiyah's
map for each $r$.  We do this by identifying our construction with Mitchell's
construction in the holomorphic category in \cite{mitchell2}.
For each $m$, we have an isomorphism
 $U_{r,d}(m) \overset \sim \longrightarrow U_{r,d+mr}(0)$ defined by $\cE \mapsto \cE(m)$.
 Therefore, it will suffice to identify our maps out of $U_{r,d}(0)$ for each $d$ with Mitchell's construction (and see that these identifications play well with twisting).

 By \eqref{gs1}, for each $r$ and $d$, our map $U_{r,d}(0) \rightarrow BGL$
is given by
$\cE \mapsto \pi_* \cE(-1)$, which agrees with Atiyah's
map (defined using index theory --- this was also observed in
\cite[Rem.~1.6.11]{lupercio}; but see Remark~\ref{r:funandprofit} for
our preferred reinterpretation), interpreted as the map $$j:  \operatorname{Hol}_{d}(\C \proj^1, BU(r))
\rightarrow BU(d)$$  of \cite[p.~802]{cls}.

Consider the following diagram in $\cG$. Each of the stacks in the top three rows are interpreted as elements of $\cG$. (Pedantically: treat them as constant sequences, and verify that they have  spatial incarnations. In the top row, they are algebraic spaces, so they are their own spatial incarnations.  For the second row, the constant sequences of $U$'s admit spatial incarnations by corresponding constant sequences of $T$'s.  In the third row, the spatial incarnations were described in \eqref{eq:BGLnsp}.)  
\begin{equation}\label{eq:aug2}
\xymatrix{T_{r,d}(0)   \ar[d]_{\text{affine bundle}}  \ar@{^(->}[r] & T_{r,d}(1)  \ar[d]^{\text{affine bundle}} \ar@{<->}[r]^-\sim & T_{r, d+r}(0) \ar[d]^{\text{affine bundle}} \\
U_{r,d}(0)   \ar[d]_{\cE \mapsto \pi_* \cE(-1)}  \ar@{^(->}[r] & U_{r,d}(1)  \ar[d]^{\cE \mapsto \pi_* \cE} \ar@{<->}[r]^-\sim & U_{r, d+r}(0) \ar[d]^{\cE \mapsto \pi_* \cE(-1)} \\
BGL(d) \ar[r] \ar[drr] & BGL(d+r) \ar[dr]  \ar@{<->}[r]^-\sim  & BGL(d+r) \ar[d] \\
& & BGL}
\end{equation}
This is clearly a commuting diagram {\em except} for the square
\begin{equation}\label{eq:oct14z}
\xymatrix{U_{r,d}(0)   \ar[d]_{\cE \mapsto \pi_* \cE(-1)}  \ar@{^(->}[r] & U_{r,d}(1)  \ar[d]^{\cE \mapsto \pi_* \cE} \\BGL(d) \ar[r] & BGL(d+r)
}\end{equation}
which does {\em not} commute. The map obtained by the composition $U_{r,d}(0) \rightarrow BGL(d) \rightarrow BGL(d+r)$ corresponds to $\cE \mapsto \pi_* \cE(-1) \oplus \oh^{\oplus r}$ (interpreted in terms of the functor of points), while the map obtained by the composition $U_{r,d}(0) \rightarrow U_{r,d}(1) \rightarrow BGL(d+r)$ corresponds to $\cE \mapsto \pi_* \cE$. These two bundles are not isomorphic! (Typically, $\pi_*\cE$ may be a non-trivial extension of $\oh^{\oplus r}$ by $\pi_* \cE(-1)$.)
Nonetheless, the larger rectangle
\begin{equation}\label{eq:oct14}
\xymatrix{T_{r,d}(0)   \ar[d]  \ar@{^(->}[r] & T_{r,d}(1)  \ar[d]
\\BGL(d) \ar[r] & BGL(d+r)
}
\end{equation}
in \eqref{eq:aug2} {\em does} commute, as we now explain.  (``Adding in the data of a splitting of the exact sequence resolves the lack of commutativity of \eqref{eq:oct14z}.'')
As before, the composition $T_{r,d}(0) \to BGL(d) \to BGL(d+r)$ corresponds to $(\cE, s_1, \ldots, s_r) \mapsto \pi_* \cE(-1) \oplus \oh^{\oplus r}$ while the composition 
$T_{r,d}(0) \to T_{r,d}(1) \to BGL(d+r)$
 corresponds to $(\cE, s_1, \ldots, s_r) \mapsto \pi_* \cE$. The difference now is that since $(\cE, s_1, \ldots, s_r)$ is an object of $T_{r,d}(0)$, the sections $s_1, \ldots, s_r$ define a splitting of the Str\o mme sequence: $\pi_*\cE \cong \pi_* \cE(-1) \oplus \oh^{\oplus r}$ (as in the proof of Lemma \ref{l:xgla}). 

Thus applying $\top$
and to  \eqref{eq:aug2} and \eqref{eq:oct14}, we have that the diagram
$$\xymatrix{\top(T_{r,d}(0)) \ar[dr] \ar[rr] & & \top(T_{r,d+r}(0)) \ar[dl] \\
  & BGL
}$$
commutes in $\HoTop$, and hence 
\begin{equation}\label{eq:aug22b}\xymatrix{\top(U_{r,d}(0)) \ar[dr] \ar[rr] & & \top(U_{r,d+r}(0)) \ar[dl] \\
  & BGL
}
\end{equation}
does as well.

We now explain why this is precisely  the diagram
\begin{equation} \label{eq:top}
\begin{tikzcd}
F_{r,d} \arrow{rr} \arrow{dr} &&F_{r,d+r} \arrow{dl} \\
&BU
\end{tikzcd}
\end{equation}
of Mitchell \cite[Thm.~A]{mitchell}, see  \cite[p.~802]{cls}.
(Caution:  \cite{cls} has reversed
 Mitchell's subscripts \cite{mitchell}; we have kept Mitchell's order, which agrees
with ours.)  Now $U_{r,d}(0)^{an} \cong \operatorname{Hol}_d(\proj^1, BU(r))$
by Serre's GAGA theorem:  informally, all holomorphic rank $r$ bundles of degree
$d$ are algebraic.
By \cite[Thm.~2]{cls}, there is a natural homotopy equivalence
$\operatorname{Hol}_d(\proj^1, BU(r)) \cong F_{r,d}$.  Under these identifications,
the maps of \eqref{eq:aug22b} are sent to the maps of \eqref{eq:top} (up to homotopy).  (The maps to
$BGL$ and $BU$ are, up to homotopy, just pushforward of the bundles on the families.)

The commuting diagram \eqref{eq:top} in $\HoTop$ (for all $r$ and $d$)  determines the map in the
topological incarnation of Bott periodicity  (this is part of the
content of \cite[Thm.~A]{mitchell}).
Thus we obtain  the desired commutative diagram
\eqref{eq:bp2}.
\epf

\epoint{Remark} We conclude that  Atiyah's map, originally defined using index theory, might be
profitably algebraically described as $ \Maps_d^\bullet(\proj^1,
BGL(r)) \rightarrow BGL$ given by $\cE \mapsto \pi_! \cE(-1)$ (not
just $\cE \mapsto \pi_* \cE(-1)$ on $U_{r,d}(0)$).\label{r:funandprofit}

\section{New topological consequences in finite rank $r$  (over $\C$)} \label{s:four}

This section deals with the topology of the moduli stack $\cX_{r,d}$.  In this section as in the previous one, $B=\C$.  We begin by noting that $\cX_{r,d}$ is the colimit (the infinite increasing union of open substacks) of the $U_{r,d}(m)$, so the only possible meaning of the topological type of $\cX_{r,d}$ in $\HoTop$ is as $\top(\Omega^2_{\alg, d}(BGL(r))$ (which is similarly described in terms of the $U_{r,d}(m)$).

  \tpoint{Theorem} \label{hequiv} {\em
The morphism $\lambda_{r,d}: \top(\Omega^2_{\alg, d}(BGL(r)) )
        \rightarrow \Omega_d^2 BGL(r)$ is a homotopy equivalence.
    }

    \bpf   Our strategy, as in  the proof of Theorem~\ref{t:bp2},
    is to identify our
    construction with the  geometric construction of
    \cite{mitchell2}.  The $F_{r,d}$ of \cite{mitchell2} give a
    filtration of the space $\Omega SU(r)$ \cite[Thm.~A]{mitchell2}.    This filtration corresponds
    precisely to our filtration of the space of bundles as $U_{r,d}(0) \subset U_{r,d}(1) \subset U_{r,d}(2) \subset \cdots$ (after
    applying $\top$, as described in the proof of Theorem \ref{t:bp2}).
  \epf

It follows that the cohomology rings of these ``loop spaces'' (the target
of $\lambda_{r,d}$) agree with the
cohomology rings of the respective moduli spaces (the sources of
$\lambda_{r,d}$).

\tpoint{Theorem} \label{algcoh} {\em
        The map $\lambda_{r,d}$ identifies the cohomology ring  of
        $\Omega^2_d BGL(r)$ (with integral coefficients!) with the Chow ring of $\cX_{r,d}$.
    }

\begin{proof}
Theorem \ref{hequiv} gives an isomorphism
$H^*(\Omega^2_dBGL(r), \Z) \cong H^*(\cX_{r,d}, \Z)$, 
 so we will establish that the Chow and cohomology rings of $\cX_{r,d}$ agree.
We follow a similar strategy to
\cite[\S 5]{hannah}, namely, we consider the stratification by splitting loci, which we show are each modeled by a space admitting a stratification by affine spaces.

Given a splitting type $\vec{e} = (e_1, \ldots, e_r)$,
let $\Sigma_{\vec{e}} \subset \cY_{r,d}$ be the universal splitting locus for splitting type $\vec{e}$ (a locally closed substack). 
Now $\Sigma_{\vec{e}}$ consists of a single (``stacky'') point;
$\Sigma_{\vec{e}} \cong B (G_{\vec{e}})$ where $G_{\vec{e}}$ is the
automorphism group (scheme) of the bundle $\oh(e_1) \oplus \cdots \oplus
\oh(e_r)$.  The group $G_{\vec{e}}$ is identified with block
upper-triangular matrices, where the blocks correspond to the
multiplicities $n_j$ of the parts of $\vec{e}$ (so $\sum n_j = r$).
The diagonal blocks correspond to $GL(n_j)$, and the remaining
nonzero entries are affine spaces. It follows that
the natural map $BGL(n_1) \times \cdots \times BGL(n_s) \to BG_{\vec{e}} \cong \Sigma_{\vec{e}}$ is an affine bundle.

Recall that $\psi_{r,d}: \cX_{r,d} \to \cY_{r,d}$ is a principal $GL(r)$-bundle.
Let $\widetilde{\Sigma}_{\vec{e}} = \psi_{r,d}^{-1}(\Sigma_{\vec{e}}) \subset \cX_{r,d}$, so we have a (2-)fiber diagram of stacks
\begin{center}
    \begin{tikzcd}
    GL(r)/(GL(n_1) \times \cdots \times GL(n_s)) \arrow{d} \arrow{r} & BGL(n_1) \times \cdots \times BGL(n_s) \arrow{d} \\
    \widetilde{\Sigma}_{\vec{e}} \arrow{r} & \Sigma_{\vec{e}}
    \end{tikzcd}
\end{center}
where the vertical arrows are affine bundles and the horizontal maps are principal $GL(r)$-bundles. Meanwhile, $GL(r)/(GL(n_1) \times \cdots \times GL(n_s))$ is an affine bundle over the partial flag variety $Fl := Fl(n_1, n_1 + n_2, \ldots, n_1 + \ldots + n_{s-1};r)$. It follows that \begin{equation} \label{isos}
A^*(\widetilde{\Sigma}_{\vec{e}}) \cong A^*(Fl) \cong H^*(Fl, \Z) \cong H^*(\widetilde{\Sigma}_{\vec{e}},\Z).
\end{equation}

Let $U \subset \cX_{r,d}$ be an open substack which is a finite union of splitting loci so that $\widetilde{\Sigma}_{\vec{e}}$ is closed inside $U \cup \widetilde{\Sigma}_{\vec{e}}$.
 Let $u(\vec{e}) = \codim (\widetilde{\Sigma}_{\vec{e}} \subset \cX_{r,d})$. By the last isomorphism of \eqref{isos}, we see that all odd cohomology of $\widetilde{\Sigma}_{\vec{e}}$ vanishes.
Thus, arguing as in \cite[Lem.~5.3]{hannah}, we have exact sequences
\begin{equation} \label{coheq}
\xymatrix{0 \ar[r] & H^{2(i - u(\vec{e}))}(\widetilde{\Sigma}_{\vec{e}}, \Z) \ar[r] &  H^{2i}(U \cup \widetilde{\Sigma}_{\vec{e}},\Z) \ar[r] & H^{2i}(U,\Z) \ar[r] &  0.}
\end{equation}
Meanwhile, the first higher Chow groups with mod $p$ coefficients $CH^*(\widetilde{\Sigma}_{\vec{e}},1, \Z/p\Z)$ also agree with those of the partial flag variety $Fl$.
Working over $\C$, the latter must vanish because $Fl$ admits a stratification into affine spaces (see \cite[p. 3]{hannah}). Thus, noting that the Chow groups $A^i(\widetilde{\Sigma}_{\vec{e}}) = A^i(Fl)$ are finitely-generated free $\Z$-modules, and
arguing as in \cite[Lem.~ 5.2]{hannah}, we have exact sequences
\begin{equation} \label{choweq}
\xymatrix{0 \ar[r] &  A^{i - u(\vec{e})}(\widetilde{\Sigma}_{\vec{e}}) \ar[r] & A^{i}(U \cup \widetilde{\Sigma}_{\vec{e}}) \ar[r] & A^{i}(U) \ar[r] & 0.}
\end{equation}
Now consider the natural map from \eqref{choweq} to \eqref{coheq}. By \eqref{isos} and induction on the number of strata in $U$, the outer two maps are isomorphisms. By the $5$-lemma, the middle map is too. As there are finitely many $\vec{e}$ with $u(\vec{e}) \leq i$, we see $A^*(\cX_{r,d}) \to H^{*}(\cX_{r,d})$ is an isomorphism in arbitrarily high codimension, hence an isomorphism.
\end{proof}

\epoint{Remark}  We note two consequences of the proof.  

First, the cohomology/Chow ring is torsion-free.  This allows a description of $A^*(\cX_{r,d})$ similar to \cite[Thm.~1.2]{hannah}. As $\cX_{r,d} \to \cY_{r,d}$ is a $GL(r)$-bundle, 
\[ A^*(\cX_{r,d}) = A^*(\cY_{r,d})/\langle a_1, \ldots, a_r \rangle \]
where $a_i$ are the Chern classes of the associated rank $r$ vector bundle, which in this case is
$\sigma_\infty^* \cE$.
Following the notation of \cite[Thm.~1.2]{hannah}, we may thus realize $A^*(\cX_{r,d})$ as the subring of $A^*_{\Q}(\cX_{r,d}) = \Q[a_2', \ldots, a_r']$ (here $a_i' \in A^{i-1}(\cX_{r,d})$) generated by the coefficients of $t^j$ in the power series expansion of $$\exp \left(\int  a_2' + a_3't + \cdots + a_r' t^{r-2}  \; dt \right).$$ As $\cX_{r,d}$ is homotopic to the affine Grassmannian, this ring is also the cohomology of the affine Grassmannian, which
was first found by
Bott \cite{bott1}; 
for alternative descriptions of this ring see also \cite{yz}, or
\cite[\S 6.1]{lam}.

Second, by \eqref{isos}, the cohomology/Chow ring of each stratum $\widetilde{\Sigma}_{\vec{e}} \subset \cX_{r,d}$ is nonzero only in finitely many degrees. (This is in contrast with the cohomology/Chow rings of the splitting loci $\Sigma_{\vec{e}} \subset \cY_{r,d}$, which are non-zero in all degrees.)

\appendix

\section{Compatibility with motivic homotopy theory}

\begin{center}
\chapterauthor{Benjamin Church}
\end{center}

\renewcommand{\P}{\mathbb{P}}
\newcommand{\U}{\mathfrak{U}}

\newcommand{\sq}{\mathrm{sq}}
\newcommand{\st}{\mathrm{st}}
\newcommand{\Nis}{\mathrm{Nis}}
\newcommand{\Shv}{\mathrm{Shv}}
\newcommand{\PSh}{\mathrm{PSh}}
\newcommand{\Sm}{\mathrm{Sm}}
\newcommand{\Spc}{\mathrm{Spc}}

\newcommand{\op}{\mathrm{op}}
\newcommand{\Fun}[2]{\mathrm{Fun}\left( #1, #2 \right)}

\newcommand{\SH}{\mathrm{SH}}
\newcommand{\KGL}{\mathrm{KGL}}
\newcommand{\GL}{\mathrm{GL}}
\newcommand{\BGL}{\mathrm{BGL}}
\newcommand{\Gr}{\mathrm{Gr}}

\newcommand{\mot}{\mathrm{mot}}
\newcommand{\Aconn}[1]{\A^1\text{-}#1\text{-connected}}
\renewcommand{\cV}{\mathcal{V}}

\newcommand{\Hen}{\mathrm{Hen}}

\newcommand{\fppf}{\mathrm{fppf}}

\newcommand{\X}{\mathcal{X}}
\newcommand{\V}{\mathbb{V}}
\newcommand{\cU}{\mathrm{U}}

\renewcommand{\sp}{\mathrm{sp}}
\newcommand{\Ho}{\mathrm{Ho}}
\newcommand{\Th}{\mathrm{Th}}
\newcommand{\Y}{\mathcal{Y}}

\renewcommand{\G}{\color{red} CHANGE}
\newcommand{\ot}{\otimes}
\newcommand{\Sing}[1]{\mathrm{Sing}(#1)}
\newcommand{\embed}{\hookrightarrow}
\newcommand{\onto}{\twoheadrightarrow}
\newcommand{\id}{\mathrm{id}}
\newcommand{\iso}{\TikzBiArrow[->]}
\newcommand{\cO}{\mathcal{O}}

\newcommand{\sm}{\setminus}
\newcommand{\cN}{\mathcal{N}}

\label{appendix}

\bpoint{Introduction to the appendix}  \label{appendixintroduction}The work of Voevodsky, Morel, and many others have offered a convincing answer to the question of where to situate homotopy-theoretic manipulations of algebro-geometric objects. The categories of motivic spaces and motivic spectra inherit  many of the desirable $\infty$-categorical properties of spaces and spectra while allowing for arguments over arbitrary characteristic and maintaining motivic information such as Hodge structures that are lost when passing to underlying topological manifolds. The construction of a Bott map using purely algebro-geometric constructions on classifying stacks has exactly the flavor an argument motivic spaces ought to support. However, it is not immediately clear that certain construction recover homotopy-theoretic notions such as loop spaces and highly connective maps. Larson and Vakil give a list of desiderata necessary for a category to support their construction of the Bott map  (\S \ref{sec:desiderata}). They produce some sort of ``free'' categories $\cGst$ and $\cGsp$ initial for those satisfying the desiderata. We show that motivic spaces can support Larson and Vakil's construction by exhibiting a functor from their free homotopy category to the homotopy category of motivic spaces that lifts to a functor of $\infty$-categories. 
\par 
One difficulty is that $\A^1$-connectivity is much more stringent than topological connectivity, so care is necessary in dealing with arguments relying on highly connected maps. Perhaps surprisingly, the ``obvious'' functor does not work -- because inclusions of an open subset with large complementary codimension
may not be highly $\A^1$-connected. However, this can be remedied by a partial $S^1$-stabilization (see \S 3.1 and Theorem~\ref{thm:top_bottom_comparison}) in order for Larson and Vakil's class of weak homotopy equivalences to actually be motivic equivalences. This partially stabilized functor agrees with the naive functor after fully stabilizing, {\em i.e.}, mapping to the category of $S^1$-motivic spectra. This compatibility establishes a new result native to motivic homotopy theory, namely an analogue of the main result of \cite{cls} for algebraic $K$-theory that, put imprecisely, the space of algebraic maps from $\P^1$ to $\BGL$ converges to the motivic $\P^1$-loop space of $\BGL$ (see \S\ref{appendix:section:comparison}). The arguments given here closely follow prior developments in the motivic homotopy theory literature. 

\bpoint{Primer on Motivic Homotopy Theory}

Let $B$ be an excellent Noetherian scheme of finite Krull dimension. The goal of this section will be to review the construction of the $\infty$-categories $\Spc(B)$ and $\SH(B)$ of motivic spaces over $B$ and motivic spectra over $B$ respectively. 

\epoint{Definitions}
Let \underline{$\Sm_B$} be the category of schemes smooth over $B$. Define the functor category \underline{$\PSh(\Sm_B) := \Fun{\Sm_B^\op}{\Spc}$} to be the $\infty$-category of presheaves valued in spaces.

The full subcategory \underline{$\Shv_{\Nis}(\Sm_B) \subset \PSh(\Sm_B)$} consists of those presheaves $F$ which are \underline{Nisnevich sheaves},  meaning that for any Nisnevich cover $\cU = \{ U_i \to X \}_i$, the map
\[ F(X) \to \lim_{\Prism} F(\cU) \]
is an equivalence. Above, the homotopy limit is taken over the cosimplicial diagram $\Prism \to \Spc$ given by $F$ applied to the \v{C}ech nerve of $U$. 

The full subcategory $\PSh_{\A^1}(\Sm_B) \subset \PSh(\Sm_B)$ consists of \underline{$\A^1$-invariant presheaves}, {\em i.e.} the presheaves for which the morphism 
\[ F(X \times \A^1) \to F(X) \]
is an equivalence for all $X \in \Sm_B$.

The \underline{category of motivic spaces} is defined as the full subcategory consisting of presheaves that are both \underline{Nisnevich-local} ({\em i.e.}, sheaves) and $\A^1$-invariant:
\[ \Spc(B) := \Shv_{\Nis}(\Sm_B) \cap \PSh_{\A^1}(\Sm_B). \]
Call a presheaf \underline{motivic} if it lives in this subcategory. 
 
These subcategories have the additional distinguishing property of being \underline{reflective}, meaning the inclusions into $\PSh(\Sm_B)$ have left adjoints called ``localization'' written $L_{\Nis}, L_{\A^1}$ respectively and $L_{\mot}$ for the localization into $\Spc(B)$. Since these categories are the full subcategories of $S$-local objects for appropriate choices of $S$, the existence of a left adjoint follows from \cite[5.5.4.15]{Lur09}. For $\Shv_{\Nis}(\Sm_B) \embed \PSh(\Sm_B)$ the left adjoint \underline{$L_{\Nis} : \PSh(\Sm_B) \to \Shv_{\Nis}$} is exactly sheafification (in the homotopy-coherent sense when valued in honest spaces).  

A $1$-morphism $f : F \to G$ in $\PSh(\Sm_B)$ is a \underline{Nisnevich-local equivalence} if $L_{\Nis} f$ is an equivalence.

Let $X$ be a scheme and $x \in X$ a point. Let \underline{$\Hen_{X,x}$} be the category of pairs $g : (Y, y) \to (X,x)$ where $g : Y \to X$ is \etale and $g(y) = x$ such that $g^{\#} : \kappa(x) \to \kappa(y)$ is an isomorphism. For any presheaf $F \in \PSh(\Sm_B)$ the \underline{stalk} at $x$ is defined to be
\[ F_x := \colim\limits_{U \in \Hen_{X,x}} F(U). \] 
A result, originally due to Jardine, says, in the following form, that Nisnevich local equivalences can be checked on stalks:

\tpoint{Theorem \cite[Lemma~1.11, p.100]{mv}} {\em 
A morphism $f$ in $\PSh(\Sm_B)$ is a Nisnevich local equivalence if and only if for each $X \in \Sm_B$ and $x \in X$ the morphism $F_x \to G_x$
is an equivalence.}

For the inclusion $\PSh_{\A^1}(\Sm_B) \embed \PSh(\Sm_B)$, the left adjoint \underline{$L_{\A^1} : \PSh(\Sm_B) \to \PSh_{\A^1}(\Sm_B)$} has an explicit model given by the ``singular chains'' construction. Denote by \underline{$\Delta^n \in \Sm_B$} the \underline{algebraic $n$-simplex}
\[ \Delta^n := \Spec{\Z[t_0,\dots,t_n]/(t_1 + \cdots + t_n - 1)} \]
which, along with the obvious maps, gives a cosimplicial scheme $\Delta^\bullet \in \Fun{\Prism}{\Sm}$. Then as a model of $L_{\A^1}$ we can define
\[ (L_{\A^1} F)(X) = \colim\limits_{\Prism^\op} F(X \times \Delta^\bullet) \]
where the homotopy colimit is taken over this simplicial space. Finally, the motivic localization adjoint
\[ L_{\mot} : \PSh(\Sm_B) \to \Spc(B) \]
is computed as the homotopy colimit of alternatively iterating $L_{\Nis}$ and $L_{\A^1}$.

\epoint{Homotopy Groups and Connectivity}

Since there are inconsistent usages in the literature, we recall what it means for a space or a morphism to be $n$-connected. Our conventions follow conventions common to motivic homotopy theory (e.g. \cite{HJNTY}) but which disagree with those commonly (and historically) used in homotopy theory. Note that sometimes the parallel terminology $n$-\textit{connective} (arising from the theory of spectra), meaning $(n-1)$-connected, is used instead in this context.

\epoint{Definition}
We say a \underline{space $X$ is  $n$-connected} if $\pi_i(X, x) = 0$ for $i \le n$ and all choices of basepoint $x \in X$.

We say a \underline{morphism $f : X \to Y$ is $n$-connected} if for each basepoint $x \in X$,
\begin{enumerate}
\item $f_* : \pi_i(X, x) \to \pi_i(Y, f(x))$ is an isomorphism for $i \le n$
\item $f_* : \pi_i(X, x) \to \pi_i(Y, f(x))$ is an epimorphism for $i = n+1$.
\end{enumerate}

\epoint{Remark}
$X \to *$ is $n$-connected if  and only if $X$ is $n$-connected. Likewise, if $X$ is nonempty then $X$ is $(n+1)$-connected if and only if $* \to X$ is $n$-connected. 

\epoint{Remark}
From the Puppe sequence, a morphism $f : X \to Y$ is $n$-connected if and only if the homotopy fiber $F_f$ is $n$-connected.

\epoint{Remark}
Note how this convention plays with connectivity of a pair $(X, A)$. We say that $(X,A)$ is $n$-connected if the relative homotopy groups (see \cite[p.343]{HatcherAG}) vanish, $\pi_k(X,A) = 0$, for $k \le n$. This is equivalent, by the long exact sequence, to the inclusion map $A \to X$ being $(n-1)$-connected. Our convention extends to a map of pairs $(X,A) \to (X', A')$ verbatim using the relative homotopy groups. Precisely, $(X,A) \to (X', A')$ is \underline{$n$-connected} if $\pi_k(X,A) \to \pi_k(X', A')$ is an isomorphism for $k \le n$ and an epimorphism for $k = n+1$.

\epoint{Definition}
In motivic homotopy theory, the objects are not spaces but rather presheaves valued in spaces. Hence, the natural objects replacing the homotopy groups are sheaves of  groups.

Let $F \in \PSh(\Sm_B)$ be a pointed presheaf. We define the homotopy groups $\pi_i(F)$ as the Nisnevich sheafification of the presheaf,
\[ U \mapsto \pi_i(F(U)). \]

We say \underline{a presheaf $F \in \PSh(\Sm_B)$ is $n$-connected} if $\pi_i(F) = 0$ for $i \le n$.

We say \underline{a morphism $f : F \to G$ in $\PSh(\Sm_B)$ is $n$-connected} if
\begin{enumerate}
\item $f_* : \pi_i(F) \to \pi_i(G)$ is an isomorphism for $i \le n$
\item $f_* : \pi_i(F) \to \pi_i(G)$ is an epimorphism for $i = n+1$.
\end{enumerate} 

A presheaf \underline{$F \in \PSh(\Sm_B)$ is $\Aconn{n}$} if $L_{\mot}(F)$ is $n$-connected. A morphism \underline{$f : F \to G$ in $\PSh(\Sm_B)$ is $\Aconn{n}$} if $L_{\mot}(f)$ is $n$-connected. 

\epoint{Remark}
We can furthermore define,
\[ \pi_i^{\A^1}(\X) := \pi_i(L_{\mot} \X) \]
and then $\A^1$-connectivity of a space or of a map corresponds to the expected definition in terms of $\pi_i^{\A^1}$. The importance of this definition is that the category of motivic spaces is \textit{hypercomplete}, meaning it satisfies the analog of the Whitehead theorem.

\tpoint{Theorem~\cite[Proposition 2.14, p.110]{mv}} {\em 
Let $f : X \to Y$ be a morphism in $\Spc(B)$ of pointed $\A^1$-connected spaces such that $f_* : \pi_i(X) \to \pi_i(Y)$ is an isomorphism for all $i$. Then $f$ is an equivalence.}

\epoint{Definition}
Recall that we define the \underline{loop space} and \underline{suspension} of a space or presheaf as the following pullback and pushout respectively:
\begin{center}
\begin{tikzcd}
\Omega X \pullback \arrow[d] \arrow[r] & * \arrow[d]
\\
* \arrow[r] & X
\end{tikzcd}
\quad \quad 
\begin{tikzcd}
X \arrow[r] \arrow[d] & * \arrow[d]
\\
* \arrow[r] & \pushout \Sigma X.
\end{tikzcd}
\end{center}
We may write $\Omega_{\PSh}$ and $\Sigma_{\PSh}$ to emphasize that the operations are computed in the ambient category of presheaves.

\epoint{Remark}
Since $L_{\mot}$ preserves homotopy colimits, $L_{\mot}(\Sigma_{\PSh} \X) \cong \Sigma_{\mot} L_{\mot}(\X)$. However, this $\Sigma_{\mot}$ is suspension computed in the category of motivic spaces. Even if $\X$ is motivic, $\Sigma_{\PSh} \X$, computed by viewing $\X$ as a mere presheaf, need not be motivic. From now on, if $\X$ is a motivic space, we drop the decoration and write $\Sigma \X := \Sigma_{\mot}\X \cong L_{\mot} \Sigma_{\PSh} X$. 

On the other hand, by the existence of a left adjoint, the inclusion $\Spc(B) \embed \PSh(\Sm_B)$ preserves homotopy limits meaning $\Omega_{\mot}$, the loop space internal to $\Spc(B)$, agrees with $\Omega_{\PSh}$ on (the inclusion of) motivic spaces. Hence, if $\X$ is motivic, $\Omega \X := \Omega_{\PSh} \X$ is also motivic. However, in general, $L_{\mot} \Omega_{\PSh} \X \not\cong \Omega_{\PSh} L_{\mot} \X$ if $\X$ is not motivic. This means although $\pi_i(\Omega_{\PSh} \X) = \pi_{i+1}(\X)$ is immediate from the definitions, $\pi_i^{\A^1}(\Omega_{\PSh} \X) \neq \pi_{i+1}^{\A^1}(\X)$ in general. However, if $\X$ is motivic then $\Omega \X$ is also motivic and hence
\[ \pi_i^{\A^1}(\Omega \X) = \pi_{i+1}^{\A^1}(\X) \]
from the result on the unlocalized $\pi_i$. 

The above remark has demonstrated the following.

\tpoint{Lemma}{\em 
If $\X \in \PSh(\Sm_B)$ is motivic, then $\Omega \X$ is motivic and $\pi_i^{\A^1}(\Omega \X) = \pi_{i+1}^{\A^1}(\X)$.}

Finally, it will be useful to know that connectivity can be computed on the stalks of a morphisms of presheaves. 

\tpoint{Lemma} {\em  \label{lem:homotopy_groups_stalks}
Let $F \in \PSh(\Sm_B)$. For $x \in X \in \Sm_B$ the Nisnevich stalk of $\pi_i(F)$ at $x$ is $\pi_i(F_x)$.}

\begin{proof}
We can commute the homotopy colimit with $\pi_i$,
\[ \pi_i(F)_x = \colim\limits_{Y \in \Hen_{X,x}} \pi_i(F(Y)) = \pi_i(\colim\limits_{Y \in \Hen_{X,x}} F(Y)) = \pi_i(F_x) \] 
because CW-complexes are compact objects of $\Spc_*$ so, in particular, the $n$-sphere $S^n$ is a compact object.
\end{proof}

\tpoint{Corollary} {\em \label{cor:connectivity_stalks}
To check if a (pointed) object $F$ or morphism $f : F \to G$ is $n$-connected, it suffices to check on Nisnevich stalks. 
}

\begin{proof}
Indeed, the morphism of sheaves $\pi_i(F) \to \pi_i(G)$ (with $G = *$ for connectivity of the space $F$) is an isomorphism if and only if it is an isomorphism of stalks.  Lemma~\ref{lem:homotopy_groups_stalks} shows that this is the same as checking the corresponding map of spaces $F_x \to G_x$ is $n$-connected. 
\end{proof}

\epoint{Morel's Connectivity Theorem}
\renewcommand{\Y}{\mathcal{Y}}

As in \cite{Mor12}, we let $B = \Spec{k}$ where $k$ is a perfect field throughout this section.

We have noted that if $\X$ is already motivic then its $\A^1$-homotopy sheaves agree with the homotopy sheaves computed viewing $\X$ as a mere presheaf. However, if we start with an arbitrary presheaf $\X$ there is no reason its homotopy sheaves should agree with those of $L_{\mot} \X$ (which are, by definition, $\pi_i^{\A^1}(\X)$). Morel's $\A^1$-Connectivity Theorem tells us that, although nonzero homotopy sheaves may not be preserved by motivic localization, connectivity is preserved.

\tpoint{Morel's $\A^1$-Connectivity Theorem \cite[Theorem 6.38]{Mor12}} \label{t:connectivity} {\em
Let $\X \in \PSh(\Sm_k)_*$ be a pointed presheaf and $n \ge 0$ an integer. If $\X$ is $n$-connected then it is \underline{$\Aconn{n}$}, meaning for all $i \le n$,
\[ \pi_i^{\A^1}(\X) = \pi_i(L_{\mot}(\X)) = 0. \]}

In the corresponding portion of \cite{Mor12}, this result is stated for the $L_{\A^1}$-localization rather than $L_{\mot}$. However, since $L_{\Nis}$ preserves homotopy limits, colimits, and connectivity (since it preserves stalks), it suffices to check that $L_{\A^1}$ preserves connectivity. 

\epoint{Purity}

In algebraic geometry and motivic homotopy theory, we do not have access to tubular neighborhoods. The replacement of this is the purity theorem (Theorem~\ref{t:purity}), which relates the cofiber of an open immersion $j : U \embed X$ with smooth complement to the normal bundle of its complement. 

\epoint{Remark}
In the sequel, 
given a morphism of presheaves $Y \rightarrow X$, 
we will often write $X/Y$ for $\cofib(Y \to X)$. Also, a scheme $X$ and the corresponding representable presheaf on $\Sm_B$ will be freely identified. Notice that even if $X$ and $Y$ are not pointed ({\em i.e.} do not come with a choice of base point), the cofiber $X / Y$ is canonically pointed at the image of $Y$. Therefore, we are free to write $\Sigma (X/Y)$ without choosing basepoints of $X$ or $Y$.

\epoint{Definition}
The \underline{Thom space} of a vector bundle $\E$ over $X$ is
\[ \Th(\E) = \V(\E)  / (\V(\E) \sm \iota(X)) \cong_{\mot}  \P(\E \oplus \cO_{X}) / \P(\E) \]
where $\V(\E)$ is the total space, $\iota : X \to \V(\E)$ is the zero section, and $\P(\E)$ is the projectivization.

\tpoint{Purity Theorem \cite[\S3, Theorem~2.23]{mv}} {\em 
\label{t:purity}Let $Z \subset X$ be closed subscheme with $X$ and $Z$ smooth $B$-schemes.  Then, denoting by $\cN_{Z|X}$ the normal bundle, there is an equivalence
\[ X / (X \sm Z) \cong_{\mot} \Th( \cN_{Z | X} ). \]}

\epoint{Example}
Let $X$ be a smooth $k$-scheme. If $Z = \{ p \}$ is a $k$-valued point then 
\[ X / (X \sm \{ p \}) \cong_{\mot} \A^d_k / (\A^d_k \sm \{ 0 \}) \]
where $d = \dim_p{X}$ is the local dimension.

\epoint{Motivic Spectra}

Let $\Sm_B$ the category of smooth schemes over a scheme $B$. Let $\Spc(B)_*$ be the $\infty$-category of pointed motivic spaces over $B$. The \underline{category of motivic spectra} is the localization $\SH(B) := \Spc(B)_*[(\P^1)^{-1}]$. Denote by
\[ \Sigma^\infty : \Spc(B)_* \to \SH(B) \]
the canonical functor.  It admits a right adjoint $\Omega^{\infty}$.

\tpoint{Theorem \cite[\S3.3]{bachmann}} {\em 
$\SH(B)$ is a stable $\infty$-category.}

The stabilization functor makes $\Sigma$ and $\Omega$ inverse to each other providing the \textit{shift} functors.   

\bpoint{Larson and Vakil's Categories of Spaces}

\label{s:lvcs}First we recall the categories of spaces defined in the main text. Recall (\S \ref{s:startingpoint}) that $\SmArt$ is the $2$-category of Artin stacks that are smooth, irreducible, have affine diagonal, and are finite type over a fixed excellent base ring $B$. The categories
 $\cGsp$ and $\cGst$ were defined in \S \ref{d:gsp}.
 There is an obvious functor $\cGsp \to \cGst$. Note that this functor is full but not faithful. 

\epoint{Incarnation in Motivic Spaces}
There is an evident functor
\[ \SmArt \to \PSh(\Sm_B) \]
sending an Artin stack to its functor of points restricted to schemes smooth over $B$. Explicitly, a stack is a sheaf of $1$-groupoids so it can be viewed as a $1$-truncated sheaf of spaces hence an object in the $\infty$-category $\PSh(\Sm_B)$. One must check that the $2$-categorical structure on $\SmArt$ is compatible with the higher categorical structure on $\PSh(\Sm_B)$. See \cite[Rem.~2.1]{CD23} for additional details. Composing with the functor $L_{\mot} : \PSh(\Sm_B) \to \Spc(B)$ adjoint to the inclusion of the reflexive subcategory of motivic spaces, we get a motivic realization functor
\[ M_{\SmArt} : \SmArt \to \Spc(B). \]  
This gives an evident way to realize the motivic space associated to a sequence
\[ M_{\cGst} : \cGst \to \Spc(B) \quad \quad M_{\cGst} : X_\bullet \mapsto \colim\limits_n M_{\SmArt} X_n \]
where we take the homotopy colimit over the diagram defined by $M(X_\bullet)$. The same definition produces two functors $M_{\cGsp}^{t}, M_{\cGsp}^{b} : \cGsp \to \Spc(B)$ ($t$ for top, $b$ for bottom) defined by 
\[ M_{\cGsp}^t : (Y_\bullet \to X_\bullet) \mapsto \colim\limits_n M_{\SmArt} Y_\bullet \quad \quad M_{\cGsp}^b : (Y_\bullet \to X_\bullet) \mapsto \colim\limits_n M_{\SmArt} X_\bullet. \]
The diagrams $(Y_\bullet \to X_\bullet)$ induce an obvious $1$-morphism $M^t_{\cGsp} \to M^b_{\cGsp}$ in $\Fun{\cGsp}{\Spc(B)}$. Furthermore, $M_{\cGsp}^b$ obviously agrees with the composition
\[ \cGsp \to \cGst \xrightarrow{M_{\cGst}} \Spc(B). \]
We will use -- in the subsequent discussion -- the evident pointed variants of these realization functors where a point of $X_\bullet \in \cGst$ (resp.\ $(Y_\bullet \to X_\bullet) \in G_{\sp}$) is a morphism $* \to X_\bullet$ (resp.\ $(* \to *) \to (Y_\bullet \to X_\bullet)$) where $*$ (resp.\ $(* \to *$)) is the constant sequence $*_n = B$ (resp.\ $B \to B$). 
\par
The difficulty arises when analyzing the comparison map $M_{\cGsp}^t \to M_{\cGsp}^b$. We want this to be an equivalence. Indeed, since any realization of our goal to ``interpret the construction'' inside $\Spc(B)$ will necessarily construct a functor of homotopy theories $\Ho \cGst \to \Ho \Spc(B)$ it is necessary for $M^t_{\cGsp} \to M^b_{\cGsp}$ to be an equivalence, since $Y_\bullet \to X_\bullet$ becomes an isomorphism in the homotopy category of $\cGst$. But this is unfortunately false. 
{\em The difficulty is that inclusions of an open subset with large complementary codimension may not be highly connected in the category of motivic spaces.}

\epoint{Example} \label{example:abelian_variety}
Let $X$ be a variety whose associated sheaf on $\Sm_B$ is already $\A^1$-local (for example, an abelian variety). Let $j : X \sm \{ p \} \embed X$ be the complement of a point. This has complementary codimension $\dim{X}$. However,  $j$ is not even $\Aconn{0}$. Indeed, $X \in \PSh(\Sm_k)$ is already a motivic space so $L_{\mot} X = X$. However, $X$ is $0$-truncated, meaning it is valued in sets. The discreteness implies that $\pi_0^{\A^1}(X \sm \{ p \}) \to \pi_0^{\A^1}(X)$ is the complement of a point on its sections over $k$. 

\point 
We are about to remedy this problem in \S \ref{s:remedy}, at the cost of some information, but with the advantage that the functor applies to the entire categories $\cGsp$ and $\cG$.
But for the purposes of Bott periodicity (and the study of classifying spaces such as $\BGL$) we can instead focus on a subcategory  $\cGsp^{\A^1}$ of $\cGsp$ on which all of the desiderata of \S\ref{sec:desiderata} do hold and the realization to motivic spaces does not require any stabilization (see \S \ref{section:refinement}).

\point
\label{s:remedy} The remedy will be straightforward, if perhaps a little surprising:  {\em use the functor 
$\Sigma^2 M_{\cGst}$ instead of
$M_{\cGst}$.} After a double suspension, an open immersion with large codimension boundary will become highly connected. In particular, this means if we ``stabilize all the way to infinity'' and realize all the objects as motivic spectra our desiderata will be achieved. However, it is somewhat remarkable that uniformly only a double suspension is necessary rather than subsequent higher and higher suspensions for more complicated objects. 

\epoint{Remark}  {Recall the infinite suspension $\Sigma^\infty : \Spc(B)_* \to \SH(B)$ commutes with suspension. So if we value our motivic realization in motivic spectra via the composition 
\[ \Sigma^{\infty} M_{\cGst} : \cGst \to \SH(B) \]
then the addition of $\Sigma^2$ becomes nothing more than a shift in the stable category:
\[ \Sigma^2 \Sigma^{\infty} M_{\cGst} \cong \Sigma^{\infty} \Sigma^2 M_{\cGst} \cong (X_\bullet \mapsto \colim\limits_\bullet \Sigma^{\infty} \Sigma^2 M_{\SmArt} X_{\bullet}). \]
}

\bpoint{Connectivity of Open Immersions}

In this section, we first recall the connectivity with respect to the Euclidean topology for open immersions $U \to X$ whose complement has large codimension. The rest of the section will be devoted to understanding how much of these results carry over to the setting of $\A^1$-homotopy theory. Throughout, we will try to give a unified perspective centered on the Blakers-Massey Theorem. For the sake of giving a self-contained treatment for algebraic geometers, we start with a review of this piece of classical topology.

\tpoint{The Blakers-Massey Theorem} {\em \label{t:bm}
Consider a homotopy pushout diagram of topological spaces,
\begin{center}
\begin{tikzcd}
U \arrow[r] \arrow[d] & B \arrow[d]
\\
A \arrow[r] & \pushout X 
\end{tikzcd}
\end{center}
Suppose that the map $U \to A$ is $m$-connected and $U \to B$ is $n$-connected. Then
\begin{enumerate}
\item the map of pairs $(A, U) \to (X, B)$ is $(m + n + 1)$-connected 
\item the map to the homotopy pullback $U \to A \times^h_X B$ is $(m+n)$-connected.
\end{enumerate} }

For a proof, the reader may consult the original article by Blakers and Massey \cite{BM} or a modern treatment such as \cite{ABFJ} and the unpublished work \cite{BMproof}.

\epoint{Remark:  equivalence of (1) and (2) in Theorem~\ref{t:bm}}
There is a natural equivalence between certain homotopy fibers
\[ \fib(U \to A \times^h_X B) \iso \fib(\fib(U \to A) \to \fib(B \to X)) \]
see for example \cite[Prop.~6.3.6]{Hovey}.
This gives the equivalence between (1) and (2) in the Blakers-Massey Theorem~\ref{t:bm}. 

\epoint{Connectivity, cofibers, and suspension}

\tpoint{Proposition} {\em 
Let $f : X \to Y$ be $n$-connected and $X$ be $m$-connected. Then 
\[ \fib(f) \to \Omega \cofib(f) \]
is $(n+m)$-connected. Furthermore, $\cofib(f)$ is $(n+1)$-connected.}

\begin{proof}
Consider the pushout square
\begin{center}
\begin{tikzcd}
X \arrow[r, "f"] \arrow[d] & Y \arrow[d]
\\
* \arrow[r] & \cofib(f).
\end{tikzcd}
\end{center}
The form of the Blakers-Massey Theorem~\ref{t:bm} for connectivity of a map of fibers shows that
\[ \fib(f) \to \fib(* \to \cofib(f)) = \Omega \cofib(f) \]
is $(n + m)$-connected. Therefore,
\[ \pi_k(\fib(f)) \to \pi_k(\Omega \cofib(f)) = \pi_{k+1}(\cofib(f)) \]
is surjective for $k \le n + m + 1$. For any nonempty space $X$, it is $(-1)$-connected by definition. Hence the above map is surjective for $k \le n$. However, if $f$ is $n$-connected then $\fib(f)$ is $n$-connected so the surjection shows that $\pi_{k+1}(\cofib(f)) = 0$ for $k \le n$ proving that $\cofib(f)$ is $(n+1)$-connected.
\end{proof}

\tpoint{Lemma} {\em  \label{lemma:connectivity_spaces}
Let $f : X \to Y$ be $n$-connected. Then $\cofib(f)$ is $(n+1)$-connected. The converse holds if $X$ is $1$-connected and $Y$ is $0$-connected.
}

\begin{proof}
By assumption, $f : X \to Y$ is an isomorphism on $\pi_0$ so $f$ is $(-1)$-connected. To prove $f$ is $n$-connected, we proceed by induction. Assume $f$ is $k$-connected. Then 
\[ \fib(f) \to \Omega \cofib(f) \]
is $(k + 1)$-connected by Blakers-Massey since $X \to *$ is $1$-connected.  But $\Omega \cofib(f)$ is $n$-connected, so for $k < n$ we see that $\fib(f)$ and hence $f$ is $(k+1)$-connected. Hence, by induction,  $f$ is $n$-connected. 
\end{proof}

When $X$ and $Y$ are not $1$-connected, we cannot get started with the induction. Only given information about $\cofib(f)$ we cannot conclude that $\fib(f) \to \Omega \cofib(f)$ has any connectivity in order to conclude some connectivity of $f$ without already assuming some connectivity of either $f$ or $X$ in order to apply the Blakers-Massey Theorem~\ref{t:bm}. However, using suspension to rotate the cofiber sequence $X \to Y \to \cofib(f)$, we can reposition $\cofib(f)$ so that its connectivity is usable as an input to Blakers-Massey. This trick comes at the cost of a suspension.

\tpoint{Lemma} {\em 
If $\cofib(f)$ is $m$-connected then $\Sigma f : \Sigma X \to \Sigma Y$ is $(m-1)$-connected.
}

\begin{proof}
Consider the homotopy cofiber diagram
\begin{center}
\begin{tikzcd}
X \arrow[d] \arrow[r, "f"] & Y \arrow[d, "r"]
\\
* \arrow[r] & \cofib(f).
\end{tikzcd}
\end{center}
We can ``rotate'' as in the formation of the Puppe sequence via taking pushouts to get a diagram consisting of homotopy pushout squares
\begin{center}
\begin{tikzcd}
X \arrow[d] \arrow[r, "f"] & Y \arrow[d, "r"] \arrow[r] & * \arrow[d]
\\
* \arrow[r] & \cofib(f) \arrow[r] \arrow[d] & \Sigma X \arrow[d, "-\Sigma f"] \arrow[r] & * \arrow[d]
\\
& * \arrow[r] & \Sigma Y \arrow[r] & \Sigma \cofib(f)
\end{tikzcd}
\end{center}
Applying Blakers-Massey to the ``rotated'' homotopy pushout square,
\begin{center}
\begin{tikzcd}
\cofib(f) \arrow[r] \arrow[d] & \Sigma X \arrow[d, "-\Sigma f"] 
\\
* \arrow[r] & \Sigma Y
\end{tikzcd}
\end{center}
we conclude using that $\cofib(f) \to \Sigma X$ is $(-1)$-connected (since $X$ is nonempty).
\end{proof}

\epoint{Remark}
This result is optimal.  Indeed, consider $S^1 \to *$, which has cofiber $\Sigma S^1 = S^2$, which is $1$-connected. Then $S^1 \to *$ is $0$-connected but not $1$-connected.

\epoint{Open Immersions in the Euclidean Topology}

We next address the  connectivity of the inclusion map $j : U \to X$ of a open set of a smooth variety with complement of  large codimension.

\tpoint{Proposition} {\em 
Let $X$ be a manifold and $Z \subset X$ a submanifold of (real) codimension $r$. Let $U = X \sm Z$. Then the inclusion $j : U \to X$ is $(r-2)$-connected.}

\begin{proof}
We apply the Blakers-Massey Theorem~\ref{t:bm} to the homotopy pushout diagram
\begin{center}
\begin{tikzcd}
N_Z \cap U \arrow[d] \arrow[r] & U \arrow[d]
\\
N_Z \arrow[r] & \pushout X
\end{tikzcd}
\end{center}
where $N_Z$ is a disk bundle variant of the normal bundle which is isomorphic to a tubular neighborhood of $Z$ inside $X$. Notice that the intersection $C := N_Z \cap U$ is homotopy equivalent to the sphere bundle $S(N_Z)$ of the normal bundle. 

Note that the projection $N_Z \to Z$ is a homotopy equivalence so we may replace the map $N_Z \cap U \to N_Z$ by the projection $N_Z \cap U \to Z$. Since $N_Z$ has rank $r$, the bundle $S(N_Z) \to Z$ is an $S^{r-1}$-bundle and hence is $(r-2)$-connected. Furthermore, $S(N_Z) \to U$ is $(-1)$-connected because both are connected. Therefore, $(N_Z, S(N_Z)) \to (X, U)$ is $(r-2)$-connected meaning
\[ \pi_k(N_Z, S(N_Z)) \to \pi_k(X,U) \]
is an isomorphism for $k \le r - 2$ and surjective for $k = r - 1$. However, the pair $(N_Z, S(N_Z))$ is also $(r-1)$-connected (since $S(N_Z) \to Z$ is $(r-2)$-connected) meaning $\pi_k(N_Z, S(N_Z)) = 0$ for $k \le r - 1$.  Thus $\pi_k(X, U) = 0$ for $k \le r - 1$.  
\end{proof}

\epoint{Remark}
This is optimal for each $r$. Indeed, let $U \embed S^n$ be the complement of one point. Then $U \cong \RR^n$ so the inclusion is exactly $(n-2)$-connected.

\tpoint{Corollary} {\em 
Let $X$ be a complex manifold (e.g.\ a smooth $\CC$-variety) and $Z \subset X$ a submanifold of (complex) codimension $r$. Let $U = X \sm Z$. Then the inclusion $j : U \to X$ is $(2r-2)$-connected.}

\epoint{Remark}
This is also optimal. For example, for $r = 1$ consider $j : E \sm \{ p \} \embed E$ the complement of a point in an elliptic curve. This is surjective on $\pi_1$ but not an isomorphism. For $\A^1 \to \P^1$, the complement of the point at $\infty$, the induced map on $\pi_2$ is not surjective.

\tpoint{Corollary} {\em 
Let $X$ be a smooth $\CC$-variety and $Z \subset X$ a subvariety (not necessarily smooth). Let $U = X \sm Z$. Then the inclusion $j : U \to X$ is $(2r - 2)$-connected.}

\begin{proof} (Cf.\ the proof of Lemma~\ref{l:twofour}.) 
Consider a filtration
\[ Z = Z_0 \supsetneq Z_1 \supsetneq \cdots \supsetneq Z_k = \varnothing \]
where $Z_i \sm Z_{i+1}$ is smooth. Such a filtration exists by generic smoothness. We proceed by induction on $r$. For $r = \dim{X}$ the subvariety consists of a finite set of points and hence is automatically smooth. Now we factor
\[ \pi_k(X \sm Z) \to \pi_k(X \sm Z_1) \to \pi_k(X). \]
Since $Z_1$ is codimension $\ge r + 1$ in $X$ by the induction hypothesis $X \sm Z_1 \to X$ is $2r$-connected. Now $X \sm Z \to X \sm Z_1$ is the complement of the smooth subvariety $Z \sm Z_1$ so by our previous results it is $(2r-2)$-connected. Hence, both maps are isomorphisms for $k \le 2r-2$ and surjective for $k = 2r - 1$.
\end{proof}

\epoint{Connectivity for Open Immersions in Motivic Spaces}

Unfortunately, $\A^1$-connectivity results for $j : U \embed X$ do not work analogously to connectivity results for complex manifolds. One manifestation is that for varieties defined over $\RR$, the motivic spaces have an $\RR$-realization, so the best one could hope for is that if the codimension of $Z = X \sm U$ is $r$ then $j$ is $(r-2)$-connected. However, $j$ may not even be $0$-connected in the sense of motivic homotopy theory (see Example~\ref{example:abelian_variety}). The reason the proof fails is that we don't have access to tubular neighborhoods or sphere bundles in motivic homotopy theory. Instead, we have the purity isomorphism which will tell us that $\cofib(U \to X)$ is highly $\A^1$-connected. Here we follow the exposition of \cite[\S8]{HJNTY}. For the remainder of this section, let $k$ be a perfect field. In the following lemma, the slightly weakened connectivity hypotheses for the converse, compared to \cite[Lemma 8.8]{HJNTY}, will be necessary in \S\ref{section:refinement}.

\tpoint{Lemma~\cite[Lemma 8.8]{HJNTY}} {\em  \label{lem:motivic_connectivity} 
Let $f : X \to Y$ be a morphism in $\PSh(\Sm_k)$ and $n \ge -1$ an integer. If $f$ is $\A^1$-$n$-connected, then $\cofib(f)$ is $\A^1$-$(n+1)$-connected. The converse holds if $X$ is $\A^1$-$1$-connected and $Y$ is $\A^1$-$0$-connected.
}

\begin{proof}
Since $L_{\mot}$ is a left-adjoint, it preserves homotopy colimits. Hence $L_{\mot}(\cofib(f)) = \cofib_{\mot}(L_{\mot}(f))$. However, these homotopy colimits are computed in different categories (in particular, $\cofib_{\mot}$ computed in motivic spaces may not be compatible with the cofiber in spaces computed stalkwise) so we cannot just apply Blakers-Massey to $L_{\mot}(f)$. 
\par 
Instead, let $C$ be the cofiber of $L_{\mot}(f)$ \textit{computed in the sheaf category} $\Shv_{\Nis}(\Sm_k)$. By Lemma~\ref{lemma:connectivity_spaces} (applied to stalks), since $L_{\mot}(f)$ is $n$-connected, by assumption, we conclude that $C$ is $(n+1)$-connected. By Morel's $\A^1$-Connectivity Theorem~\ref{t:connectivity}, this implies that $L_{\mot} C$ is also $(n+1)$-connected. Now $C = L_{\Nis} \cofib(L_{\mot} f)$ so $L_{\mot} C = L_{\mot} \cofib(L_{\mot} f) = \cofib_{\mot}(L_{\mot} f) = L_{\mot} \cofib(f)$ and thus $\cofib(f)$ is $\Aconn{(n+1)}$.
\par 
Conversely, suppose that $X$ and $Y$ are $\Aconn{1}$ and that $L_{\mot} \cofib(f)$ is $(n+1)$-connected.  Consider $C = L_{\Nis} \cofib(L_{\mot} f)$.  We will prove that $C$ is $(n+1)$-connected and $L_{\mot} f$ is $n$-connected. This follows exactly the inductive strategy of Lemma~\ref{lemma:connectivity_spaces}. Let $F = L_{\Nis} \fib(L_{\mot} f)$ be the homotopy fiber computed in $\PSh(\Sm_k)$. Since $L_{\mot} X$ is $1$-connected and $L_{\mot} Y$ is $0$-connected, $L_{\mot} f$ is $(-1)$-connected, meaning $F$ is $(-1)$-connected, and $C$ is $0$-connected. The issue is that $C$ is \textit{not} motivic since it is the cofiber computed in $\PSh(\Sm_k)$ not in $\Spc(k)$. Therefore, our assumption only gives connectivity of $L_{\mot} C$ not of $C$ itself so we cannot directly run the induction. We need the partial converse to Morel's connectivity theorem given in \cite[Lemma~6.60]{Mor12} to run the induction. The induction on $k$ will assume $L_{\mot} f$ is $k$-connected and $C$ is $(k+1)$-connected starting with $k = -1$. 
\par 
Assuming $L_{\mot} f$ is $k$-connected, the Blakers-Massey Theorem~\ref{t:bm} shows that $F \to \Omega C$ is $(k+1)$-connected since $L_{\mot} X$ is $1$-connected. This gives an isomorphism $\pi_{k+1}(F) \iso \pi_{k+2}(C)$. For $k = -1$, we also know from the long exact sequence $\pi_1^{\A^1}(Y) \iso \pi_0(F)$ so $\pi_0(F)$ is strongly $\A^1$-invariant\footnote{Ayoub has found an example \cite{Ayoub:example} showing that $\pi_0(\X)$ for $\X$ motivic may not even be $\A^1$-invariant so this remark is necessary.}. For $k \ge 0$, the isomorphism $\pi_{k+1}(F) \iso \pi_{k+2}(C)$ shows that $\pi_{k+2}(C)$ is strictly (for $k = 0$ the isomorphism implies it is abelian) $\A^1$-invariant since $F$ is motivic because the inclusion $\Spc(k) \embed PSh(\Sm_k)$ preserves homotopy limits. If we also assume $C$ is $(k+1)$-connected, then $\pi_{i}(C) = 0$ for $i < k+2$. By \cite[Lemma~6.60]{Mor12} $\pi_{k+2}(C) = \pi_{k+2}^{\A^1}(C) = \pi_{k+2}^{\A^1}(L_{\mot} C)$. But $L_{\mot} C$ is $(n+1)$-connected so $\pi_{k+1}(F) = \pi_{k+2}(C) = 0$ as long as $k + 2 \le n + 1$ proving the claim by induction. Therefore, we conclude that $L_{\mot} f$ is $n$-connected. 
\end{proof}

\tpoint{Corollary} \label{cor:suspension_connectivity} {\em 
If $X \in \PSh(\Sm_k)$ is $\Aconn{n}$ then $\Sigma X$ is $\Aconn{(n+1)}$ meaning the motivic space $\Sigma L_{\mot} X \cong L_{\mot} \Sigma X$ is $(n+1)$-connected.}

\begin{proof}
$\Sigma X$ is the cofiber of the map $X \to *$ which is $\Aconn{(n+1)}$ by assumption. 
\end{proof}

We similarly  need a lemma ensuring that the cofiber of a high codimension open embedding is highly $\A^1$-connected. When the complement is smooth, this follows from the purity isomorphism, Theorem~\ref{t:purity}.  We continue to rely on the arguments of \cite[\S 8]{HJNTY}.

\tpoint{Lemma} \label{lem:thom_connectivity} {\em Let $\cE$ be a vector bundle of rank $r$ on a smooth $k$-scheme $X$. Then $\Th(\cE)$ is $\Aconn{(r-1)}$. 
}

\begin{proof}
Nisnevich (indeed, Zariski) locally on $X$, the vector bundle $\cE$ is trivial. By Corollary~\ref{cor:connectivity_stalks}, connectivity is Nisnevich-local on the target (one can also argue by replacing $X$ by the \v{C}ech nerve $U_\bullet$ of a Nisnevich cover $U \to X$ since $U_\bullet \to X$ becomes an equivalence in $\Spc(B)$). For the trivial bundle, $\Th(\cE) \cong_{\mot} \Sigma^r X_+$. Then applying Corollary~\ref{cor:suspension_connectivity} $r$ times to the $(-1)$-connected space $X_+$, we conclude that $\Th(\cE)$ is $\Aconn{(r-1)}$.
\end{proof}

\tpoint{Lemma  \cite[Lemma 8.9]{HJNTY}} {\em  \label{lem:cofib_open_embedding_connected}
Let $X$ a smooth $k$-scheme, and $Z \subset X$ a closed subscheme of codimension $\ge r$. Then $\Sigma(X/(X \sm Z))$ is $\Aconn{r}$. }

\begin{proof}
$X$ is a disjoint union of quasi-compact smooth schemes so we may assume that $X$ is quasi-compact. If $Z$ is smooth then $X / (X \sm Z)$ is $\Aconn{(r-1)}$ by the purity isomorphism, Theorem~\ref{t:purity}, and Lemma~\ref{lem:cofib_open_embedding_connected}.
For the general case, we first reduce to the case where $Z$ is reduced. Indeed, the sheaves represented by $Z$ and $Z_{\red}$ on $\Sm_k$ are the same since a morphism from a reduced scheme factors through the reduction of the target.
\par 
Since $k$ is perfect, generic smoothness implies that there exists a filtration
\[ \varnothing = Z_0 \subsetneq Z_1 \subsetneq \cdots \subsetneq Z_n = Z \]
of closed subschemes such that $Z_j \sm Z_{j-1}$ is smooth. We prove the result by induction on $n$. For $n = 0$, we verify that $X / X  = *$ is indeed $r$-connected. Consider the cofiber sequence,
\[ \frac{X \sm Z_{n-1}}{X \sm Z} \to \frac{X}{X \sm Z} \to \frac{X}{X \sm Z_{n-1}} \]
Therefore, because homotopy colimits commute, we get a cofiber sequence,
\[ \Sigma \left( \frac{X \sm Z_{n-1}}{X \sm Z} \right) \to \Sigma \left( \frac{X}{X \sm Z} \right) \to \Sigma \left( \frac{X}{X \sm Z_{n-1}} \right). \]
By the induction hypothesis, $\Sigma X / (X \sm Z_{n-1})$ is $\Aconn{r}$. If we can show that the first two terms are $\Aconn{1}$ then the converse part of Lemma~\ref{lem:motivic_connectivity} proves that the morphism,
\[ \Sigma \left( \frac{X \sm Z_{n-1}}{X \sm Z} \right) \to \Sigma \left( \frac{X}{X \sm Z} \right)  \]
is $\Aconn{(r-1)}$. Indeed $X / (X \sm Z)$ is $\Aconn{0}$ by \cite[Lemma 6.1.4]{Mor05} for any closed embedding $Z \embed X$ where $X$ is smooth and irreducible. Hence, the first two terms are $\Aconn{1}$.

Since, $X \sm Z \to X \sm Z_{n-1}$ has smooth complement $Z \sm Z_{n-1}$, we can apply the smooth case to show that the first term, $\Sigma ((X \sm Z_{n-1})/(X \sm Z))$ is $\Aconn{r}$ so connectivity of the map implies that $\Sigma( X / (X \sm Z))$ is $\Aconn{r}$.
\end{proof}

It will also be useful to record a similar result trading the partial stabilization in the conclusion for an additional connectivity hypothesis. 

\tpoint{Lemma} \label{lem:connectivity_inclusion_with_connectivity_assumption} {\em Let $j : U \to X$ be an open embedding of a smooth $k$-scheme. Suppose that $U$ is $\Aconn{1}$ and the complement $Z$ of $j$ has codimension at least $r$.  Then $j$ is $\Aconn{(r-2)}$. }

\bpoint{Realizing the class ($k$-conn) in motivic homotopy theory} 

From now on, we set our base $B$ to be the spectrum of a perfect field $k$. (One might hope the constructions performed in $\cGst$ over a base scheme $B = \Spec{\Z}$ would extend to similar constructions in $\Spc(\Z)$. However, because various lemmas rely on Morel's $\A^1$-Connectivity Theorem~\ref{t:connectivity} which --- to the author's knowledge --- is known only when the base is a perfect field, we have not realized this hope.)

\epoint{Recalling definitions}
First, we recall some definitions from the main part of the article.
An $\A^n$-bundle $\pi : E \to X$ is a (smooth-locally-trivial) torsor for a vector bundle of rank $n$ over $X$. In particular, it is a smooth morphism whose fibers are isomorphic to $\A^n$.
The definition when 
a morphism $X \to Y$ in $\SmArt$ is in the class (nicely-$k$-conn) or
 ($k$-conn) was given in \S \ref{s:niceconn}.
 
\tpoint{Lemma} {\em 
Let $f : X \to Y$ be a ($k$-conn) morphism in $\SmArt$. Then $\Sigma^2 M_{\SmArt}(f)$ is $\Aconn{(k-1)}$. In fact, it induces an isomorphism on $\pi_i$ for all $i \le k$ (but may fail to be $\Aconn{k}$ if the induced map on $\pi_{k+1}^{\A^1}$ is not an epimorphism).
} 

\begin{proof}
Let $I_{k}$ be the class of $1$-morphism $f : X \to Y$ in $\Spc(k)$ inducing isomorphisms $\pi_i(X) \to \pi_i(Y)$ for all $i \le k$ (and all choices of base point $* \to X$). The reason to consider this property, rather than the more natural property of $(k-1)$-connectivity is that it enjoys a two-out-of-three property by functoriality of $\pi_i$ while connectivity does not (because the two-out-of-three property fails for epimorphisms).

By the definition of ($k$-conn) it suffices to show that both affine bundles and $k$-codimension open subsets are in $I_{k}$. Since both are examples of (nicely-$k$-conn) morphisms (the first with $i = \id$, the second with $a = \id$) it will be convenient to prove by hand that $\Sigma^2 M_{\SmArt}$ of any (nicely-$k$-conn) morphism is actually $k$-connected. This is done in the following lemma.
\end{proof}

\tpoint{Lemma (cf.\ \cite[Proposition 8.10]{HJNTY})} {\em  \label{lem:M_comparison}
Let $f : X \to Y$ be a morphism of $\SmArt$ in (nicely-$k$-conn).  Then $\Sigma^2 M_{\SmArt}(f)$ is $k$-connected.}

\begin{proof}
The proof follows the same strategy as \cite[Proposition 8.10]{HJNTY} except we also must handle the case of arbitrary $\A^n$-bundles rather than only vector bundles.
Every presheaf $Y \in \PSh(\Sm_k)$ can be written as the filtered homotopy colimit of representable presheaves, say by some collection $\{ Y_i \}_{i \in I}$ with $Y_i$ smooth $k$-schemes. Since filtered homotopy colimits commute with products we see that
\[ \colim\limits_{i \in I} (X \times_Y Y_i) = X \times_Y \colim\limits_{i \in I} Y_i = X \]
and since $X \to Y$ is smooth and representable (as it is a composition of an open embedding an $\A^n$-bundle) the presheaves $X \times_Y Y_i$ are represented by smooth schemes. Hence $f : X \to Y$ is the homotopy colimit of $f_i : X_i \to Y_i$ with $X_i := X \times_Y Y_i$ a smooth $k$-scheme and $f_i$ a (nicely-$k$-conn) morphism. Since $L_{\mot}$ commutes with homotopy colimits and filtered homotopy colimits preserve connectivity (checking on stalks and using that filtered homotopy colimits commute with homotopy groups) we reduce to the case that $Y \in \Sm_k$. Write $X \to Y$, factored as
\[ X \embed V \to Y \]
where $X \embed V$ is an open embedding which has complementary codimension $\ge k$ on each fiber and $\pi : V \to Y$ is an affine space bundle. Since $Y$ is a scheme, $\pi : V \to Y$ is Zariski-locally trivial. Pick a Nisnevich cover $W \to Y$ trivializing $V \to Y$. Then the \v{C}ech nerve $W_\bullet$ of $W \to V$ gives a Nisnevich-local equivalence $W_\bullet \to V$. Since $V \times_Y W_\bullet \to W_\bullet$ is a trivial affine space bundle it is an $\A^1$-equivalence. Hence $V \to Y$ is a motivic equivalence. By assumption, $X \embed V$ has complement of codimension $\ge k$ so by Lemma~\ref{lem:cofib_open_embedding_connected} we have $\Sigma (V/X)$ is $\Aconn{k}$ and hence $\Sigma^2 (V/X)$ is $\Aconn{(k+1)}$. Furthermore, $\Sigma^2 X$ is $1$-connected hence $\Aconn{1}$ by Morel's $\A^1$-Connectivity Theorem~\ref{t:connectivity}. Therefore, we can apply the converse part of Lemma~\ref{lem:motivic_connectivity} to the cofiber sequence
\[ \Sigma^2 X \to \Sigma^2 V \to \Sigma^2 (V/X) \]
to conclude that $\Sigma^2 X \to \Sigma^2 V$ is $\Aconn{k}$. 
\end{proof}

\tpoint{Theorem} {\em \label{thm:top_bottom_comparison}
The $1$-morphism $\Sigma^2 M_{\cGsp}^t \to \Sigma^2 M_{\cGsp}^b$ in $\Fun{\cGsp}{\Spc(k)}$ is an equivalence. Furthermore, for any $X_\bullet \in \cGst$ there is an map
$X_{\ell} \to X_\bullet$ in $\cGst$ represented by the diagram
\begin{center}
    \begin{tikzcd}
        \cdots \arrow[r] & \varnothing \arrow[r] \arrow[d] & X_{\ell} \arrow[d] \arrow[r] & X_{\ell} \arrow[d] \arrow[r] & X_{\ell} \arrow[d] \arrow[r] & \cdots
        \\
        \cdots \arrow[r] & X_{\ell-1} \arrow[r] & X_{\ell} \arrow[r] & X_{\ell+1} \arrow[r] & X_{\ell + 2} \arrow[r] & \cdots
    \end{tikzcd}
\end{center}
Then $\Sigma^2 M_{\cGst} X_{\ell} \to \Sigma^2 M_{\cGst} X_\bullet$ is $k(\ell)$-connected. }

\begin{proof}
It suffices to show that for each $(Y_\bullet \to X_\bullet)$ the morphism 
\[ \Sigma^2 M_{\cGsp}^t (Y_\bullet \to X_\bullet) \to \Sigma^2 M_{\cGsp}^b (Y_\bullet \to X_\bullet) \]
is an equivalence. Unwinding the definitions, we need to prove that
\[ \colim\limits_{i \in I} \Sigma^2 M_{\SmArt} Y_i \to  \colim\limits_{i \in I} \Sigma^2 M_{\SmArt} X_i \]
is an equivalence in $\Spc(k)$. Since $L_{\mot}$ is a left adjoint, we can compute the homotopy colimit in $\PSh(\Sm_k)$ and then apply $L_{\mot}$. By  Morel's $\A^1$-Connectivity Theorem~\ref{t:connectivity}, it suffices to prove that the homotopy colimit computed in presheaves is $n$-connected for all $n$, hence an equivalence by hypercompleteness. From the definition of $\cGsp$ there is a weakly increasing unbounded function $k : \Z^+ \to \Z^+$ such that $Y_\ell \to X_{\ell}$ is nicely $k(\ell)$-connected. By Lemma~\ref{lem:M_comparison} the map 
\[ \Sigma^2 M_{\SmArt} Y_{\ell} \to \Sigma^2 M_{\SmArt} X_{\ell} \]
is $\Aconn{k(\ell)}$. Denote the stalk at $x \in W \in \Sm_k$ of this morphism by $F_i \to G_i$ and likewise $F \to G$ for the homotopy colimit (which equals the stalk of the colimit computed in sheaves). Since filtered homotopy colimits commute with homotopy groups, we see that $\pi_j(F) \to \pi_j(G)$ is the homotopy colimit of $\pi_j(F_i) \to \pi_j(G_i)$. These are isomorphisms for $i$ large enough that $k(i) \ge j$. Therefore, $\pi_j(F) \to \pi_j(G)$ is an isomorphism. 

To show that $\Sigma^2 M_{\cGst} X_{\ell} \to \Sigma^2 M_{\cGst} X$ is highly connected, we repeat the same argument and reduce to showing that if $\pi_j(F_{\ell}) \to \pi_j(F_i)$ is an isomorphism for all $i \ge \ell$, then $\pi_j(F_{\ell}) \to \pi_j(F)$ is an isomorphism. This is clear because $\pi_j$ commutes with homotopy colimits so $\pi_j(F_{\ell}) \to \pi_j(F) = \colim\limits_{i \ge \ell} \pi_j(F_i)$ is an isomorphism.
\end{proof}

\bpoint{Refinement for the $\mathbb{A}^1$-connected subcategory} \label{section:refinement}

The stabilization $\Sigma^2$ is really only necessary to deal with disconnected spaces. When working with classifying spaces for $\GL_n$ the objects of interest are certain open sets of Grassmannians which are indeed highly connected. In this section, we show that there is a subcategory $\cGsp^{\A^1} \subset \cGsp$ of $\Aconn{1}$ spatial realizations that are directly realizable in $\Spc(k)$ without the double stabilization. 

Let $\cGsp^{\A^1}$ be the full subcategory of $\cGsp$ whose objects are diagrams $Y_\bullet \to X_\bullet$ so that for all $\ell \ge \ell_0$ larger than some cutoff $\ell_0$ (depending on $Y_\bullet \to X_\bullet)$ the motivic space $M_{\SmArt}(Y_{\ell})$ is $1$-connected {\em i.e.} $\pi_i^{\A^1}(Y_{\ell}) = 0$ for $i = 0,1$.

Now we can use the motivic realizations naively as desired in \S \ref{appendixintroduction}, we write $M_{\cGsp^{\A^1}}^{t}$ (resp.\ $M_{\cGsp^{\A^1}}^{b}$) for the restriction of $M_{\cGsp}^{t}$ (resp. $M_{\cGsp}^{b}$) to not cause any confusion with the source of the functors. 

\tpoint{Theorem} \label{thm:A1_top_bottom_comparison}
{\em The $1$-morphism $M_{\cGsp^{\A^1}}^t \to M_{\cGsp^{\A^1}}^b$ in $\Fun{\cGsp^{\A^1}}{\Spc(k)}$ is an equivalence. Furthermore, consider the morphism
$X_{\ell} \to X_\bullet$ in $\cGst^{\A^1}$.  Then $M_{\SmArt} X_{\ell} \to \colim M_{\SmArt} X_\bullet$ is $(k(\ell)-3)$-connected. }

\begin{proof}
The proof is identical to Theorem~\ref{thm:top_bottom_comparison} except now we must modify the application of Lemma~\ref{lem:M_comparison} to use the extra assumptions on connectivity. As before, it suffices to show that $M_{\SmArt} Y_{\ell} \to M_{\SmArt} X_{\ell}$ is $\Aconn{(k(\ell)-2)}$. Recall by assumption $Y_{\ell}$ is $\Aconn{1}$ connected and $Y_{\ell} \to X_{\ell}$ factors as a closed embedding into an $\A^n$-bundle whose complement fiberwise is smooth of codimension $\ge k(\ell)$. By reducing to representable presheaves, it suffices to show that if $U \embed Z$ is an open immersion of smooth varieties whose complement has codimension $\ge k$ and $U$ is $\Aconn{1}$ then $U \embed Z$ is $\Aconn{(k-2)}$. This is exactly the content of Lemma~\ref{lem:connectivity_inclusion_with_connectivity_assumption}.
\end{proof}

\bpoint{The Motivic Spectrum $\KGL$}

For a sheaf of groups $G$, there are two sorts of classifying spaces usually considered in motivic homotopy theory. Broadly speaking, these classify Nisnevich-locally or \etale-locally trivial $G$-torsors. Luckily, by descent, these coincide for $G = \GL_n$ (cf.\ Hilbert's theorem 90).
\par 
Let us first recall the constructions of these classifying spaces in $\Spc(B)$. 

\epoint{Definition}
Let $G$ be a presheaf of groups on $\Sm_B$ and $B_\bullet G \in \PSh(\Sm_B)$ be the simplical presheaf obtained via the bar construction. Explicitly, $B_\bullet G$ as a simplicial object in presheaves of sets is 
\[ B_\bullet G : \Prism^{\op} \to \Fun{\Sm_B^\op}{\mathbf{Set}} \]
given by
\[ [n] \mapsto G^{\times n} \]
with face maps $d_i : G^{\times n} \to G^{\times (n-1)}$ given by
\[ d_i(g_1, \dots, g_n) = 
\begin{cases}
(g_2, \dots, g_n) & i = 0
\\
(g_1, \dots, g_i g_{i+1}, g_{i+2}, \dots, g_n) & 1 \le i \le n - 1
\\
(g_1, \dots, g_{n-1}) & i = n
\end{cases} \]
and degeneracies $s_i : G^{\times n} \to G^{\times (n+1)}$ given by
\[ s_i(g_1, \dots, g_b) = \begin{cases}
(e, g_1, \dots, g_n) & i = 0
\\
(g_1, \dots, g_i, e, g_{i+1}, \dots, g_n) & 1 \le i \le n
\end{cases}\]
Notice that $B_\bullet G$ is exactly the \v{C}ech nerve of the morphism of stacks $S \to \mathbf{B}_{S, \tau} G$ where $\mathbf{B}_{S, \tau} G$ is the classifying stack of $G$-torsors over $S$ for a Grothendieck topology $\tau$. We may then view this simplicial object as an element of $\PSh(\Sm_B)$ by taking the realization of the simplicial set as a space.

\tpoint{Proposition  \cite[Lemma 2.3.2]{AHW18}} {\em
The induced morphism $B_\bullet G \to \mathbf{B}_{\tau} G$ is a $\tau$-local equivalence.}

In particular, $B_\bullet G \to \mathbf{B}_{\Nis} G$ is a Nisnevich local equivalence and hence \textit{a fortiori} a motivic equivalence. If we define $B_{\mot} G := L_{\mot} B_{\bullet} G$ then the natural map $B_{\mot} G \to L_{\mot} \mathbf{B}_{\Nis} G$ is an equivalence. 
\par 
Since we only include Nisnevich locally trivial torsors, the Nisnevich stack $\mathbf{B}_{\Nis} G$ is not an algebraic stack in the usual sense (it does not even satisfy \etale gluing). We could also consider $B_{\et} G := M_{\SmArt}(B G)$ where $BG \in \SmArt$ denotes the usual classifying stack of a smooth group scheme (hence $G$-torsors are \etale locally trivial). Morel and Voevodsky have a detailed discussion of this object in motivic homotopy theory \cite[\S4.2]{mv}. Thankfully, for $G = \GL_n$ the situation is made very simple by descent since every $\fppf$-torsor is Zariski-locally trivial. In particular $\mathbf{B}_{\Nis} \GL_n = B_{\et} \GL_n = \BGL_n$.
\par 
From the equivalence $M^t_{\cGsp^{\A^1}} \to M^b_{\cGsp^{\A^1}}$ applied to the spatial realization in $\cGsp^{\A^1}$
\begin{center}
\begin{tikzcd}
 \Gr(n, n+1) \arrow[d] \arrow[r, hook] & \Gr(n, n+2) \arrow[d] \arrow[r] & \Gr(n, n+3) \arrow[d] \arrow[r] & \cdots 
 \\
 \BGL_n \arrow[r, equals] & \BGL_n \arrow[r, equals] & \BGL_n \arrow[r] & \cdots 
\end{tikzcd} 
\end{center}
whose downward maps classify the universal $n$-space bundle over $\Gr(n, n+k)$, we obtain a motivic equivalence in $\Spc(k)$
\[ \Gr(n, \infty) := \colim\limits_k L_{\mot} \Gr(n, n+k) \iso \BGL_n \]
recovering \cite[Proposition~3.7, p.138]{mv}.

Furthermore, there are stabilization maps $\BGL_n \to \BGL_{n+1}$ given by $E \mapsto E \oplus \cO$ and we define in $\Spc(k)$ the motivic space $\BGL := \colim\limits_n \BGL_n$.

\tpoint{Corollary} {\em 
Recall the  spatial realization 
\eqref{eq:BGLsp}
of ``$\BGL$''.
The equivalence $M^t_{\cGsp^{\A^1}} \iso M^b_{\cGsp^{\A^1}}$ induces an equivalence in $\Spc(k)$
\[ \colim\limits_n \Gr(n, 2n) \iso \BGL. \]
}

One of the primary motivations for the category $\Spc(B)$ of motivic spaces is the representability of algebraic $K$-theory by, not quite a motivic space, but rather a motivic spectrum. 

\tpoint{Theorem} {\em Let $B$ be a regular Noetherian scheme of finite Krull dimension.
There exists a motivic spectrum $\KGL \in \SH(B)$ representing algebraic $K$-theory in the sense that $[\Sigma^n \Sigma^{\infty} X_+, \KGL]_{\SH(B)} \cong K_n(X)$ for any $X \in \Sm_B$. Furthermore, there is a motivic equivalence $\Z \times \BGL \iso \Omega^{\infty} \KGL$.  Therefore, by adjunction, $K$-theory is representable by a motivic space
\[ K_n(X) \cong [\Sigma^n \Sigma^{\infty} X_+, \KGL]_{\SH(B)} = [\Sigma^n X_+, \Z \times \BGL]_{\Spc(B)_*}. \]
}

See  \cite[\S 2.3]{bachmann} for a proof and discussion of this result.

\epoint{Remark}
As usual, $K = \Omega^{\infty} \KGL$ is pointed at $0 \in K(S)$ and therefore we may also define the reduced $K$-theory $\tilde{K}_n(X)$ of a pointed space via 
\[ \tilde{K}_n(X) := [\Sigma^n \Sigma^{\infty} X, \KGL]_{\SH(B)} = [\Sigma^n X, K]_{\Spc(B)_*}\]
As with topological $K$-theory, this corresponds to classes that are trivial over the basepoint.  (For connected spaces this is equivalent to classes of rank $0$.)

\epoint{Remark}
We need the regularity hypotheses on $B$ for algebraic $K$-theory to be connective. Otherwise, we need to consider the so-called \textit{connective} $K$-theory for the above statements to hold.

\bpoint{Comparison of Bott Maps} \label{appendix:section:comparison}

\renewcommand{\E}{\mathcal{E}}

We set $\Omega_Y$ to be the right adjoint to $Y \wedge (-)$ where $Y$ is a pointed motivic space considered as an endofunctor of the category of pointed motivic spaces. Note that  $\Omega_Y$ exists because $\mathrm{Mor}_{\Spc(B)_*}(Y \wedge (-), X)$ commutes with homotopy colimits and forms a pointed $\A^1$-invariant sheaf when evaluated on objects of $\Sm_B$ for any pointed motivic space $X$. Set $K = \Omega^{\infty} \KGL$,  which is the space representing algebraic $K$-theory as we saw in the previous section. The motivic space $K$ comes equipped with a natural map
\[ \beta : K \to \Omega_{\P^1} K \]
called the Bott map. Since $\KGL$ is a ring spectrum, such a map is induced by the map
\[ \Z = [*_+, K]_* \to [*_+, \Omega_{\P^1} K]_* = [\P^1, K]_* = \tilde{K}_0(\P^1) \]
defined by mapping $1$ to the element $\beta = (\gamma - 1)$  (here $\gamma$ is the class of the tautological bundle $\cO_{\P^1}(-1)$) which we call the Bott element. 
\par 
Indeed, we can produce the Bott map using tensor product and adjunction. The map
\[ \beta : K_0(X) \to K_0(\P^1 \times X) \]
given by the difference of the $- \ot \gamma$ map and the pullback induces, by Yoneda and adjunction, a map
\[ K \to \Omega_{\P^1_+} K \]
since $\P^1_+ \wedge X = \P^1 \times X$. This is itself, adjoint to a map
\[ \P^1_+ \wedge K \to K. \]
To show it factors through $\P^1 \wedge K \to K$ it suffices to show that $*_+ \wedge K \to K$ is equivalent to $0$. This is exactly the statement that restricting $\gamma - 1$ to $* \embed \P^1$ gives a trivial class. 

\tpoint{Theorem} {\em Let $\BGL \in \cGst$ represent the stabilized classifying stack as in Definition~\ref{d:BGL} and $\Omega^2_{\alg}(\BGL)$  be the disjoint union over $d$ of one of the sequences described in Proposition~\ref{p:final} as described in Definition~\ref{d:final}. Let $\beta_{LV} : \Z \times \BGL \to \Omega^2_{\alg}(\BGL)$ be the union over $d$ of the inverses of the maps $\mu : \Omega^2_{ \alg,d}(\BGL) \to \BGL$ described in Theorem~\ref{thm:bpg}. Then there exists a square in the category of motivic spaces
\begin{center}
\begin{tikzcd}
M_{\cGst}(\Z \times \BGL) \arrow[r, "\beta_{LV}"] \arrow[d] & M_{\cGst} \Omega^2_{\mathrm{alg}}(\BGL) \arrow[d]
\\
K \arrow[r, "-\beta"] & \Omega_{\P^1} K
\end{tikzcd}
\end{center}
commutative up to homotopy whose morphisms are all motivic equivalences.}

\begin{proof}
The diagram arises from the following diagram of presheaves on $\Sm_k$
\begin{equation}\label{eq:benaug30}\xymatrix{& T'_{r,d}(m) \ar[dr] \ar[dl]_{\al} \\ 
BGL(d+mr)\ar[d]^{[-]-mr} & & U_{r,d}(m) \ar[d]^{[-]-r} \\ K \ar[rr]^{-\be} & & \Omega_{\proj^1} K }\end{equation}
Here $[-]$ is the map taking a vector bundle to a $K$-theory class using that the values $[X, K]_{\PSh(\Sm_k)} = K_0(X)$ on smooth schemes of $K$ is algebraic $K$-theory. To explain the map $U_{r,d}(m) \to \Omega_{\P^1} K$, notice that the value of a pointed preaheaf $P$ on a smooth $k$-scheme $X$ is the same as the space of pointed maps $[X_+, P]_{\PSh(\Sm_k)_*}$. We need to show that a bundle on $\P^1 \times X$ induces an element of $\Omega_{\P^1} K$ valued at $X$ which is exactly an element of the space $[\P^1 \wedge X_+, K]_*$. A vector bundle on $\P^1 \times X$ zero along $\sigma_{\infty}$ defines a $K$-theory class on the presheaf $\P^1 \wedge X_+$, because $[\P^1 \wedge X_+, K]_*$ is the fiber of
\[ K_0(\P^1 \times X) \xrightarrow{\sigma_{\infty}^*} K_0(X) \]
since $\P^1 \wedge X_+$ is the cofiber of $\sigma_{\infty} : X \to \P^1 \times X$. Finally, $\alpha$ is the same map as (the vertical maps) in \eqref{gs1}  which is defined as 
\[ \E \mapsto \pi_* \E(m-1). \]
The idea behind these maps is that the map to $K$-theory should be $\E \mapsto [\pi_! \E(-1)]$.  (Aside:  $\E \mapsto [\pi_! \E]$ is ``equally good'', but our choice makes the quiver description, and some of the extensions in \cite{bryanvakil}, cleaner.)  Because $\E$ is trivialized along $\sigma_\infty$ we have $[\pi_! \E(-1)] = [\pi_! \E(m-1)] - mr$. Moreover $\E$ lives in $U_{r,d}(m)$ so $\E(m-1)$ has no higher direct images. Together this gives an equality of $K$-classes
\[ [\pi_! \E(-1)] = [\pi_* \E(m-1)] - mr. \]
To check commutativity of this diagram, we need to relate $[\E]$ to $[\pi_* \E(m-1)]$. This is done by Str\o{}mme's Lemma~\ref{l:stromme}
implying that 
\[ 0 \to (\pi^* \pi_* \E(m-1))(-1) \to \pi^* \pi_* \E(m) \to \E(m) \to 0 \]
is short exact. Therefore, in $K$-theory
\[ [\E(m)] = \pi^* [\pi_* \E(m)] - \pi^* [\pi_* \E(m-1)] \gamma. \]
Furthermore, the $m$-regularity of $\E$ implies the exactness of 
\[ 0 \to \pi_* \E(m-1) \to \pi_* \E(m) \to \sigma^*_{\infty} \E(m) \to 0. \]
Using the trivialization at infinity, 
\[ [\pi_* \E(m)] = [\pi_* \E(m-1)] + r \]
and therefore
\[ [\E(m)] = (r + \pi^* [\pi_* \E(m-1)] - \pi^* [\pi_* \E(m-1)] \gamma) = r + (1 - \gamma) \pi^* [\pi_* \E(m-1)].  \]
The trivialization proves that $[\E(m)] = (1 - \gamma) mr + [\E]$ and therefore
\[ [\E] - r = (1- \gamma) (\pi^* [\pi_* \E(m-1)] - mr) \]
proving that the diagram \eqref{eq:benaug30} commutes. Notice that this diagram commutes with the stabilization maps $T'_{r,d}(m) \to T'_{r+1,d}(m)$ given by $\E \mapsto \E \oplus \cO$ and $\BGL(d + rm) \to \BGL(d + r(m+1))$ via $\mathcal{F} \mapsto \mathcal{F} \oplus \cO^{\oplus m}$. Here we used that $\pi_* \cO(m-1) \cong \cO^{\oplus m}$. Therefore, it defines a diagram of sequences,
\begin{center}
\begin{tikzcd}
 & T'_{\bullet,d}(m) \arrow[rd] \arrow[ld, "\alpha"']
 \\
 \BGL(d + m \bullet) \arrow[d, "{[-] - m\bullet}"] & & U_{\bullet,d}(m) \arrow[d, "{[-] - \bullet}"]
 \\
 K \arrow[rr, "-\beta"] & & \Omega_{\P^1} K 
\end{tikzcd}
\end{center}
Hence applying $L_{\mot}$ and passing to the colimit we obtain,
\begin{center}
\begin{tikzcd}
 & M_{\cGsp^{\A^1}}^t T'_{\bullet,d}(m) \arrow[rd] \arrow[ld, "\alpha"']
 \\
 M_{\cGsp^{\A^1}}^b \BGL(d + m \bullet) \arrow[d, "{[-] - m\bullet}"] \arrow[rr, dashed, "\beta_{LV}"] & & M_{\cGsp^{\A^1}}^b U_{\bullet,d}(m) \arrow[d, "{[-] - \bullet}"]
 \\
 K \arrow[rr, "-\beta"] & & \Omega_{\P^1} K 
\end{tikzcd}
\end{center}
where we use that $T'_{\bullet, d}(m) \to \BGL(d + m\bullet)$ is a spatial realization of $\BGL$ and $T'_{\bullet,d}(m) \to U_{\bullet,d}(m)$ is a spatial realization of $\Omega_{\alg,d}^2 \BGL$ by Definition~\ref{d:final}. By Theorem~\ref{thm:A1_top_bottom_comparison} these maps are equivalences; the Bott map $\beta_{LV}$ in $\Ho \, \cGst$ is, by definition, the dashed map induced by this roof of equivalences. Therefore, passing to the disjoint union over all $d$ we obtain that the desired diagram commutes. By Theorem~\ref{thm:bpg}, $\beta_{LV}$ is an equivalence. Furthermore, we saw the map $\Z \times \BGL \to K$ is a motivic equivalence in the previous section. Hence, taking the union over each degree $d$, we obtain that the top and left maps are equivalences. Finally, the Bott map $\beta$ is a motivic equivalence (for example see \cite[Theorem 3.16]{bachmann}).
This implies that the remaining map
\[ \Omega^2_{\alg} \BGL \to \Omega_{\P^1} K \]
is also a motivic equivalence. 
\end{proof}

As a consequence, we have shown that algbraic loops $\Omega^2_{\alg}$, in the sense of the stack of bundles over $\P^1$, correctly converges to $\Omega_{\P^1}$ for these classifying spaces in the motivic category. This should be viewed as an analogue of the theorem of Cohen--Lupercio--Segal \cite[Theorem 1]{cls} (see \S \ref{pf:trad}) in the context of motivic homotopy theory. The author is unaware of a proof of this equivalence independent of Bott periodicity for algebraic $K$-theory analogous to the proof given in \cite{cls} which does not depend on Bott periodicity for topological $K$-theory.

}  

\end{document}

%% file: RVpreamble.tex
\documentclass[12pt,letterpaper]{amsart} 
\usepackage{palatino, euler, epic,eepic, amssymb, xypic, floatflt, microtype}
\usepackage{verbatim,color}
\usepackage[linktocpage,hidelinks]{hyperref}
\hypersetup{
    colorlinks,
    linkcolor={blue!50!black},
    citecolor={green!50!black},
    urlcolor={red!80!black}
}

\usepackage{tikz-cd}

\input xy
\xyoption{all}

 \newlength{\baseunit}               
 \newcount{\numlines}                
\newcommand{\bpf}{\noindent {\em Proof.  }}
\newcommand{\epf}{\qed \vspace{+10pt}}

\newcommand{\point}{\vspace{1mm}\par \noindent \refstepcounter{subsection}{\thesubsection.} }
\newcommand{\tpoint}[1]{\vspace{1mm}\par \noindent \refstepcounter{subsection}{\thesubsection.} 
  {\bf #1. ---} }
\newcommand{\epoint}[1]{\vspace{1mm}\par \noindent \refstepcounter{subsection}{\thesubsection.} 
  {\em #1.} }
\newcommand{\bpoint}[1]{\vspace{1mm}\par \noindent \refstepcounter{subsection}{\thesubsection.} 
  {\bf #1.} }

\newcommand{\tU}{\widetilde{U}}

\newcommand{\proj}{\mathbb P}

\newcommand{\C}{\mathbb{C}}

\newcommand{\G}{\mathbb{G}}

\newcommand{\Z}{\mathbb{Z}}
\newcommand{\A}{\mathbb{A}}
\newcommand{\E}{\mathbb{E}}
\renewcommand{\L}{\mathbb{L}}
\newcommand{\Q}{\mathbb{Q}}

\def\CC{\mathbb{C}}

\def\RR{\mathbb{R}}

\newcommand{\cA}{{\mathscr{A}}}
\newcommand{\cB}{{\mathscr{B}}}
\newcommand{\cC}{{\mathscr{C}}}
\newcommand{\cE}{{\mathscr{E}}}
\newcommand{\cF}{{\mathscr{F}}}
\newcommand{\cG}{{\mathscr{G}}}

\newcommand{\cV}{{\mathscr{V}}}

\newcommand{\cX}{{\mathscr{X}}}
\newcommand{\cY}{{\mathscr{Y}}}

\newcommand{\oh}{{\mathscr{O}}}

\newcommand{\al}{\alpha}

\newcommand{\be}{\beta}

\newcommand{\si}{\sigma}

\renewcommand{\top}{\rm{top}}
\newcommand{\red}{\rm{red}}

\newcommand{\coloneq}{:=}

\newcommand{\propernormal}{%
  \mathrel{\ooalign{$\lneq$\cr\raise.22ex\hbox{$\lhd$}\cr}}}

\newcommand{\colim}{\operatorname{colim}}

\renewcommand{\span}{\operatorname{span}}
\newcommand{\coker}{\operatorname{coker}}

\newcommand{\Hom}{\operatorname{Hom}}

\newcommand{\Maps}{\operatorname{Maps}}

\newcommand{\codim}{\operatorname{codim}}

\newcommand{\Spec}{\operatorname{Spec}}

\newcommand{\lremind}[1]{{\bf[label:  #1]}}
\newcommand{\notation}[1]{}
\renewcommand{\lremind}[1]{{}}

\newcommand{\cut}[1]{}



%% file: BCpreamble.tex
\makeatletter
\newcommand{\chapterauthor}[1]{%
  {
    \parindent0pt\vspace*{-10pt}%
    \linespread{1}\scshape by #1%
  }
}
\makeatother

\usepackage{relsize}
\usepackage[bbgreekl]{mathbbol}
\usepackage{amsfonts}
\DeclareSymbolFontAlphabet{\mathbb}{AMSb} 
\DeclareSymbolFontAlphabet{\mathbbl}{bbold}
\newcommand{\Prism}{{\mathlarger{\mathbbl{\Delta}}}} 

\DeclareMathOperator{\fib}{\mathrm{fib}}
\DeclareMathOperator{\cofib}{\mathrm{cofib}}

\theoremstyle{definition}

\usepackage{xspace}
\newcommand{\etale}{\'{e}tale\xspace}

\usetikzlibrary{arrows.meta}

\newcommand*{\TikzBiArrow}[1][]{\mathbin{\tikz [baseline=-0.25ex,-latex, #1] \draw [#1] (0pt,0.5ex) -- (1.3em,0.5ex) node[midway, above, label={[label distance=-0.4cm]$\sim$}] {} ;}}

\newsavebox{\pullbacksymbol}
\sbox\pullbacksymbol{%
\begin{tikzpicture}%
\draw (0,0) -- (1ex,0ex);%
\draw (1ex,0ex) -- (1ex,1ex);%
\end{tikzpicture}}

\newsavebox{\pushoutsymbol}
\sbox\pushoutsymbol{%
\begin{tikzpicture}%
\draw (0,0) -- (-1ex,0ex);%
\draw (-1ex,0ex) -- (-1ex,-1ex);%
\end{tikzpicture}}

\newcommand{\pullback}{\arrow[dr, phantom, "\usebox\pullbacksymbol" , very near start, color=black]}

\newcommand{\pushout}{\arrow[lu, phantom, "\usebox\pushoutsymbol" , very near start, color=black]}